\documentclass[11pt]{amsart}
\usepackage[margin=1in]{geometry}
\usepackage[T1]{fontenc}
\IfFileExists{newtxtext.sty}{\usepackage{newtxtext}\usepackage[varvw]{newtxmath}}{\usepackage{mathptmx}\usepackage{amssymb}}
\usepackage{amsmath,amsthm,mathtools,enumitem}
\usepackage[colorlinks=true,linkcolor=blue,citecolor=blue,urlcolor=blue]{hyperref}
\usepackage{color}
\usepackage{tikz}

\providecommand{\doi}[1]{\href{https://doi.org/#1}{\nolinkurl{doi:#1}}}
\newcommand{\Z}{\mathbb Z}\newcommand{\R}{\mathbb R}\newcommand{\T}{\mathbb T}
\newcommand{\N}{\mathbb N}\newcommand{\F}{\mathbb F}
\newcommand{\eps}{\varepsilon}\newcommand{\Hilb}{\operatorname{Hilb}}
\newcommand{\cone}{\operatorname{cone}}\newcommand{\conv}{\operatorname{conv}}
\newcommand{\supp}{\operatorname{supp}}\newcommand{\wt}{\operatorname{wt}}
\newcommand{\ord}{\operatorname{ord}}\newcommand{\Vol}{\operatorname{Vol}_{\Z}}
\DeclareMathOperator{\spanR}{span_{\R}}
\makeatletter
\def\@setaddresses{\par
  \nobreak \begingroup
\footnotesize
  \def\author##1{\nobreak\addvspace\bigskipamount}%
  \def\\{\unskip, \ignorespaces}%
  \interlinepenalty\@M
  \def\address##1##2{\begingroup
    \par\addvspace\bigskipamount\indent
    \@ifnotempty{##1}{\ignorespaces##1\unskip: }%
    {\scshape\ignorespaces##2}\par\endgroup}%
  \def\curraddr##1##2{\begingroup
    \@ifnotempty{##2}{\nobreak\indent\curraddrname
      \@ifnotempty{##1}{, \ignorespaces##1\unskip}\/:\space
      ##2\par}\endgroup}%
  \def\email##1##2{\begingroup
    \@ifnotempty{##2}{\nobreak\indent\emailaddrname
      \@ifnotempty{##1}{, \ignorespaces##1\unskip}\/:\space
      \ttfamily##2\par}\endgroup}%
  \def\urladdr##1##2{\begingroup
    \def~{\char`\~}%
    \@ifnotempty{##2}{\nobreak\indent\urladdrname
      \@ifnotempty{##1}{, \ignorespaces##1\unskip}\/:\space
      \ttfamily##2\par}\endgroup}%
  \addresses
  \endgroup
}
\makeatother
\renewcommand{\urladdrname}{\textit{ORCID}}
\theoremstyle{plain}
\newtheorem{theorem}{Theorem}[section]
\newtheorem{proposition}[theorem]{Proposition}
\newtheorem{lemma}[theorem]{Lemma}
\newtheorem{corollary}[theorem]{Corollary}
\theoremstyle{definition}

\newtheorem{examplex}[theorem]{Example}
\newenvironment{example}
  {\pushQED{\qed}\examplex}
  {\popQED\endexamplex}
\theoremstyle{remark}
\newtheorem{remark}[theorem]{Remark}
\numberwithin{equation}{section}
\numberwithin{figure}{section}
\setlist[enumerate]{itemsep=1pt,topsep=2pt}
\setlist[itemize]{itemsep=1pt,topsep=2pt}
\allowdisplaybreaks[2]
\begin{document}
\title[Arithmetic universality in disjunctions of staircase simplices]{Arithmetic Universality and Integer Decomposition\texorpdfstring{\\}{ } in Disjunctions of Three Staircase Simplices}
\author{Yushan Qu}
\address[Yushan Qu]{University of Michigan, Ann Arbor, MI, USA}
\email{yushanqu@umich.edu}
\urladdr{\href{https://orcid.org/0000-0002-2862-4224}{0000-0002-2862-4224}}
\author{Jon Lee}
\address[Jon Lee]{University of Michigan, Ann Arbor, MI, USA}
\email{jonxlee@umich.edu}
\urladdr{\href{https://orcid.org/0000-0002-8190-1091}{0000-0002-8190-1091}}
\date{\today}
\subjclass[2020]{52B20, 90C10, 90C11, 05E40, 20K01}
\keywords{Disjunctive programming, integer decomposition property, Hilbert basis, staircase simplex, Cayley polytope, finite abelian group, $0/1$ polytope}

\begin{abstract}
The staircase simplices $\{x\in[0,1]^d:x_{\sigma(1)}\le\cdots\le x_{\sigma(d)}\}$, $\sigma\in S_d\,$, are the $d!$ unimodular cells of the staircase triangulation of the unit cube. The disjunctive hull (Cayley polytope) of a selection $\mathcal P_0\,,\ldots,\mathcal P_n$ of them has the integer decomposition property (IDP) if, for all integers $\lambda_i\ge0$, every lattice point of $\sum_i\lambda_i\mathcal P_i$ is a sum of vertices, $\lambda_i$ of them from $\mathcal P_i\,$. For two simplices, IDP always holds. For three, every nontrivial finite abelian group $G$ arises, in the optimal dimension $\Theta(\log|G|)$, as the residue group of a simplex facet whose $|G|-1$ nonzero box points are Hilbert basis elements of the whole cone; $k$ times such a point is a sum of vertices exactly when its order in $G$ divides $k$. Isolating this arithmetic is costly: if a simplex face has a single Hilbert basis element $h$ that is not a vertex, then $kh$ is a sum of vertices exactly when some $q$ divides $k$, and $d\ge3q-2$, attained for every $q\ge3$; we also find the least dimension for each $G$ with the fewest such elements. Without isolation, every nonzero submonoid of $\mathbb{Z}_{\ge0}$ omitting $1$ is $\{k\ge0:kh$ is a sum of vertices$\}$ for a Hilbert basis element $h$ of degree three. For every $0/1$ polytope with $N$ vertices in $\mathbb{R}^r$, three staircase simplices in dimension $O(N(r+1))$ have a face that has IDP exactly when the polytope does. Deciding whether a lattice point $h$ is a sum of vertices is strongly NP-complete, even given vertex expressions of $2h$ and $3h$, and recognizing IDP of a designated facet is coNP-complete; for the whole hull, IDP can first fail in degree linear in $d$, and it holds whenever $d\le5$. For a family of triples defined by words over three letters, we determine the entire Hilbert basis and characterize IDP by two forbidden subsequences.
\end{abstract}

\maketitle

\section{Introduction}\label{sec-intro}

\par\medskip\noindent\emph{Disjunctions and their arithmetic.}
Disjunctive programming models a choice among polytopes $\mathcal P_0\,,\ldots,\mathcal P_n\subseteq\R^d$ by the convex hull, or \emph{disjunctive hull},
\[
 \mathcal D:=\conv\textstyle\bigcup_{i}\bigl(\mathcal P_i\times\{\mathbf e_i\}\bigr)
\]
of the polytopes lifted by indicator vectors ($\mathbf e_0:=\mathbf 0$); the extended formulations of Balas~\cite{Balas} and Jeroslow and Lowe~\cite{JL} describe exactly this hull, and its properties govern the strength of mixed-integer models of unions~\cite{Vielma}. In polyhedral geometry, $\mathcal D$ is the \emph{Cayley polytope} of $\mathcal P_0\,,\ldots,\mathcal P_n\,$. For the simplices studied here, the only lattice points of each $\mathcal P_i$ are its vertices, so the lattice points of $\mathcal D$ are its lifted vertices, and the next question is arithmetic. The cone $K$ over $\mathcal D$ consists of the points $(x,\lambda)$ with multiplicities $\lambda\ge0$ and $x\in\sum_i\lambda_i\mathcal P_i\,$; the \emph{degree} of such a point is $\lambda_0+\cdots+\lambda_n\,$. We say that a lattice point of $K$ \emph{splits} if it is a sum of lifted vertices, that is, if it lies in the \emph{vertex monoid} they generate. The disjunction has the \emph{integer decomposition property} (IDP) if every lattice point of $K$ splits; equivalently, for all integers $\lambda_i\ge0$, every lattice point of $\lambda_0\mathcal P_0+\cdots+\lambda_n\mathcal P_n$ is a sum of vertices, $\lambda_i$ of them from $\mathcal P_i\,$. This holds exactly when the lifted vertices form the \emph{Hilbert basis} $\Hilb(K)$ of the cone, the unique minimal set of lattice points generating all of them as a monoid. Hilbert bases underlie total dual integrality~\cite{GP} and integer Carath\'eodory questions~\cite{CFS,BGHMW}. We call a lattice point of the cone \emph{indecomposable} if it is a Hilbert basis element. When IDP fails, the obstructions are the \emph{indecomposable holes}: elements of $\Hilb(K)$ other than lifted vertices. Following Takemura and Yoshida~\cite{TY} and Hemmecke, Takemura and Yoshida~\cite{HTY}, the lattice points of the cone outside the vertex monoid are called \emph{holes}; an indecomposable hole is in particular a \emph{fundamental} hole in their sense (one from which no other hole can be reached by subtracting a sum of lifted vertices), but not conversely, because a fundamental hole can be a sum of two holes.

\par\medskip\noindent\emph{Clearing.} Even when a lattice point $h$ of the cone does not split, some positive-integer multiple $kh$ does. We measure this with three notions:
\begin{itemize}
\item the \emph{clearing set} $\mathcal C_K(h):=\{k\in\Z_{>0}:kh\text{ lies in the vertex monoid}\}$;
\item the \emph{least clearing integer} $\delta_K(h):=\min\mathcal C_K(h)$;
\item the \emph{clearing order} $q$, if $\mathcal C_K(h)$ consists exactly of the multiples of $q$; then $\delta_K(h)=q$.
\end{itemize}
For instance, a point with clearing order $2$ does not split, but twice it does. This clearing terminology is ours; it refines normality, which says that every clearing set contains $1$, and seminormality. Not every clearing set consists of the multiples of one number: Corollary~\ref{cor-clearsets-all} exhibits an indecomposable hole whose clearing set is $\{4,5,8,9,10,12,13,\ldots\}$, the sums of $4$'s and $5$'s; it has least clearing integer $4$ but no clearing order.

\par\medskip\noindent\emph{Staircase simplices.}
We take the polytopes to be \emph{staircase simplices} $\mathcal P_\sigma=\{x\in[0,1]^d:x_{\sigma(1)}\le\cdots\le x_{\sigma(d)}\}$, $\sigma\in S_d\,$, the $d!$ cells of the staircase (Freudenthal--Kuhn) triangulation of the cube. They are the cells of standard mixed-integer models of nonseparable piecewise-linear functions of Vielma, Ahmed and Nemhauser~\cite{VAN}; each is unimodular. Qu and Lee studied their disjunctive hulls through big-$M$ lifting~\cite{QL}, the facets of the hull of all $d!$ of them, which is the Cayley polytope of Stanley's triangulation~\cite[Proposition~3.1, Theorem~3.6]{Cayley}, and the facet normals of Minkowski sums of staircase simplices~\cite{Stair}. We use their three-class determinant bound \cite[Lemma~4.6]{Stair} and zero-sum description of facet normals \cite[Theorem~3.1]{Stair}. For two staircase simplices, IDP always holds, because their lifted vertex matrix is totally unimodular (Proposition~\ref{prop-two}). This paper asks how complex the arithmetic of a disjunction of \emph{three} staircase simplices can be. The answer is that every finite abelian group occurs, and that the complexity can be located, bounded, and in one family classified.

\par\medskip\noindent\emph{An example.}
The $d+1$ vertices of $\mathcal P_\sigma$ are $0/1$ vectors forming a chain: starting from $0$, each vertex switches on one more coordinate, until all coordinates are $1$. The order in which the coordinates are switched on, namely $\sigma(d),\sigma(d-1),\ldots,\sigma(1)$, is the \emph{growth order} of $\mathcal P_\sigma\,$; the vertices are the indicator vectors of its \emph{prefixes} (initial segments, including the empty one). For example, the growth order $(1,2,4,5,6,3)$ describes the simplex with vertices $0$, $\mathbf e_1\,$, $\mathbf e_1+\mathbf e_2\,$, $\mathbf e_1+\mathbf e_2+\mathbf e_4\,$, \dots, the all-ones vector $\mathbf e$, that is, $\{x\in[0,1]^6:x_3\le x_6\le x_5\le x_4\le x_2\le x_1\}$. In $\R^6$, take the growth orders
\[
 \mathcal P_{\sigma_0}:(1,2,4,5,6,3),\qquad\mathcal P_{\sigma_1}:(1,3,4,6,5,2),\qquad\mathcal P_{\sigma_2}:(2,3,5,6,4,1).
\]
Let $a\in\mathcal P_{\sigma_0}\,$, $b\in\mathcal P_{\sigma_1}\,$, $c\in\mathcal P_{\sigma_2}$ be the vertices with supports $\{1,2,4,5\}$, $\{1,3,4\}$, $\{2,3,5\}$, respectively. The point $x=(1,1,1,1,1,0)$ with multiplicities $\lambda=(1,1,1)$ is half the sum of the three origins and $a$, $b$, $c$, so $h=(x,\lambda)$ is a lattice point of the cone. A splitting of $h$ would write $x=u+v+w$ with vertices $u\in\mathcal P_{\sigma_0}\,$, $v\in\mathcal P_{\sigma_1}\,$, $w\in\mathcal P_{\sigma_2}\,$, none containing coordinate $6$, whose supports partition $\{1,\ldots,5\}$. Because $v_5=1$ would force $v_6=1$, we have $u_5=1$ or $w_5=1$. If $u_5=1$, then $u=a$, and neither $v$ nor $w$ can supply $3$ without repeating a coordinate of $a$; if $w_5=1$, then $w=c$, and neither $u$ nor $v$ can supply $4$ without repeating a coordinate of $c$. Thus $h$ does not split, while $2h$ splits into the three origins and $a$, $b$, $c$; Theorem~\ref{thm-words} shows that $h$ is the only indecomposable hole, and by Proposition~\ref{prop-dim5} three staircase simplices with $d\le5$ always have IDP, so six is minimal.

\begin{figure}[t]
\centering
\begin{tikzpicture}[x=0.62cm,y=0.62cm,font=\small,scale=1.349]
 \foreach \p/\l in {1/A,2/B,3/C,4/A,5/B} {\node at (\p+1.6,4.25) {$\l$};}
 \node at (7.6,4.25) {\scriptsize guard};
 \foreach \p in {1,...,6} {\node[gray] at (\p+1.6,3.65) {\scriptsize$\p$};}
 \foreach \row/\lab/\ord/\k in {2.6/{\mathcal P_{\sigma_0}\ (A,B)}/{1,2,4,5,6,3}/4, 1.6/{\mathcal P_{\sigma_1}\ (A,C)}/{1,3,4,6,5,2}/3, 0.6/{\mathcal P_{\sigma_2}\ (B,C)}/{2,3,5,6,4,1}/3} {
   \node[anchor=east] at (2.0,\row) {$\lab$};
   \foreach \c [count=\i] in \ord {
     \pgfmathsetmacro{\fillc}{\i<=\k ? 1 : 0}
     \ifdim\fillc pt>0.5pt \fill[blue!18] (\i+1.15,\row-0.4) rectangle (\i+2.05,\row+0.4);\fi
     \draw (\i+1.15,\row-0.4) rectangle (\i+2.05,\row+0.4);
     \node at (\i+1.6,\row) {$\c$};
   }
 }
 \node[anchor=west] at (8.3,2.6) {$a$}; \node[anchor=west] at (8.3,1.6) {$b$}; \node[anchor=west] at (8.3,0.6) {$c$};
 \begin{scope}[xshift=6.9cm,yshift=0.2cm]
  \coordinate (A0) at (0,0); \coordinate (B0) at (1,2); \coordinate (C0) at (2,0);
  \coordinate (A1) at (3,2); \coordinate (B1) at (4,0); \coordinate (C1) at (5,2);
  \coordinate (A2) at (6,0); \coordinate (B2) at (7,2);
  \fill[blue!18] (A0)--(B0)--(C0)--cycle; \fill[blue!18] (C1)--(A2)--(B2)--cycle;
  \draw (A0)--(B0)--(C0)--(A1)--(B1)--(C1)--(A2)--(B2);
  \draw (A0)--(C0)--(B1)--(A2); \draw (B0)--(A1)--(C1)--(B2);
  \foreach \v/\l/\pos in {A0/A_0/below,B0/B_0/above,C0/C_0/below,A1/A_1/above,B1/B_1/below,C1/C_1/above,A2/A_2/below,B2/B_2/above} {\fill (\v) circle (1.6pt); \node[\pos] at (\v) {$\l$};}
  \node at (1,0.7) {$T_0$}; \node at (6,1.3) {$T_5$};
 \end{scope}
\end{tikzpicture}
\caption{The smallest example. \emph{Left:} the growth orders, read off from the word $ABCAB$ (coordinate $p$ carries the $p$th letter; $6$ is the guard); the shaded prefixes are the vertices $a,b,c$; $h$ is half the sum of the three shaded lifted vertices and the three lifted origins. \emph{Right:} the triangle strip of Theorem~\ref{thm-words-graph}. The triangles $T_0$ (before the word) and $T_5$ (after it) are disjoint and joined by no edge; this pair is the indecomposable hole $h$.}
\label{fig-example}
\end{figure}
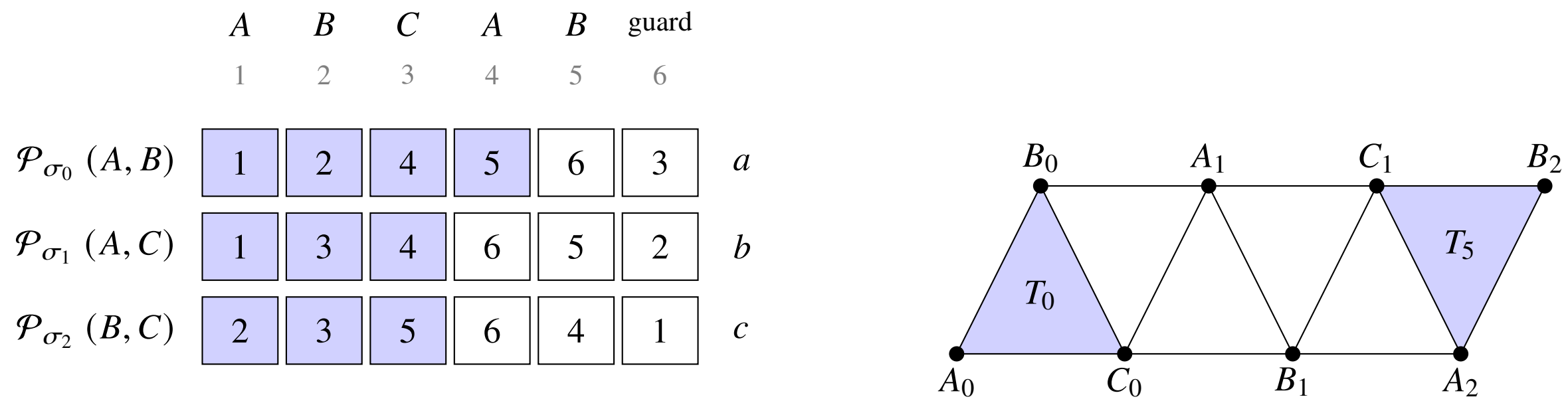

These orders come from the word $ABCAB$ (Section~\ref{sec-words}): each simplex switches on the coordinates carrying its two letters in order, then a separator, then the rest in reverse; the indecomposable hole reflects the $C$ between both $A$'s and both $B$'s (Figure~\ref{fig-example}).

\par\medskip\noindent\emph{Main results.}
For a face $F$ of $\mathcal D$ that is a simplex, write $G_F$ for its \emph{residue group}: the lattice points in the linear span of the cone over $F$ modulo the lattice generated by its lifted vertices. Each residue $g$ has a unique representative $H_g$ in the fundamental parallelepiped of $F$; these are the \emph{box points} of $F$, with $H_0=0$. We call $K$ and $\Hilb(K)$ the \emph{ambient} cone and Hilbert basis, as opposed to those of a face.

\smallskip\noindent\hypertarget{res-A}{\emph{(A) Arithmetic universality}} (Theorem~\ref{thm-main}). Every nontrivial finite abelian group $G$ is $G_F$ for a simplex facet $F$ of three staircase simplices in dimension $d=O(\log|G|)$. The $|G|-1$ nonzero $H_g$ lie in the ambient Hilbert basis, so they are indecomposable holes; they have degree between $ad$ and $d$ for an absolute constant $a>0$, and $H_g$ has clearing order $\ord(g)$. Determinant and group-rank bounds show that $\Theta(\log|G|)$ is optimal. One triple in dimension $O(m)$ serves all groups of order at most $2^m$ at once (Theorem~\ref{thm-simul}), so clearing every lattice point with one multiplier can need $2^{\Theta(d)}$ bits, although every single lattice point clears at a multiplier at most $3^{(d-1)/3}$ (Corollary~\ref{cor-uniform}). The same construction gives inclusion-minimal non-IDP simplex facets with cyclic residue group of every order $q$ in dimension $O(\log q)$, all of whose box points are indecomposable of the same degree (Theorem~\ref{thm-minfacet}), and $2^{\Theta(d)}$ inclusion-minimal non-IDP faces (Corollary~\ref{cor-manymin}).

\smallskip\noindent\hypertarget{res-B}{\emph{(B) Clearing sets and faces with few indecomposable holes}} (Section~\ref{sec-cyclic}). A simplex face with $r$ indecomposable holes has a residue group of order at most $\binom{\lfloor(d-1)/3\rfloor+r}{r}$ (Theorem~\ref{thm-atoms}). If the cone of \emph{any} face has a single indecomposable hole $h$, with least clearing integer $q$, then $d\ge3q-2$, or $d\ge3q-3$ for $\mathcal D$ itself (Theorem~\ref{thm-onehole-any}), and this is attained for every $q\ge3$ by simplex facets on which $h$ has clearing order $q$ (Theorem~\ref{thm-onehole}). With only $\mu(G)$ indecomposable holes, $\mu(G)$ the minimal number of generators of a nontrivial residue group $G$, the least possible $d$ is exactly $3\Sigma(G)-2$, for $\Sigma(G)$ the least sum of $\mu(G)$ cyclic orders composing $G$, except that it is $6$ for $\Z/2\Z$ (Theorem~\ref{thm-khole}). For $d\ge7$ the largest residue group of a simplex face with at most $r\ge1$ indecomposable holes has order $e^{\Theta(r\log(1+d/r))}$ (Corollary~\ref{cor-tradeoff}), and an indecomposable hole of clearing order $q$ on a simplex face has degree at most $\lfloor((q-1)(d+1)+1)/q\rfloor$, which is sharp for $q\ge3$, $d\ge3q-1$ and for $q=2$, $d\ge7$ (Theorem~\ref{thm-degq}). When other indecomposable holes may share the face, this rigidity disappears: clearing order $q$ occurs at degree three on a simplex face in dimension $O(\log q)$ (Corollary~\ref{cor-sort}), and on faces that need not be simplices every nonzero additive submonoid of $\Z_{\ge0}$ omitting $1$ is $\{0\}\cup\mathcal C_K(h)$ for an indecomposable hole $h$ of the minimal degree three (Theorem~\ref{thm-clearsets} and Corollary~\ref{cor-clearsets-all}).

\smallskip\noindent\emph{The two/three threshold} (Theorem~\ref{thm-threshold}). Two staircase simplices always give an IDP hull; for each fixed $n\ge2$, the maximum number of indecomposable holes of $n+1$ staircase simplices in $\R^d$ and the largest least clearing integer of a lattice point are both $2^{\Theta(d)}$. Explicitly, the binary family of Theorem~\ref{thm-binary} has $2^k-1$ indecomposable holes in dimension $4k+2$.

\smallskip\noindent\hypertarget{res-C}{\emph{(C) Universality for IDP}} (Theorem~\ref{thm-idp-transfer}). For every $0/1$ polytope $\mathcal Q$ with $N$ vertices in $\R^r$, three staircase simplices in dimension $O(N(r+1))$ have a face $F_{\mathcal Q}$ such that every lattice point of the cone over $F_{\mathcal Q}$ is a sum of lifted vertices on $F_{\mathcal Q}$ if and only if $\mathcal Q$ has IDP. The indecomposable holes of the cone over $\mathcal Q$ correspond bijectively to those of the cone over $F_{\mathcal Q}\,$, by an explicit linear map sending each lifted vertex $(v,1)$ of $\mathcal Q$ to a sum of lifted vertices on $F_{\mathcal Q}$; corresponding indecomposable holes have the same clearing sets, degrees are multiplied by $6(r+1)$, and the indecomposable holes on $F_{\mathcal Q}$ lie in the ambient Hilbert basis. For full-dimensional $\mathcal Q$, the face can be a facet.

\smallskip\noindent\hypertarget{res-D}{\emph{(D) Complete Hilbert bases for word disjunctions}} (Theorem~\ref{thm-words-summary}). For triples whose growth orders are read off from a \emph{word} (a sequence) over $A,B,C$, we determine the \emph{entire ambient} Hilbert basis. All indecomposable holes have degree three and clearing order $2$; IDP amounts to avoiding the subsequences $XYZXY$ and $XYZYX$ for distinct letters $X,Y,Z$; the minimal non-IDP faces are five-dimensional simplices of normalized volume two; and the maximum number of indecomposable holes of a word disjunction in dimension $d\ge6$ is $\binom{d-4}2$. Every lattice point of the cone has clearing set $2\Z_{>0}$ or $\langle2,\tau\rangle\setminus\{0\}$ with $\tau$ odd (Theorem~\ref{thm-words-clearing}); an odd clearing threshold $\tau=q\ge3$ first occurs in dimension $5q-4$, so seminormality (every lattice point $h$ with $2h$ and $3h$ in the vertex monoid lies in it) first fails in dimension $11$ within the word family (Corollary~\ref{cor-words-clearing}). Given the word, IDP is decided in linear time.

\smallskip\noindent\hypertarget{res-E}{\emph{(E) Splitting and recognition}} (Theorem~\ref{thm-complexity} and Corollary~\ref{cor-recognition}). On initial faces where each coordinate can be switched on by at most two simplices, deciding whether a given lattice point $h$ splits is a linear program plus $2$-SAT; in general it is strongly NP-complete, even when $2h$ and $3h$ are given as vertex sums, fixed-parameter tractable in the number of coordinates available to all three simplices on the face, and W[1]-hard in the least multiplicity. Through~\hyperlink{res-C}{(C)}, recognizing IDP of a designated facet, one given by an exposing form, is coNP-complete, and although IDP of a simplex facet is a determinant test, deciding whether a designated simplex facet has an indecomposable hole of degree at most $T$ is strongly NP-complete, even when every nonzero box point clears at two (Theorem~\ref{thm-holedeg-hard}). This hardness concerns a designated facet: a facet can be unimodular while the hull has many indecomposable holes (Corollary~\ref{cor-words-unimodular}). For the whole hull, certificates of degree below the rank place the recognition of IDP in coNP (Lemma~\ref{lem-cert}).

\smallskip\noindent\hypertarget{res-F}{\emph{(F) Where integer decomposition fails}} (Section~\ref{sec-failures}). IDP can first fail in large degree: the least degree of an indecomposable hole of the whole hull can be any $k\ge3$ in dimension at most $8k+12$, even when twice every lattice point splits (Theorem~\ref{thm-firsthole}), while doubling forces $k\le\lfloor(d+3)/2\rfloor$ (Proposition~\ref{prop-double}). So testing lattice points up to degree $o(d)$ cannot certify IDP, and the certificates of Lemma~\ref{lem-cert} are tight up to a constant factor. A multiple $kh$ lying in the lattice generated by the rays of the minimal face of $h$ can fail to split for $k$ as large as $2^{\Theta(d)}$, and this order is sharp (Theorem~\ref{thm-delay}). Six dimensions are needed for any failure (Proposition~\ref{prop-dim5}), and away from the designated faces of our constructions, an indecomposable hole can lie in no simplex face at all (Proposition~\ref{prop-offface-y}).

\par\medskip\noindent\emph{Technique.}
Results~\hyperlink{res-A}{(A)}--\hyperlink{res-C}{(C)} rest on one recipe in three steps. \emph{(i) From a matrix to a face.} Let $A$ be a $0/1$ matrix with $s$ columns, split into three classes with pairwise disjoint supports inside each class, and let $\alpha$ be a weight vector with $A^T\alpha=0$ and suitable signs. Each class defines a growth order that switches on its supports one after another, and a single linear form built from $\alpha$ exposes exactly the vertices at the ends of the supports (Lemmas~\ref{lem-growth} and~\ref{lem-shared}). \emph{(ii) From the face to a group.} Within each simplex, consecutive support endpoints differ by a column of $A$. For a simplex face, the residue group is therefore the \emph{torus kernel} $\Gamma_A:=\{f\in\T^s:Af\equiv0\}$, $\T:=\R/\Z$, with one scalar per column, and the box points have explicit coefficients; for the more general faces in~\hyperlink{res-C}{(C)}, the torus kernel of the ray matrix identifies the induced lattice and vertex monoid through the \emph{series lift}, which replaces each vertex of $\mathcal Q$ by a chain of generators summing to it. \emph{(iii) Programming the group.} Each row of $A$ meets each class at most once, so it is an equation among at most three scalars; a row whose columns carry $\theta_u\,$, $\theta_w$ and $-\theta_v$ encodes $\theta_v=\theta_u+\theta_w\,$. For~\hyperlink{res-A}{(A)}, one obstacle remains: a box point is not a Hilbert basis element if it is the sum of two nonzero box points, which happens when their coefficients add without a carry, i.e.\ no sum reaches $1$. We force a carry by making the supports of any two nonzero solutions intersect, using an additive encoder built from sparse neighbor-sum maps, in the spirit of linear-time codes~\cite{Spielman}, followed by a random ``support connector'' graph. Result~\hyperlink{res-D}{(D)} uses a unimodular graph representation and half-integral edge covers, and~\hyperlink{res-E}{(E)} uses the \emph{profiles} of vertex expressions on initial faces, which record the total coefficient of the prefixes containing each coordinate; their half-integrality gives the $2$-SAT criteria, and in general they encode integer programs. The first indecomposable holes of large degree in~\hyperlink{res-F}{(F)} come from a signed odd-cycle model realized by explicit growth orders.

\par\medskip\noindent\emph{Related work.}
Every finite abelian group is the residue group of some lattice simplex~\cite[Theorem~2.3]{BH}. New in \hyperlink{res-A}{(A)} are the realization on a \emph{facet} of three staircase simplices, membership in the ambient Hilbert basis, and joint control of degrees and clearing sets.

The closest relative of our results is the theorem of Billera and Sarangarajan~\cite{BS}: every $0/1$ polytope is affinely equivalent to a face of an asymmetric traveling salesman polytope $T_s\,$. Its cube directions are $\pm1$ vectors on disjoint cycles, so it preserves induced lattices and, after homogenization, vertex monoids, Hilbert bases and clearing sets. Our distinction in~\hyperlink{res-C}{(C)} is realization by three staircase simplices, any two of which have a totally unimodular lifted vertex matrix, in dimension polynomial in $N$ and $r$, with a facet option; \hyperlink{res-A}{(A)} further gives the optimal group-dependent dimension. In a similar spirit, Adiprasito and Padrol~\cite{AP} showed that every polytope with algebraic coordinates is projectively equivalent to a face of a projectively unique polytope, encoding field arithmetic by von Staudt constructions; ours encode group arithmetic by additions.

Two-way transportation matrices are totally unimodular, whereas De Loera and Onn~\cite{DLO} represent every rational polytope by a slim three-way transportation polytope: a similar two/three threshold. The series lift used for \hyperlink{res-C}{(C)} is a generalized Lawrence construction~\cite{PTV1,PTV2}; we also identify induced lattices, indecomposable holes and clearing sets on a face. Intersecting codes~\cite{BSS} give a variant of~\hyperlink{res-A}{(A)} for elementary abelian groups.

For \hyperlink{res-D}{(D)}, odd-cycle normality criteria for edge rings and relatives~\cite{OH,SVV,BL,FM} supply the normalization mechanism. Our contribution is a unimodular realization of every triangle-strip graph cone as the \emph{ambient} cone of three staircase simplices, with its subword and face descriptions.

For \hyperlink{res-E}{(E)}, deciding whether a set of integer vectors is a Hilbert basis of its cone is coNP-complete~\cite{Pap}; for edge polytopes of connected graphs, IDP is the odd cycle condition~\cite{OH}, whose failure is NP-complete to detect~\cite{GKPT}.

\par\medskip\noindent\emph{Organization.}
Section~\ref{sec-foundations} fixes notation and proves the basic upper bounds, including total unimodularity for two simplices and the bound on the residue order of a face with few indecomposable holes. Section~\ref{sec-realization} turns ordered disjoint supports into faces, introduces the incidence construction, realizes every clearing set at degree three, and encodes bounded linear systems by profiles of vertex expressions. Section~\ref{sec-cyclic} treats cyclic residue groups, faces with few indecomposable holes, and the two/three threshold. Section~\ref{sec-universality} proves arithmetic universality, Section~\ref{sec-idp} transfers IDP of $0/1$ polytopes, and Section~\ref{sec-ambient} studies initial faces and the complexity of splitting. Section~\ref{sec-failures} studies how IDP fails in the whole hull: in large degree, after long delays, from dimension six on, and away from simplex faces. Section~\ref{sec-words} classifies word disjunctions.

\section{Preliminaries and basic bounds}\label{sec-foundations}
This section fixes notation for disjunctive hulls of staircase simplices and their cones, and proves the upper bounds used throughout: total unimodularity for two simplices, bounds on determinants, degrees and residue groups, and the bound on the residue order of a simplex face in terms of its number of indecomposable holes.

\subsection{Disjunctive hulls, lattices, and clearing sets}
\begin{figure}[t]
\centering
\begin{tikzpicture}[x=0.72cm,y=0.72cm,scale=1.598]
 \begin{scope}[x=0.97cm,y=0.52cm]
 \foreach \x in {0,...,3} \foreach \y in {0,...,3} {\fill[gray!55] (\x,\y) circle (0.9pt);}
 \fill[blue!12] (0,0)--(1,0)--(2,3)--(1,3)--cycle;
 \draw[thick] (0,0)--(1,0); \draw[thick] (0,0)--(1,3);
 \draw[dashed] (1,0)--(2,3)--(1,3);
 \draw[->,gray] (1,0)--(3.2,0); \draw[->,gray] (1,3)--(1.2,3.6);
 \fill (0,0) circle (1.6pt); \fill (1,0) circle (1.8pt); \fill (1,3) circle (1.8pt);
 \node[below,font=\scriptsize] at (1,0) {$r_1$}; \node[left,font=\scriptsize] at (1,3) {$r_2$};
 \fill[blue!70!black] (1,1) circle (2.2pt); \fill[blue!70!black] (1,2) circle (2.2pt);
 \node[left,font=\scriptsize] at (0.98,1) {$H_1$}; \node[right,font=\scriptsize] at (1.02,2) {$H_2$};
 \node[below left,font=\scriptsize] at (0,0) {$0$};
 \end{scope}
\node[anchor=east,font=\scriptsize] at (6.45,1.35) {$kH_1$ splits};
\draw (6.75,1.35) circle (2.3pt);
\draw (7.17,1.35) circle (2.3pt);
\fill[blue!70!black] (7.59,1.35) circle (2.3pt);
\draw (8.01,1.35) circle (2.3pt);
\draw (8.43,1.35) circle (2.3pt);
\fill[blue!70!black] (8.85,1.35) circle (2.3pt);
\draw (9.27,1.35) circle (2.3pt);
\draw (9.69,1.35) circle (2.3pt);
\fill[blue!70!black] (10.11,1.35) circle (2.3pt);
\draw (10.53,1.35) circle (2.3pt);
\draw (10.95,1.35) circle (2.3pt);
\fill[blue!70!black] (11.37,1.35) circle (2.3pt);
\draw (11.79,1.35) circle (2.3pt);
\draw (12.21,1.35) circle (2.3pt);
\fill[blue!70!black] (12.63,1.35) circle (2.3pt);
\node[anchor=east,font=\scriptsize] at (6.45,0.30) {$\langle4,5\rangle$};
\draw (6.75,0.30) circle (2.3pt);
\node[font=\tiny,gray] at (6.75,-0.08) {1};
\draw (7.17,0.30) circle (2.3pt);
\node[font=\tiny,gray] at (7.17,-0.08) {2};
\draw (7.59,0.30) circle (2.3pt);
\node[font=\tiny,gray] at (7.59,-0.08) {3};
\fill[blue!70!black] (8.01,0.30) circle (2.3pt);
\node[font=\tiny,gray] at (8.01,-0.08) {4};
\fill[blue!70!black] (8.43,0.30) circle (2.3pt);
\node[font=\tiny,gray] at (8.43,-0.08) {5};
\draw (8.85,0.30) circle (2.3pt);
\node[font=\tiny,gray] at (8.85,-0.08) {6};
\draw (9.27,0.30) circle (2.3pt);
\node[font=\tiny,gray] at (9.27,-0.08) {7};
\fill[blue!70!black] (9.69,0.30) circle (2.3pt);
\node[font=\tiny,gray] at (9.69,-0.08) {8};
\fill[blue!70!black] (10.11,0.30) circle (2.3pt);
\node[font=\tiny,gray] at (10.11,-0.08) {9};
\fill[blue!70!black] (10.53,0.30) circle (2.3pt);
\node[font=\tiny,gray] at (10.53,-0.08) {10};
\draw (10.95,0.30) circle (2.3pt);
\node[font=\tiny,gray] at (10.95,-0.08) {11};
\fill[blue!70!black] (11.37,0.30) circle (2.3pt);
\node[font=\tiny,gray] at (11.37,-0.08) {12};
\fill[blue!70!black] (11.79,0.30) circle (2.3pt);
\node[font=\tiny,gray] at (11.79,-0.08) {13};
\fill[blue!70!black] (12.21,0.30) circle (2.3pt);
\node[font=\tiny,gray] at (12.21,-0.08) {14};
\fill[blue!70!black] (12.63,0.30) circle (2.3pt);
\node[font=\tiny,gray] at (12.63,-0.08) {15};
\end{tikzpicture}
\caption{\emph{Left:} the cone over $r_1=(1,0)$, $r_2=(1,3)$ has residue group $\Z/3\Z$ and box points $H_1=\frac23r_1+\frac13r_2\,$, $H_2=\frac13r_1+\frac23r_2$ (in the shaded half-open parallelepiped), both indecomposable holes; $2H_1=H_2+r_1$ and $3H_1=2r_1+r_2\,$, so $H_1$ has clearing order $3$. \emph{Right:} the $k\le15$ for which $kh$ splits (filled), for $h=H_1$ and for an indecomposable hole with clearing set $\langle4,5\rangle\setminus\{0\}$ (Theorem~\ref{thm-clearsets}), which has no clearing order.}
\label{fig-defs}
\end{figure}
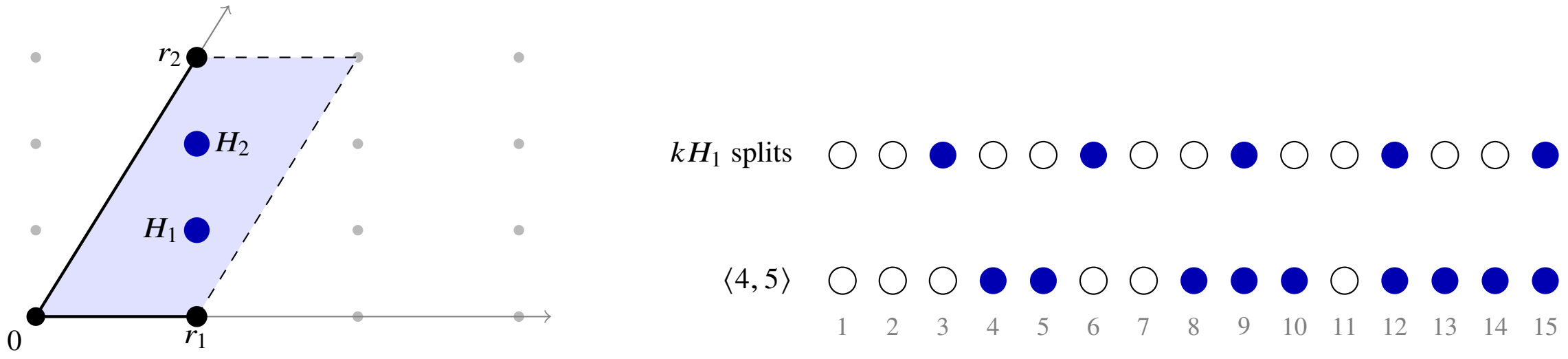

For $\sigma\in S_d\,$, the staircase simplex
\[
 \mathcal P_\sigma=\{x\in[0,1]^d:x_{\sigma(1)}\le\cdots\le x_{\sigma(d)}\}
\]
has as vertices the indicator vectors of the suffixes $\{\sigma(j),\ldots,\sigma(d)\}$ of $\sigma$, including the empty suffix. Equivalently, its \emph{growth order} $\gamma:=(\sigma(d),\sigma(d-1),\ldots,\sigma(1))$, the reverse of $\sigma$, switches the coordinates on one at a time, and the vertices are the indicator vectors of the \emph{prefixes} of $\gamma$; we use both descriptions. The coordinates of $x$ are the \emph{spatial} coordinates, and $d$ is the spatial dimension. We sometimes call these simplices \emph{full} staircase simplices, to stress that every growth order contains all $d$ coordinates, in contrast with the prefix chains on shorter lists in Section~\ref{sec-ambient}. Given $n+1$ such simplices $\mathcal P_0\,,\ldots,\mathcal P_n\,$, called the \emph{constituents}, put $\mathbf e_0:=\mathbf 0\in\R^n$ and let $\mathbf e_1\,,\ldots,\mathbf e_n$ be the unit vectors. Their disjunctive hull, lifted vertices, and homogenized cone are
\[
 \mathcal D:=\conv\bigcup_{i=0}^n(\mathcal P_i\times\{\mathbf e_i\}),\qquad
 V:=\{(v,\mathbf e_i\,,1):v\in\operatorname{vert}(\mathcal P_i)\},\qquad K:=\cone V.
\]
We use the natural coordinates $(x,z,t)$: $z_i$ indicates $\mathcal P_i$ and $t-\sum_i z_i$ indicates $\mathcal P_0\,$. Write $S_K:=\Z_{\ge0}V$ for the \emph{vertex monoid} and $\overline S_K:=K\cap\Z^{d+n+1}$ for its saturation. A real (integral) \emph{vertex expression} of a point of $K$ writes it as a nonnegative real (integral) combination of lifted vertices; a lattice point of $K$ \emph{splits} if it has an integral vertex expression, that is, if it lies in $S_K\,$. The lattice points of $\mathcal D$ are the points $(v,\mathbf e_i)$ with $v$ a vertex of $\mathcal P_i$ (Lemma~\ref{lem-full}), whose homogenizations $(v,\mathbf e_i\,,1)$ form $V$; thus $\mathcal D$ has the integer decomposition property (IDP) if and only if $S_K=\overline S_K\,$. The Hilbert basis $\Hilb(K)$ is the finite set of nonzero points of $\overline S_K$ that are not sums of two nonzero points of $\overline S_K\,$. We call a point of $\overline S_K$ \emph{indecomposable} if it lies in $\Hilb(K)$; the elements of $\Hilb(K)\setminus V$ are the \emph{indecomposable holes}. We reserve ``primitive'' for integer vectors and forms. All cones here are polyhedral, so every face $F$ of $K$ or $\mathcal D$ is exposed: it is where some linear form that is nonnegative on $V$ vanishes, an \emph{exposing form} of $F$. For a face $F$ of $\mathcal D$, let $K_F$ be the cone generated by the lifted vertices on $F$, its \emph{rays}, taken in its induced lattice; the \emph{ray lattice} of $F$ is the lattice generated by its rays. An \emph{indecomposable hole of $F$} is an element of $\Hilb(K_F)$ other than a lifted vertex on $F$; by Lemma~\ref{lem-face} it is an indecomposable hole of $K$. Degree means $t$. For $h\in\overline S_K$ define
\[
 \mathcal C_K(h):=\{q\in\Z_{>0}:qh\in S_K\},\qquad
 \delta_K(h):=\min\mathcal C_K(h),
\]
the \emph{clearing set} and the \emph{least clearing integer} of $h$. If $\mathcal C_K(h)=q\Z_{>0}$ for some $q$, we call $q$ the \emph{clearing order} of $h$; then $\delta_K(h)=q$. The vertex monoid $S_K$ is \emph{seminormal} if every $h\in\overline S_K$ with $2h,3h\in S_K$ lies in $S_K\,$.
For a simplex face $F$ with independent homogenized rays $r_1\,,\ldots,r_s\,$, its \emph{residue group} is
\[
 G_F:=\Lambda_F/\sum_j\Z r_j\,,\qquad
 \Lambda_F:=\spanR\{r_j\}\cap\Z^{d+n+1}.
\]
Every point of the span of $F$ is $\sum_ja_jr_j$ for unique \emph{ray coefficients} $a_j\,$, and each residue $g$ has a unique representative $H_g$ in the half-open box $\{\sum_j a_jr_j:0\le a_j<1\}$. We call $H_g$ the \emph{box point} of $g$, and the coefficients $a_j$ its \emph{box coefficients}. The normalized volume $\Vol(F)$ of a simplex face, measured in its induced lattice, equals $|G_F|$, and $F$ is \emph{unimodular} if $G_F$ is trivial. We write $\{u\}:=u-\lfloor u\rfloor$ for the fractional part of a real number $u$. All face lattices below are induced lattices. Figure~\ref{fig-defs} illustrates these notions.

For $n\ge0$ and $d\ge1$, let $H_n(d)$ be the maximum number of indecomposable holes for $n+1$ staircase simplices $\mathcal P_0\,,\ldots,\mathcal P_n$ in $\R^d$, and let $\Delta_n(d)$ be the maximum of $\delta_K(y)$ over the lattice points of their cones. Repeated constituents are allowed.
Put $\T:=\R/\Z$ and $\N:=\Z_{\ge0}\,$. The change to multiplicities
\begin{equation}\label{eq-multiplicity}
 \lambda_0:=t-\sum_{i=1}^n z_i\,,\qquad \lambda_i:=z_i\quad(1\le i\le n)
\end{equation}
is an integral linear bijection, with inverse $z_i=\lambda_i$ and $t=\sum_i\lambda_i\,$; we call $(x,\lambda)$ the \emph{multiplicity coordinates}. In these coordinates the lifted vertices are $(v,\eps_i)$, where $\eps_i$ is a unit vector of $\R^{n+1}$. The cone condition is
$\lambda\ge0,\qquad x\in\sum_{i=0}^{n}\lambda_i\mathcal P_i\,.$
This change preserves induced lattices, semigroups, and Hilbert bases. We will prove total unimodularity directly in the natural coordinates, because arbitrary unimodular row changes need not preserve it.

\begin{lemma}[Generation and lattice points]\label{lem-full}
The vectors $V$ generate $\Z^{d+n+1}$ as a group. The cone $K$ is full dimensional, $\dim \mathcal D=d+n$, and the lattice points of $\mathcal D$ are precisely the points $(v,\mathbf e_i)$ with $v$ a vertex of $\mathcal P_i\,$.
\end{lemma}
\begin{proof}
The zero vertex of $\mathcal P_0$ gives $(0,0,1)$. Subtracting it from the zero vertex in each other constituent gives the unit vectors in the $z$-coordinates. Differences between consecutive suffix vertices give every spatial unit vector. This proves the lattice and dimension assertions. A lattice point of $\mathcal D$ has $z\ge0$ and $\sum_i z_i\le1$, hence $z=\mathbf e_i$ for some $i\in\{0,\ldots,n\}$. Its spatial part is then a lattice point of $\mathcal P_i\,$, hence a suffix indicator.
\end{proof}

\begin{proposition}[Two simplices]\label{prop-two}
For at most two staircase simplices, the matrix with columns $V$ is totally unimodular. Every triangulation of $\mathcal D$ using its lifted vertices is unimodular, and $\mathcal D$ has IDP.
\end{proposition}
\begin{proof}
IDP is known more generally for Cayley sums of two order polytopes \cite[Proposition~6.1]{HOT}. The network-matrix argument below is the one used in \cite[Theorem~4.7(a)]{Stair} to show that the facet normals of a sum of two staircase simplices have entries in $\{0,\pm1\}$; see also \cite[Corollary~6.2]{Stair}.
Consider any square submatrix. Within each constituent, order its selected columns by increasing suffix inclusion and subtract each preceding selected column from its successor, working from last to first. These operations preserve the determinant. The resulting spatial supports are disjoint within each constituent. In every $z$-row and in the $t$-row, only the first selected column in a constituent can retain a nonzero entry. Negate the second constituent's columns. Each row now has at most one $+1$ and at most one $-1$. A row with at most one nonzero entry permits determinant expansion and induction. If all rows have two entries, the all-ones vector is in the kernel. Thus every square minor is $0$ or $\pm1$.

By Lemma~\ref{lem-full}, nonzero maximal determinants are normalized simplex volumes. Every vertex triangulation is therefore unimodular, and its cones decompose every lattice point of the cone into integer ray sums.
\end{proof}

\begin{lemma}[Faces]\label{lem-face}
Let $K'$ be a face of $K$, with its induced lattice. Then $\Hilb(K')\subseteq\Hilb(K)$, and every lattice point $h$ of $K'$ has the same clearing set in $K'$ as in $K$. If its rays $r_1\,,\ldots,r_s$ are independent and extend to a basis $B$ of columns from $V$, inclusion induces an injection
\[
 \bigl(\spanR\{r_j\}\cap\Z^{d+n+1}\bigr)\big/\sum_j\Z r_j
 \ \hookrightarrow\ \Z^{d+n+1}/B\Z^{d+n+1}.
\]
\end{lemma}
\begin{proof}
A nonnegative exposing form vanishing at $h=u+v$ must vanish at both summands. Thus any decomposition, or expression by vertices, stays in $K'$, so exposure also preserves clearing sets. For the injection, an integer combination of columns of $B$ lying in the face span has zero coefficients on all the other basis columns, by independence.
\end{proof}

\subsection{Determinant, degree, and rank bounds}
We now bound the arithmetic that three staircase simplices can carry: the determinants of bases of lifted vertices, the degrees of indecomposable holes, and the order and rank of residue groups of simplex faces.
For $r\ge0$, let $B(r)$ be the largest product of positive integers with sum at most $r$ (the empty product being one). It is nondecreasing in $r$, and $B(r)\le3^{r/3}$, because $\prod_ja_j=\prod_j(a_j^{1/a_j})^{a_j}$ and $a^{1/a}\le3^{1/3}$ for every positive integer $a$.
\begin{lemma}[Three-class determinant bound {\cite[Lemma~4.6]{Stair}}]\label{lem-threeclass}
Let $A_0$ be a nonsingular $k\times k$ $0/1$ matrix, $k\ge1$, whose columns are partitioned into three classes such that the supports of the columns within each class are pairwise disjoint. Then $|\det A_0|\le B(k-1)$.
\end{lemma}
\begin{proof}
Within a class the supports are disjoint, so the sum of the columns of a class is the indicator of the set of rows it covers. If two classes both covered all $k$ rows, their column sums would agree, a linear dependence; so at most one class covers all rows, and we may choose a class $J$ covering a set of at most $k-1$ rows. Let its columns have supports $S_1\,,\ldots,S_s\,$, of sizes $a_1\,,\ldots,a_s\ge1$ (a zero column would make $A_0$ singular), so $\sum_ja_j\le k-1$. Expand $\det A_0$ by the generalized Laplace expansion along the columns of $J$: it is a signed sum, over $s$-sets $R$ of rows, of $\det A_0[R,J]\cdot\det A_0[R^c,J^c]$. Because the supports $S_j$ are disjoint, $\det A_0[R,J]$ vanishes unless $R$ contains exactly one row of each $S_j\,$, in which case it is $\pm1$; there are $\prod_ja_j$ such $R$. The columns of $J^c$ form two classes with disjoint supports in each, so every row of $A_0[R^c,J^c]$ has at most one entry one in each class. After negating the columns of one class, every row has at most one $+1$ and at most one $-1$, so this matrix is totally unimodular by the row argument in the proof of Proposition~\ref{prop-two}, and its determinant is $0$ or $\pm1$. Therefore $|\det A_0|\le\prod_ja_j\le B(k-1)$.
\end{proof}
\begin{proposition}[Determinant bounds]\label{prop-det}
A nonsingular $(d+n+1)\times(d+n+1)$ matrix $B$ of columns in $V$ satisfies
\[
 |\det B|\le (n+1)^{d/2}.
\]
For three constituents,
\[
 |\det B|\le B(d-1)\le3^{(d-1)/3}<1.4423^{\,d}.
\]
\begin{equation}\label{eq-det24}
 |G_F|\le3^{(d-1)/3}\quad\text{for every simplex face $F$ of three staircase simplices.}
\end{equation}
\end{proposition}
\begin{proof}
Use \eqref{eq-multiplicity}. Each constituent occurs among the columns, or its multiplicity row would vanish. Keep one reference column in each constituent and replace later selected columns by consecutive differences. Expansion in the multiplicity rows reduces the determinant to a $d\times d$ $0/1$ matrix $A_0$ whose rows have at most $n+1$ ones, because the difference columns within one constituent are indicators of pairwise disjoint intervals of its coordinate order. Hadamard's inequality gives the first bound, as for facet normals in \cite[Theorem~4.7(d)]{Stair}. For $n=2$ the columns of $A_0$ fall into three classes with pairwise disjoint supports within each class, so Lemma~\ref{lem-threeclass} gives $|\det A_0|\le B(d-1)\le3^{(d-1)/3}$. For a simplex face, extend its independent rays to a full basis from $V$. Lemma~\ref{lem-face} embeds its residue group into the full basis quotient, whose order is the determinant, proving~\eqref{eq-det24}.
\end{proof}
\begin{proposition}[Counting, clearing and degree bounds]\label{prop-bounds}
For $n+1$ staircase simplices in $\R^d$, with $\beta_n(d):=(n+1)^{d/2}$ for $n\ne2$ and $\beta_2(d):=B(d-1)$,
\begin{equation}\label{eq-count-bound}
 H_n(d)\le\binom{(n+1)(d+1)}{d+n+1}\bigl(\lfloor\beta_n(d)\rfloor-1\bigr),
 \qquad \delta_K(y)\le\lfloor\beta_n(d)\rfloor.
\end{equation}
For $n=2$, $\beta_2(d)\le3^{(d-1)/3}$. Every indecomposable hole has degree in $\{3,\ldots,d\}$.
\end{proposition}
\begin{proof}
Conical Carath\'eodory expresses a lattice point of the cone using independent rays, which can be extended to a full basis $B$ from $V$. For an indecomposable point other than a ray, all coefficients are below one: otherwise subtracting that ray leaves a nonzero lattice point of the cone. Such a point is in the half-open parallelepiped of $B$, which has $|\det B|$ lattice points, including zero. There are at most $\binom{(n+1)(d+1)}{d+n+1}$ bases, proving the count.

If $k_i$ basis rays come from constituent $i$, then $k_i\ge1$, $\sum_i k_i=d+n+1$, and the integral multiplicities of a box point satisfy $0\le\lambda_i\le k_i-1$. Thus its degree is at most $d$. Degree-one points are rays. Degree-two points use at most two constituents and decompose by Proposition~\ref{prop-two}.

Finally, for integer $y\in\cone B$, Cramer's rule gives $|\det B|B^{-1}y\in\Z_{\ge0}^{d+n+1}$. Hence $\delta_K(y)\le|\det B|$, and Proposition~\ref{prop-det} applies.
\end{proof}

\begin{proposition}[Group-rank bound]\label{prop-grouprank}
Let $F$ be a simplex face of a disjunctive hull of three staircase simplices in $\R^d$, with independent rays $r_1\,,\ldots,r_s\,$, of which $s_i$ come from constituent $i$, and with nontrivial residue group $G_F\,$. For every prime $p$, the $p$-rank $r_p(G_F)=\dim_{\F_p}G_F[p]$ satisfies
\[
 r_p(G_F)\le s_i-1\qquad(i=0,1,2).
\]
Consequently $s\ge3\mu(G_F)+3$ and $d\ge3\mu(G_F)+1$, where $\mu(G_F)=\max_pr_p(G_F)$ is the minimal number of generators. In particular, the box point of every nonzero residue has positive multiplicity in each constituent.
\end{proposition}
\begin{proof}
Identify $G_F$ with the group of box points, and project a box point to its vector of ray coefficients in constituent $i$. This projection is injective: a representative with all coefficients zero in constituent $i$ is a lattice point in the cone of independent rays from at most two constituents, hence an integer combination of them by Proposition~\ref{prop-two}, and coefficients that are integral and below one vanish. The coefficients in constituent $i$ sum to the integer multiplicity $\lambda_i\,$, so the projection lands in the sum-zero subtorus of $\T^{s_i}$, which is isomorphic to $\T^{s_i-1}$ and has $p$-torsion of $\F_p$-dimension $s_i-1$. This proves the inequality, and because $G_F\ne0$ every $s_i$ is positive. Summing over $i$ gives $s\ge3\mu(G_F)+3$. A proper face of $K$ has at most $\dim K-1=d+2$ independent rays, and $F$ is proper because $\mathcal D$ has $3(d+1)>d+3$ vertices; hence $d\ge s-2\ge3\mu(G_F)+1$.
\end{proof}

For three constituents we collect these bounds.

\begin{proposition}[Upper bounds]\label{prop-upper}
With $B$ as above, for three staircase simplices in $\R^d$:
\begin{enumerate}
\item every basis $U\subseteq V$ of $\R^{d+3}$ has $|\det U|\le B(d-1)\le3^{(d-1)/3}$; hence $|G_F|\le3^{(d-1)/3}$ for every simplex face $F$, and $\delta_K(h)\le3^{(d-1)/3}$ for every $h\in\overline S_K$;
\item if $G_F$ is nontrivial and $F$ has $s_i$ rays from constituent $i$, then the $p$-rank of $G_F$ is at most $s_i-1$ for every prime $p$ and every $i$; hence $d\ge3\mu(G_F)+1$, where $\mu$ is the minimal number of generators, and every nonzero $H_g$ has positive multiplicity in every constituent;
\item every indecomposable hole has degree between $3$ and $d$.
\end{enumerate}
\end{proposition}
\begin{proof}
Part~(1) is Propositions~\ref{prop-det} and~\ref{prop-bounds}, part~(2) is Proposition~\ref{prop-grouprank}, and part~(3) is Proposition~\ref{prop-bounds}.
\end{proof}

The number of indecomposable holes of a face also limits its residue group: isolating arithmetic costs dimension.
\begin{theorem}[Atoms versus residue order]\label{thm-atoms}
Let $F$ be a simplex face of three staircase simplices in $\R^d$ with $s_i$ rays from constituent $i$, residue group of order $Q>1$, and exactly $r$ indecomposable holes. Put $M_F:=\min_is_i-1$. Then
\[
 Q\le\binom{M_F+r}{r}\le\binom{\lfloor(d-1)/3\rfloor+r}{r}.
\]
In particular, if $F$ has a single indecomposable hole then $d\ge3Q-2$, and if $d\le C\log_2Q$ then $r\ge c_C\log_2Q$ for a constant $c_C>0$ depending only on $C$.
\end{theorem}
\begin{proof}
A box point has ray coefficients in $[0,1)$, so its integral multiplicity in constituent $i$ is at most $s_i-1$. A decomposition of a box point into Hilbert basis elements of $K_F$ contains no ray, because a ray summand would give a coefficient at least one, so all summands are indecomposable holes of $F$, which are themselves box points of nonzero residues and have positive multiplicity in every constituent by Proposition~\ref{prop-grouprank}. Hence such a decomposition has at most $M_F$ summands. Choosing one decomposition for each of the $Q$ box points gives an injective map to vectors in $\Z_{\ge0}^r$ with sum at most $M_F\,$, of which there are $\binom{M_F+r}r$. A simplex face is proper, because $\mathcal D$ has $3(d+1)>d+3$ vertices, so $\sum_is_i\le d+2$ and $M_F\le\lfloor(d-1)/3\rfloor$. For $r=1$ this reads $Q\le\lfloor(d-1)/3\rfloor+1$. For the last assertion put $L=\log_2Q$ and $z=r/L$; the estimate $\binom{M_F+r}{r}\le(e(M_F+r)/r)^r$ gives $1\le z\log_2(e(C/3+z)/z)$, whose right side tends to zero as $z\to0$.
\end{proof}
The logarithmic order of $r$ is attained by Proposition~\ref{prop-binary-normal}, and the bound $d\ge3Q-2$, for every $Q\ge3$, by Theorem~\ref{thm-onehole}. A single indecomposable hole is in fact very rigid.
\begin{corollary}[Faces with one indecomposable hole]\label{cor-onehole-rigid}
Suppose a simplex face $F$ as in Theorem~\ref{thm-atoms} has exactly one indecomposable hole $h$. Then $G_F$ is cyclic of order $Q$ generated by the class of $h$, $\mathcal C_K(h)=Q\Z_{>0}\,$, every nonzero ray coefficient of $h$ equals $1/Q$, and the rays in its support number $Q\deg h\ge3Q$. The face spanned by them is the unique minimal face of $F$ failing IDP, and all its proper faces are unimodular.
\end{corollary}
\begin{proof}
The saturated face monoid is generated by the rays and $h$, and every box point is a sum of these without a ray summand, so the box points are exactly the multiples $nh$ that lie in the half-open box; hence $G_F$ is generated by $[h]$, and its order $Q$ is the clearing order of $h$ by Lemma~\ref{lem-face} and independence. If $h$ has a coefficient $a/Q$ with $a\ge2$, fewer than $Q$ multiples of $h$ lie in the box; so all nonzero coefficients are $1/Q$, and the degree, which is the sum of the coefficients, equals the size of the support divided by $Q$. On that support every nonzero residue has all coefficients nonzero, so deleting a ray kills the residue group, and a face failing IDP must contain $h$ and hence the whole support. Finally $\deg h\ge3$ by Proposition~\ref{prop-bounds}.
\end{proof}
\section{From disjoint supports to simplex faces}\label{sec-realization}
This section develops the main tool of the paper: it turns a $0/1$ matrix with three classes of disjoint column supports into a face of three staircase simplices whose residue group, box points and clearing sets can be read off from the matrix. Subsection~\ref{sec-ordering} orders the supports and identifies the face. Subsection~\ref{sec-incidence} introduces scalar presentations and the incidence construction, which produces such matrices from systems of additions; square presentations with zero real kernel give simplex facets, and Subsection~\ref{sec-extended} allows the arbitrary kernels needed in Section~\ref{sec-idp}. Without the rank condition, the same exposure gives every clearing set at degree three (Subsection~\ref{sec-clearsets}), and profiles of vertex expressions encode bounded linear systems (Subsection~\ref{sec-bounded}).

\subsection{Ordering supports}\label{sec-ordering}
The following lemma summarizes the construction; Lemmas~\ref{lem-growth} and~\ref{lem-shared} below state and prove it in detail.

\begin{lemma}[From supports to simplex faces]\label{lem-from-supports}
Let $A\in\{0,1\}^{p\times s}$ have full column rank and columns in three nonempty ordered classes with disjoint supports within each class. Let $\alpha\in\Z^p$ satisfy $A^T\alpha=0$, taking both signs on every support and, for each class, a positive value outside the union of its supports. For each class, list the outside nonnegative coordinates (zeros first), then its supports in order (negative, zero, positive weights within each), then the remaining negative coordinates. These three lists are growth orders of three staircase simplices in $\R^p$ with the following properties.
\begin{enumerate}
\item The support-boundary prefixes are exactly the vertices on a face $F$; its $s+3$ rays are independent, and $G_F\cong\Gamma_A:=\{f\in\T^s:Af\equiv0\}$; if $p=s+1$, $F$ is a facet.
\item For $f\in\Gamma_A$ with ordered values $f_{i,1}\,,\ldots,f_{i,m_i}\in[0,1)$ on the $m_i$ columns of class $i$, the box coefficients of $H_f$ in class $i$ are $\{-f_{i,1}\},\ \{f_{i,b}-f_{i,b+1}\}\ (1\le b<m_i),\ f_{i,m_i}\,$, and $\mathcal C_K(H_f)=\ord(f)\Z_{>0}\,$.
\item With weights in $\{-1,0,1\}$, one positive and one negative coordinate per support, one outside positive coordinate per class and no remaining negative coordinates, the exposing form is primitive and takes only the values $0,1$ on $V$.
\end{enumerate}
\end{lemma}
\par
\begin{figure}[b]
\centering
\begin{tikzpicture}[x=0.47cm,y=0.75cm,scale=1.375]
\begin{scope}[xshift=0.00cm]
 \draw[gray!60] (0,0)--(7,0); \draw[gray!60] (0,0)--(0,1.25);
 \node[anchor=east,font=\scriptsize] at (0,0) {$0$}; \node[anchor=east,font=\scriptsize] at (0,1) {$1$};
 \draw[gray!40,dashed] (0,1)--(7,1);
 \draw[thick] (0,0) -- (1,1) -- (2,0) -- (3,1) -- (4,0) -- (5,0) -- (6,1) -- (7,0);
 \filldraw[fill=white,draw=black] (0,0) circle (1.6pt);
 \fill[blue!70!black] (1,1) circle (2.2pt);
 \filldraw[fill=white,draw=black] (2,0) circle (1.6pt);
 \fill[blue!70!black] (3,1) circle (2.2pt);
 \filldraw[fill=white,draw=black] (4,0) circle (1.6pt);
 \filldraw[fill=white,draw=black] (5,0) circle (1.6pt);
 \fill[blue!70!black] (6,1) circle (2.2pt);
 \filldraw[fill=white,draw=black] (7,0) circle (1.6pt);
 \node[font=\scriptsize] at (0.5,-0.42) {$Q_1$};
 \node[font=\scriptsize] at (1.5,-0.42) {$R_2$};
 \node[font=\scriptsize] at (2.5,-0.42) {$P_2$};
 \node[font=\scriptsize] at (3.5,-0.42) {$P_1$};
 \node[font=\scriptsize] at (4.5,-0.42) {$R_1$};
 \node[font=\scriptsize] at (5.5,-0.42) {$g$};
 \node[font=\scriptsize] at (6.5,-0.42) {$Q_2$};
 \draw[gray] (0.080000000000000,-0.75)--(0.080000000000000,-0.85)--(0.920000000000000,-0.85)--(0.920000000000000,-0.75);
 \draw[gray] (1.080000000000000,-0.75)--(1.080000000000000,-0.85)--(2.920000000000000,-0.85)--(2.920000000000000,-0.75);
 \draw[gray] (3.080000000000000,-0.75)--(3.080000000000000,-0.85)--(5.920000000000000,-0.85)--(5.920000000000000,-0.75);
 \draw[gray] (6.080000000000000,-0.75)--(6.080000000000000,-0.85)--(6.920000000000000,-0.85)--(6.920000000000000,-0.75);
 \node[font=\small] at (3.5,1.65) {$\mathcal P_{\sigma_0}$};
\end{scope}
\begin{scope}[xshift=4.0cm]
 \draw[gray!60] (0,0)--(7,0); \draw[gray!60] (0,0)--(0,1.25);
 \node[anchor=east,font=\scriptsize] at (0,0) {$0$}; \node[anchor=east,font=\scriptsize] at (0,1) {$1$};
 \draw[gray!40,dashed] (0,1)--(7,1);
 \draw[thick] (0,0) -- (1,1) -- (2,0) -- (3,1) -- (4,0) -- (5,0) -- (6,1) -- (7,0);
 \filldraw[fill=white,draw=black] (0,0) circle (1.6pt);
 \fill[blue!70!black] (1,1) circle (2.2pt);
 \filldraw[fill=white,draw=black] (2,0) circle (1.6pt);
 \fill[blue!70!black] (3,1) circle (2.2pt);
 \filldraw[fill=white,draw=black] (4,0) circle (1.6pt);
 \filldraw[fill=white,draw=black] (5,0) circle (1.6pt);
 \fill[blue!70!black] (6,1) circle (2.2pt);
 \filldraw[fill=white,draw=black] (7,0) circle (1.6pt);
 \node[font=\scriptsize] at (0.5,-0.42) {$g$};
 \node[font=\scriptsize] at (1.5,-0.42) {$P_1$};
 \node[font=\scriptsize] at (2.5,-0.42) {$Q_1$};
 \node[font=\scriptsize] at (3.5,-0.42) {$Q_2$};
 \node[font=\scriptsize] at (4.5,-0.42) {$R_1$};
 \node[font=\scriptsize] at (5.5,-0.42) {$P_2$};
 \node[font=\scriptsize] at (6.5,-0.42) {$R_2$};
 \draw[gray] (0.080000000000000,-0.75)--(0.080000000000000,-0.85)--(0.920000000000000,-0.85)--(0.920000000000000,-0.75);
 \draw[gray] (1.080000000000000,-0.75)--(1.080000000000000,-0.85)--(2.920000000000000,-0.85)--(2.920000000000000,-0.75);
 \draw[gray] (3.080000000000000,-0.75)--(3.080000000000000,-0.85)--(5.920000000000000,-0.85)--(5.920000000000000,-0.75);
 \draw[gray] (6.080000000000000,-0.75)--(6.080000000000000,-0.85)--(6.920000000000000,-0.85)--(6.920000000000000,-0.75);
 \node[font=\small] at (3.5,1.65) {$\mathcal P_{\sigma_1}$};
\end{scope}
\begin{scope}[xshift=8.0cm]
 \draw[gray!60] (0,0)--(7,0); \draw[gray!60] (0,0)--(0,1.25);
 \node[anchor=east,font=\scriptsize] at (0,0) {$0$}; \node[anchor=east,font=\scriptsize] at (0,1) {$1$};
 \draw[gray!40,dashed] (0,1)--(7,1);
 \draw[thick] (0,0) -- (1,1) -- (2,0) -- (3,1) -- (4,0) -- (5,0) -- (6,1) -- (7,0);
 \filldraw[fill=white,draw=black] (0,0) circle (1.6pt);
 \fill[blue!70!black] (1,1) circle (2.2pt);
 \filldraw[fill=white,draw=black] (2,0) circle (1.6pt);
 \fill[blue!70!black] (3,1) circle (2.2pt);
 \filldraw[fill=white,draw=black] (4,0) circle (1.6pt);
 \filldraw[fill=white,draw=black] (5,0) circle (1.6pt);
 \fill[blue!70!black] (6,1) circle (2.2pt);
 \filldraw[fill=white,draw=black] (7,0) circle (1.6pt);
 \node[font=\scriptsize] at (0.5,-0.42) {$P_2$};
 \node[font=\scriptsize] at (1.5,-0.42) {$Q_2$};
 \node[font=\scriptsize] at (2.5,-0.42) {$g$};
 \node[font=\scriptsize] at (3.5,-0.42) {$R_2$};
 \node[font=\scriptsize] at (4.5,-0.42) {$R_1$};
 \node[font=\scriptsize] at (5.5,-0.42) {$Q_1$};
 \node[font=\scriptsize] at (6.5,-0.42) {$P_1$};
 \draw[gray] (0.080000000000000,-0.75)--(0.080000000000000,-0.85)--(0.920000000000000,-0.85)--(0.920000000000000,-0.75);
 \draw[gray] (1.080000000000000,-0.75)--(1.080000000000000,-0.85)--(2.920000000000000,-0.85)--(2.920000000000000,-0.75);
 \draw[gray] (3.080000000000000,-0.75)--(3.080000000000000,-0.85)--(5.920000000000000,-0.85)--(5.920000000000000,-0.75);
 \draw[gray] (6.080000000000000,-0.75)--(6.080000000000000,-0.85)--(6.920000000000000,-0.85)--(6.920000000000000,-0.75);
 \node[font=\small] at (3.5,1.65) {$\mathcal P_{\sigma_2}$};
\end{scope}
\end{tikzpicture}
\caption{Lemmas~\ref{lem-growth} and~\ref{lem-shared} on the facet $F_3$ of Example~\ref{ex-q3}. Each panel plots $\alpha^Tx$ along one growth order, with coordinates below and blocks bracketed. The value rises to $1$ on the initial block, and each zero-sum block dips and returns to $1$ at its end, so $t-\alpha^Tx$ vanishes exactly at the nine filled prefixes, the rays of $F_3\,$. Consecutive filled prefixes differ by a block indicator, a column of $A$.}
\label{fig-supports}
\end{figure}
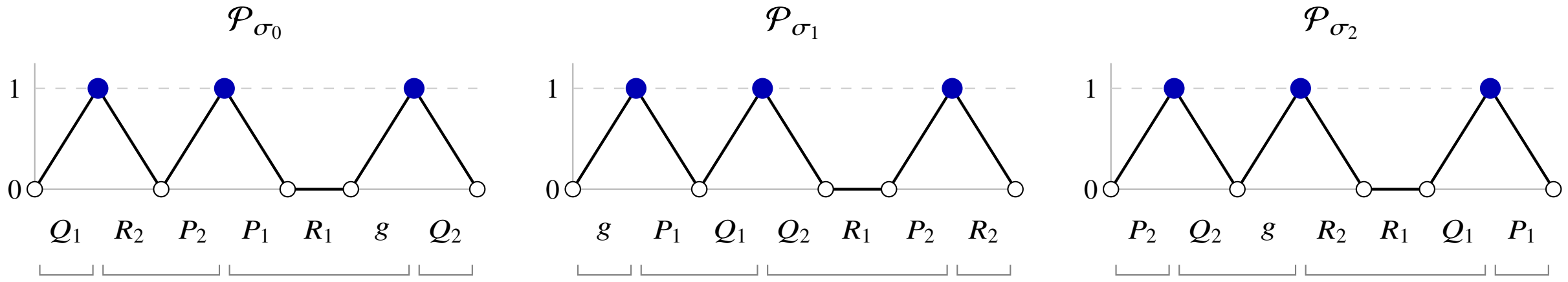

\medskip\noindent\emph{Idea of the proof.} Let $M_i$ be the total outside positive weight of class $i$. Along the growth order of class $i$, the value $\alpha^Tx$ rises to $M_i$ on the outside positive coordinates. Each support then first decreases the value strictly and returns to $M_i$ exactly at its last coordinate, because the support has total weight zero; the remaining negative coordinates decrease it again. Hence $\alpha^Tx\le M_0\lambda_0+M_1\lambda_1+M_2\lambda_2$ is valid on $V$ and tight exactly at the support boundaries. Replacing consecutive rays of $F$ by their differences and shearing by the first prefix on $F$ of each class turns the ray matrix into $\operatorname{diag}(A,I_3)$, which gives independence and $G_F\cong(\spanR A\cap\Z^p)/A\Z^s\cong\Gamma_A\,$. Writing a face point in ray coordinates, the difference columns carry the tail sums of the ray coefficients, and inverting gives~(2). Clearing sets follow from exposure (Lemma~\ref{lem-face}) and independence. Figure~\ref{fig-supports} shows the argument on a small facet.

\begin{proof}
Lemma~\ref{lem-growth} gives the face and part~(3); Lemma~\ref{lem-shared} gives independence, the residue group, facet status and part~(2).
\end{proof}

In detail, the first lemma turns a $0/1$ matrix with three classes of disjoint supports into three growth orders and an exposing form; the second identifies the face it exposes, its residue group and its box points.
\begin{lemma}[Ordering supports]\label{lem-growth}
Partition the columns of $A\in\{0,1\}^{p\times s}$ into three nonempty ordered classes $J_0\,,J_1\,,J_2\,$, with disjoint supports within each class. Suppose $\alpha\in\Z^p$ satisfies $A^T\alpha=0$, every support has both a positive- and a negative-weight coordinate, and every class has a positive-weight coordinate outside its supports.

For class $i$, list the outside zero-weight coordinates, the outside positive-weight coordinates, its ordered supports (negative, then zero, then positive within each support), and finally the remaining negative coordinates. Reverse this list to define $\sigma_i\,$. Let $M_i>0$ be the total outside positive weight. Then
\begin{equation}\label{eq-expose}
 \alpha^Tx\le M_0(t-z_1-z_2)+M_1z_1+M_2z_2
\end{equation}
exposes exactly the support endpoints, with adjacent support endpoints identified.

If, additionally, each support has exactly one positive- and one negative-weight coordinate, both of absolute weight one, all weights lie in $\{0,1,-1\}$, each class has exactly one outside positive coordinate, and every negative coordinate lies in a support of every class, then $M_i=1$ and $\ell=t-\alpha^Tx$ is primitive with $\ell(v)\in\{0,1\}$ for every $v\in V$.
\end{lemma}
\begin{proof}
Suffix indicators grow by adding coordinates in the displayed list order. The objective increases to $M_i$ at the last outside positive coordinate. Traversing a support first decreases it strictly, then returns to $M_i$ exactly at the last positive coordinate, because that support has total weight zero. The remaining negative coordinates decrease it. These are precisely the maximizing vertices. Under the additional hypotheses there is one downward and one upward unit step per support, the outside positive step is one, and there are no remaining negative coordinates. All objective values are therefore zero or one. The coefficient of $t$ in $\ell$ is one.
\end{proof}
This chooses coordinate orders to realize prescribed zero-$\alpha$-sum blocks, reversing the maximizing-block description of \cite[Lemma~2.1]{Stair}; see also \cite[Lemmas~2.4 and~2.5]{Cayley}. Figure~\ref{fig-supports} shows the argument on the facet of Example~\ref{ex-q3}.

\begin{lemma}[Shared endpoints]\label{lem-shared}
Under the weight hypotheses of Lemma~\ref{lem-growth}, assume $A$ has full column rank $s$. Put $\Gamma_A:=\{f\in\T^s:Af=0\text{ in }\T^p\}$, the \emph{torus kernel} of $A$. Write the ordered columns of class $i$ as $A_{i,1}\,,\ldots,A_{i,q_i}$ and its first maximizing suffix as $u_i\,$. The exposed rays, in multiplicity coordinates, are
\[
 r_{i,b}=\left(u_i+\sum_{j\le b}A_{i,j}\,,\eps_i\right),\qquad 0\le b\le q_i\,.
\]
They are independent, define a simplex face of dimension $s+2$, and have residue group $\Gamma_A\,$. If $p=s+1$, the face is a facet. For $f\in\Gamma_A$ with values $f_{i,1}\,,\ldots,f_{i,q_i}\in[0,1)$ on the ordered columns of class $i$, its box point has ray coefficients
\begin{equation}\label{eq-shared-coeff}
 c_{i,0}=\{-f_{i,1}\},\quad
 c_{i,b}=\{f_{i,b}-f_{i,b+1}\}\ (1\le b<q_i),\quad
 c_{i,q_i}=f_{i,q_i}\,.
\end{equation}
Every integer face point is uniquely a nonnegative integral ray sum plus one box point $H_f\,$, and a face point with box point $H_f$ has clearing set $\ord_{\Gamma_A}(f)\Z_{>0}$ in the full cone.
\end{lemma}
\begin{proof}
Lemma~\ref{lem-growth} identifies the exposed vertices. Replace successive rays by consecutive differences and shear $x\mapsto x-\sum_i\lambda_i u_i\,$. These integral operations reduce the ray matrix to $\operatorname{diag}(A,I_3)$, up to column order. Its rank is $s+3$ and its saturation quotient is
$(\spanR A\cap\Z^p)/A\Z^s\cong\Gamma_A\,,$
where $f\mapsto A\tilde f+A\Z^s$, for any real lift $\tilde f\in\R^s$ of $f$, induces the isomorphism. This proves independence, the group assertion, and facet status when $p=s+1$.

For ray coefficients $c_{i,b}$ the coefficients of the difference columns are $F_{i,j}=\sum_{b\ge j}c_{i,b}$ and $\lambda_i=\sum_bc_{i,b}\,$. At a lattice point, $f=\{F\}$ belongs to $\Gamma_A\,$. Integral multiplicities and these congruences give exactly \eqref{eq-shared-coeff}; conversely those coefficients produce a lattice point. Flooring the unique ray coefficients gives the stated decomposition. Exposure confines every semigroup expression for a multiple of a box point to these rays. Its coefficients become integral exactly when the multiple kills $f$.
\end{proof}
\begin{remark}[Relation to facet normals of Minkowski sums]\label{rem-stair}
For $p=s+1$, the supports in Lemma~\ref{lem-shared} satisfy the three-color zero-sum criterion for Minkowski-sum facet normals \cite[Theorem~3.1]{Stair}; inequality~\eqref{eq-expose} is the corresponding full optimal big-$M$ lifting, which by \cite[Proposition~3.1]{Cayley} describes the face, and for $p=s+1$ a facet, of the disjunctive hull. The additional information here is the induced residue group, box coefficients and clearing sets.
\end{remark}
\subsection{Scalar presentations and the incidence construction}\label{sec-incidence}

A \emph{scalar presentation} consists of labeled scalar variables, its \emph{nodes}, and equations
of the following forms:
\[
 \theta_v=\theta_u+\theta_w\,,\qquad
 \theta_u+\theta_w=0,\qquad \theta_u=0.
\]
The equations are \emph{gates}, \emph{terminal constraints} and \emph{zero constraints}, respectively; a designated zero node is written $\theta_\star\,$. The kernel is the solution group over the displayed coefficient group: the \emph{real kernel} over $\R$ and the \emph{torus kernel} over $\T=\R/\Z$. For an integer matrix $R$ with $m$ columns we likewise write $\ker_\T R:=\{f\in\T^m:Rf=0\}$ for its torus kernel.

\begin{lemma}[Common incidence construction]\label{univ-incidence}
Let a scalar presentation have $Q\ge1$ nodes and $Q$ equations of the allowed forms, with repeated operands permitted. Assume that its real kernel is zero, and write $G$ for its finite torus kernel, possibly trivial. There is a $0/1$ matrix $A$ with $9Q+1$ rows and $9Q$ columns in three classes of $3Q$ columns each, with disjoint supports within every class. It has full column rank and torus kernel $G$. Each class has two columns with values $+\theta_j\,,-\theta_j$ for each scalar node $j$, and $Q$ further designated columns whose torus values vanish identically.

For \emph{any} chosen order within each class, the support-order construction gives three staircase simplices in spatial dimension $d=9Q+1$, with a simplex facet $F$ having $9Q+3$ independent homogenized rays. The induced-lattice residue group is $G$ and $\Vol(F)=|G|$. A primitive integral exposing form takes only the values zero and one on the ambient lifted vertices.

For any chosen orders, the box coefficients are given by~\eqref{eq-shared-coeff} with $q_i=3Q$. 
Every lattice point of the facet cone is a unique nonnegative integral ray sum plus one box point. Every nonzero representative $H_g$ has exact ambient clearing set $\mathcal C_K(H_g)=\ord(g)\Z_{>0}\,$.
\end{lemma}
\begin{proof}
\emph{The incidence configuration.}
For node $j$, introduce positive graph vertices $p_{0,j}\,,p_{1,j}\,,p_{2,j}$ and
negative graph vertices $n_{0,j}\,,n_{1,j}\,,n_{2,j}\,$. Form the following properly
three-colored six-cycle. The signs identify the values of its edge variables
in the incidence kernel:
\[
\begin{array}{c|cc}
 \text{color}&\text{edge of value }\theta_j&\text{edge of value }-\theta_j\\\hline
 A&p_{0,j}n_{0,j}&p_{2,j}n_{1,j}\\
 B&p_{2,j}n_{2,j}&p_{1,j}n_{0,j}\\
 C&p_{1,j}n_{1,j}&p_{0,j}n_{2,j}\,.
\end{array}
\]
Join $p_{1,j}$ to $n_{2,j+1}$ by a color-$A$ bridge for $1\le j<Q$.
Add a positive leaf vertex $L$ joined to $n_{1,1}$ by a color-$B$ edge.
There are now $6Q+1$ vertices and $7Q$ edges. Add $2Q$ completion edges:
\[
\begin{array}{c|l}
 A&p_{1,Q}n_{2,1}\\
 B&p_{0,j}n_{1,j}\quad(2\le j\le Q)\\
 C&p_{2,j}n_{0,j}\quad(1\le j\le Q).
\end{array}
\]
The coloring remains proper. Every negative vertex has all three colors. The
only missing positive incidences are colors $A,C$ at $L$ and color $B$ at
$p_{0,1}\,$. In particular, each color class contains $3Q$ edges, and the graph
has $s=9Q$ edges.

Let $A$ be the $0/1$ matrix with one column for each edge. Its rows are the
unsigned graph incidence rows, one singleton row for each completion edge, and
one row for each scalar equation. A gate
$\theta_v=\theta_u+\theta_w$ uses the positive color-$A$ copy of $u$, the
positive color-$C$ copy of $w$, and the negative color-$B$ copy of $v$.
A terminal constraint uses the first two copies; each zero constraint uses the positive color-$A$ copy of its scalar. These columns are distinct even for repeated operands.
Every row uses at most one column of each color, so column supports within a
color class are disjoint. There are
$p=(6Q+1)+2Q+Q=9Q+1=s+1$
rows.

\emph{The kernel.}
Over $\R$ or $\T$, singleton rows first force the completion-edge variables
to vanish. With those edges removed, each bridge is a cut edge. On one side of
its cut, sum the positive incidence equations minus the negative incidence
equations. Internal edges cancel, leaving the bridge variable zero. The leaf
equation forces its edge to zero as well. The remaining six-cycle equations
give alternating edge values $\theta_j\,,-\theta_j\,$. The last $Q$ rows are
exactly the scalar presentation. Therefore $A$ has full column rank, and
$\Gamma_A:=\{f\in\T^s:Af=0\}\cong G$.
In every torus solution the bridge, leaf and completion edges are zero. There
are exactly $Q$ such columns in each color: the $Q-1$ bridges and one
completion in $A$, the leaf and $Q-1$ completions in $B$, and the $Q$
completions in $C$.

\emph{Arbitrary orders and the facet.}
Give positive incidence rows weight $+1$, negative incidence rows weight $-1$, and all other rows weight zero. Let $\alpha$ be this weight vector. Every column has one positive and one negative incidence, so $A^T\alpha=0$. Every class misses exactly one positive row and no negative row. For any order of its columns, apply Lemma~\ref{lem-growth}, using a fixed row order for ties. Its strengthened hypotheses hold, giving full coordinate permutations and the primitive $0/1$ exposing form
$\ell(x,z,t)=t-\alpha^Tx$.
Because $A$ has $s=9Q$ independent columns and $p=s+1$ rows, Lemma~\ref{lem-shared} gives a simplex facet with $s+3=9Q+3$ rays, residue group $G$ and volume $|G|$. It also gives the box coefficients, flooring decomposition and exact ambient clearing sets.
\end{proof}

\subsection{Arbitrary scalar kernels}\label{sec-extended}
Lemma~\ref{univ-incidence} requires a square presentation with zero real kernel. The IDP transfer of Section~\ref{sec-idp} needs presentations with arbitrary kernels; the same construction still gives a face, whose rays may now be dependent.
\begin{proposition}[Extended incidence construction]\label{prop-extended-incidence}
Let a scalar presentation have $Q\ge1$ nodes and $E$ equations of the forms allowed in Lemma~\ref{univ-incidence}. No finiteness or rank assumption is imposed. The incidence construction gives a $0/1$ matrix $A$ with $9Q$ columns in three disjoint-support classes and $8Q+E+1$ rows, whose torus kernel is the scalar kernel. For any order in each class, it gives three staircase simplices in dimension $8Q+E+1$ with a face having the $9Q+3$ support endpoints as lifted vertices and a primitive exposing form taking values zero and one on $V$. The torus kernel of the endpoint matrix is obtained from the ordered column values of $\ker_\T A$ by the backward-difference formula~\eqref{eq-shared-coeff}. The endpoints need not be independent.
\end{proposition}
\begin{proof}
Use the same graph, completion edges, and scalar-equation rows as in Lemma~\ref{univ-incidence}. Only the equation-row count changes, giving $(6Q+1)+2Q+E=8Q+E+1$ rows. The kernel argument first kills the completion edges, bridges, and leaf, then recovers the alternating scalar values on each six-cycle. It did not use finiteness, so the torus kernel is precisely the scalar kernel. The weight vector and proper coloring are unchanged, hence Lemma~\ref{lem-growth} exposes the same endpoint vertices with primitive $0/1$ slack. Consecutive endpoint differences and the shear $x\mapsto x-\sum_i\lambda_i u_i$ reduce the endpoint matrix to $\operatorname{diag}(A,I_3)$. The inverse column transformation gives the backward-difference formula, without requiring independence.
\end{proof}

The support orders in Lemma~\ref{lem-shared} are free, and the choice matters. Ordering every class by nonincreasing values of a fixed nonzero $f\in\Gamma_A$ makes the coefficients in~\eqref{eq-shared-coeff} sum to one in each class, so $H_f$ has degree three (Proposition~\ref{prop-sort}); with an addition chain, this gives every clearing order $q\ge2$ at degree three on a simplex face in dimension at most $5\ell(q)+1$, where $\ell(q)$ is the length of a shortest \emph{addition chain} $1=a_0<a_1<\cdots<a_\ell=q$, in which every term after $a_0$ is a sum of two earlier ones (Corollary~\ref{cor-sort}). Without the rank assumption, the same exposure gives every clearing set.

\subsection{Every clearing set at degree three}\label{sec-clearsets}
We use multiplicity coordinates $(x,\lambda_0\,,\lambda_1\,,\lambda_2)$ for three staircase simplices with growth orders $\gamma_0\,,\gamma_1\,,\gamma_2\,$, and cumulative profiles (used again in Subsection~\ref{sec-dichotomy}): for a real vertex expression of $(x,\lambda)$, the profile $f_i(j)$ is the total coefficient of the prefixes of $\gamma_i$ containing coordinate $j$. Profiles are nonincreasing along $\gamma_i\,$, take values in $[0,\lambda_i]$, satisfy $\sum_if_i(j)=x_j\,$, and determine the expression by consecutive differences, including $\lambda_i-f_i$ at the first coordinate and $f_i$ at the last. Integral expressions correspond to integral profiles. Every lattice point of degree at most two splits (Proposition~\ref{prop-two}), so a nonsplitting lattice point of degree three is an ambient Hilbert basis element.

\begin{lemma}[Ordered supports without a rank condition]\label{lem-ordered}
Let $M$ be a $0/1$ matrix with $r$ rows and no zero column, whose columns are split into three nonempty ordered classes with pairwise disjoint supports within each class. There are three staircase simplices in $\R^{2r+1}$ with a face whose lifted vertices in class $i$ are
\[
 r_{i,b}=\Bigl(u_i+\sum_{j=1}^bA_{i,j}\,,\ \eps_i\Bigr)\quad(0\le b\le n_i),\qquad A=\begin{pmatrix}M\\ M\\ 0\end{pmatrix},
\]
where $n_i$ is the number of columns of class $i$, $A_{i,j}$ is its $j$th column and $u_i$ is a $0/1$ vector. For $b\in\Z^r$ and $k\ge0$, the vertex expressions on this face of
\[
 \Bigl(k\textstyle\sum_iu_i+(kb,kb,0),\ k,k,k\Bigr)
\]
correspond bijectively and affinely to the solutions of $Mf=kb$ whose entries in each class are nonincreasing in the column order and lie in $[0,k]$; for integral $k$ the correspondence preserves integrality.
\end{lemma}
\begin{proof}
Give the first $r$ rows of $A$ weight $+1$, the next $r$ weight $-1$ and the last row weight $+1$. In class $i$, list the positive coordinates outside its column supports, then the supports in their order, each with its negative coordinates before its positive ones, and then the remaining negative coordinates. The prefix weight first reaches its maximum just before the supports, at the prefix $u_i\,$; in each support it drops strictly below the maximum and returns to it at the last coordinate, because each support has total weight zero. Hence, with $\alpha$ the weight vector and $m_i=\alpha^\top u_i\,$, the form $\sum_im_i\lambda_i-\alpha^\top x$ is nonnegative on $K$ and vanishes exactly at the displayed endpoints. For endpoint coefficients $c_{i,0}\,,\ldots,c_{i,n_i}$ put $f_{i,j}=\sum_{b\ge j}c_{i,b}\,$. Then $x=k\sum_iu_i+Af$ and $(c_{i,0}\,,\ldots,c_{i,n_i})=(k-f_{i,1}\,,f_{i,1}-f_{i,2}\,,\ldots,f_{i,n_i})$, which gives the correspondence. The same identities hold for signed coefficients.
\end{proof}

\begin{theorem}[Every clearing set at degree three]\label{thm-clearsets}
Let $s\ge1$ and $q_1\,,\ldots,q_s\in\Z_{\ge2}\,$. There are three staircase simplices in dimension $O(\sum_j\log q_j)$ and an ambient Hilbert basis element $h$ with multiplicities $(1,1,1)$ such that
\[
 \mathcal C_K(h)=\Bigl\{\textstyle\sum_jq_ja_j:a\in\Z_{\ge0}^s\Bigr\}\setminus\{0\}.
\]
For every $k\ge1$ the real vertex expressions of $kh$ are in affine bijection with $\{a\in\R_{\ge0}^s:\sum_jq_ja_j=k\}$, integral expressions corresponding to integral $a$. The order of $h$ modulo the lattice generated by the rays of its minimal face is $\gcd(q_1\,,\ldots,q_s)$.
\end{theorem}
\par\medskip\noindent\emph{Idea of the proof.} As in Lemma~\ref{lem-from-supports}, but without full rank, the expressions of $kh$ are the solutions of the row equations at scale $k$, nonincreasing in each class with values in $[0,k]$ (Lemma~\ref{lem-ordered}). Rows following addition chains for the $q_j$ make path increments multiples of $a_j\,$; a terminal row forces $\sum_jq_ja_j=k$, and monotonicity is $a\ge0$.

\begin{proof}
For block $i$ put $u_i'=\lfloor q_i/2\rfloor$ and $v_i'=q_i-u_i'$. Choose an addition chain from $1$ to $u_i'$ (the one-term chain $(1)$ if $u_i'=1$), let $I_i$ be its set of levels together with $0$, and put $J_i=I_i\cup\{v_i'\}$; the binary method gives $|I_i|=O(\log q_i)$. The columns of $M$ are variables in three classes $A,B,C$: in block $i$, $A_{i,b}$ for $b\in I_i$ and $B_{i,b}\,,C_{i,b}$ for $b\in J_i\,$. In classes $A$ and $C$ identify the last variable of each block with the first of the next, $A_{i,u_i'}=A_{i+1,0}$ and $C_{i,v_i'}=C_{i+1,0}\,$; write $A_*\,,C_*$ for the last variables of these paths. Class $B$ uses distinct variables in every block and a terminal variable $B_*\,$. The rows of $M$, with right-hand sides at scale $k$, are
\begin{align}
 A_{1,0}=C_{1,0}=B_*&=0,\label{eq-cs-zero}\\
 A_{i,0}+B_{i,0}+C_{i,0}&=k,\label{eq-cs-base}\\
 A_{i,b}+B_{i,b}+C_{i,0}=A_{i,0}+B_{i,b}+C_{i,b}&=k\quad(b\in I_i\setminus\{0\}),\label{eq-cs-synch}\\
 A_{i,p}+B_{i,b}+C_{i,r}&=k\quad(b=p+r\text{ in the chain}),\label{eq-cs-gate}\\
 A_{i,u_i'}+B_{i,v_i'}+C_{i,1}=A_{i,0}+B_{i,v_i'}+C_{i,v_i'}&=k\quad(q_i\text{ odd}),\label{eq-cs-odd}\\
 A_*+B_*+C_*&=k,\label{eq-cs-term}
\end{align}
each equality in~\eqref{eq-cs-zero}, \eqref{eq-cs-synch} and~\eqref{eq-cs-odd} being a separate row. To complete class $A$, append a reverse copy of the whole $C$ path, each $C$ replaced by a new variable $\overline C$ of class $A$ with a row $\overline C+C=k$, identifying $\overline C_*$ with $A_*\,$; complete class $C$ symmetrically by the reversed $A$ path, identifying $\overline A_*$ with $C_*\,$, and include the central row $A_*+C_*=k$ once. Order classes $A$ and $C$ in reverse of their completed paths, and class $B$ by blocks, increasing in level, with $B_*$ last. Every row uses at most one variable per class, every column occurs in a row, and there are $O(\sum_i\log q_i)$ rows.

Let $a_i=A_{i,1}-A_{i,0}\,$, $S_i^A=\sum_{j<i}u_j'a_j$ and $S_i^C=\sum_{j<i}v_j'a_j\,$. By~\eqref{eq-cs-base} and~\eqref{eq-cs-synch} the $A$ and $C$ increments above the block starts agree at every level, by~\eqref{eq-cs-gate} minus~\eqref{eq-cs-synch} they are additive along the chain, and~\eqref{eq-cs-odd} gives the extra $C$ level of an odd block. With~\eqref{eq-cs-zero} and the shared endpoints,
\[
 A_{i,b}=S_i^A+ba_i\,,\qquad C_{i,b}=S_i^C+ba_i\,,\qquad B_{i,b}=k-S_i^A-S_i^C-ba_i\,,
\]
and~\eqref{eq-cs-term} becomes $\sum_iq_ia_i=k$; the appended variables are determined by the copy rows. Conversely these formulas solve all rows when $\sum_iq_ia_i=k$. If all $a_i\ge0$, the completed $A$ and $C$ paths increase from $0$ to $k$ and the $B$ path decreases from $k$ to $0$, because $v_i'<q_i$ and block $i$ lowers $B$ by $q_ia_i\,$; so the order and bounds of Lemma~\ref{lem-ordered} hold. Conversely the order of $A_{i,0}\,,A_{i,1}$ forces $a_i\ge0$, and the values are integral exactly when the $a_i$ are, as $a_i=A_{i,1}-A_{i,0}\,$.

Apply Lemma~\ref{lem-ordered} with $b$ the $0/1$ vector of right-hand sides at $k=1$. The point $h$ at $k=1$ lies on the face, because $\sum_iq_ia_i=1$ has a nonnegative real solution, and exposure confines every expression of $kh$ to the face. This gives the bijections. As $q_i\ge2$, $h$ does not split; having degree three, it is an ambient Hilbert basis element.

In each class order the first and last values are $k$ and $0$, so the first and last endpoint coefficients vanish; taking every $a_i>0$ at $k=1$ makes all other endpoint coefficients positive, so these other endpoints are exactly the rays of the minimal face of $h$. Signed integer expressions on these rays satisfy the same linear equations, hence correspond to integer $a$ with $\sum_iq_ia_i=k$, and conversely. Thus $kh$ lies in the ray lattice exactly when $\gcd(q_1\,,\ldots,q_s)$ divides $k$.
\end{proof}

\begin{corollary}[All clearing sets]\label{cor-clearsets-all}
The sets $\{0\}\cup\mathcal C_K(h)$, for indecomposable holes $h$ of degree three, are exactly the nonzero submonoids of $\Z_{\ge0}$ that omit $1$. The vertex monoid need not be seminormal, and the largest $k$ with $kh\notin S_K$ can be $q^2-q-1$ in dimension $O(\log q)$ for an indecomposable hole that is trivial modulo the ray lattice of its minimal face.
\end{corollary}
\begin{proof}
Clearing sets are closed under addition, and nonempty because a rational expression can be scaled to an integral one; they omit $1$ for an indecomposable hole. A submonoid of $\Z_{\ge0}$ with least positive element $q$ is generated by $q$ and the least element of each other residue class modulo $q$ that it meets, a finite list of integers at least two; apply Theorem~\ref{thm-clearsets}. For $\langle2,3\rangle$ the indecomposable hole $h$ satisfies $h\notin S_K$ and $2h,3h\in S_K\,$, while $V$ generates the full lattice, so $S_K$ is not seminormal~\cite{BLR}. For $\langle q,q+1\rangle$, the integer $aq+b$ with $0\le b<q$ lies in the monoid exactly when $a\ge b$, so its largest gap is $q^2-q-1$, and $\gcd(q,q+1)=1$.
\end{proof}

\subsection{Bounded systems as profile equations}\label{sec-bounded}
The profiles of Subsection~\ref{sec-clearsets} also encode general bounded linear systems with at most three terms per equation. In this construction a coordinate that lies in exactly one of the three growth orders before a final guard coordinate $g$ is a \emph{singleton coordinate}, and one that lies in exactly two is a \emph{pair coordinate}; the lemma is used in Sections~\ref{sec-universality} and~\ref{sec-ambient}.

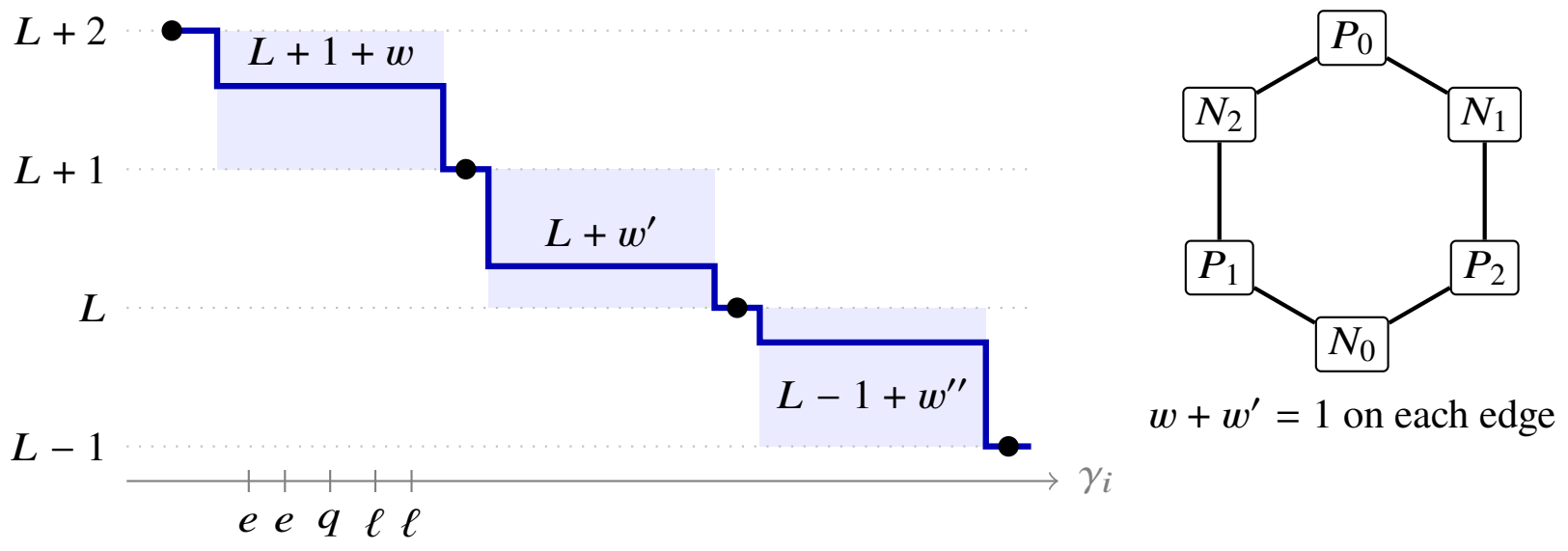
\begin{figure}[t]
\centering
\begin{tikzpicture}[x=0.62cm,y=0.95cm,font=\small]
 \fill[blue!8] (1,2) rectangle (3.5,3);
 \fill[blue!8] (4,1) rectangle (6.5,2);
 \fill[blue!8] (7,0) rectangle (9.5,1);
 \foreach \y/\l in {3/L+2,2/L+1,1/L,0/L-1} {\draw[gray!50,dotted] (0,\y)--(10,\y); \node[anchor=east,font=\scriptsize] at (0,\y) {$\l$};}
 \draw[->,gray] (0,-0.25)--(10.3,-0.25) node[right,font=\scriptsize] {$\gamma_i$};
 \draw[very thick,blue!70!black] (0.5,3)--(1,3)--(1,2.6)--(3.5,2.6)--(3.5,2)--(4,2)--(4,1.3)--(6.5,1.3)--(6.5,1)--(7,1)--(7,0.75)--(9.5,0.75)--(9.5,0)--(10,0);
 \foreach \x/\y in {0.5/3,3.75/2,6.75/1,9.75/0} {\fill (\x,\y) circle (2pt);}
 \foreach \x/\l in {1.35/e,1.75/e,2.25/q,2.75/\ell,3.15/\ell} {\draw[gray] (\x,-0.33)--(\x,-0.17); \node[font=\scriptsize] at (\x,-0.55) {$\l$};}
 \node[font=\scriptsize] at (2.25,2.85) {$L+1+w$};
 \node[font=\scriptsize] at (5.25,1.55) {$L+w'$};
 \node[font=\scriptsize] at (8.25,0.38) {$L-1+w''$};
 \begin{scope}[xshift=8.4cm,yshift=1.75cm,x=1cm,y=1cm]
  \foreach \n/\a/\l/\v in {P0/90/P_0/z_v,N1/30/N_1/1-z_v,P2/-30/P_2/z_v,N0/-90/N_0/1-z_v,P1/-150/P_1/z_v,N2/150/N_2/1-z_v} {
    \node[draw,rounded corners=1pt,inner sep=2pt,font=\scriptsize,fill=white] (\n) at (\a:1.05) {$\l$};
  }
  \draw[thick] (P0)--(N1)--(P2)--(N0)--(P1)--(N2)--(P0);
  \node[font=\scriptsize] at (0,-1.55) {$w+w'=1$ on each edge};
 \end{scope}
\end{tikzpicture}
\caption{Lemma~\ref{lem-bounded}. \emph{Left:} the profile of one class along its growth order, near three consecutive blocks. Singleton coordinates (dots) pin the profile to integer values between blocks, so a block followed by $L$ blocks takes values in $[L,L+1]$ (shaded). Early ($e$) and late ($\ell$) pair coordinates hold the profile constant across the block at $L+w$, and an equation coordinate ($q$) inside the block reads off $w$. \emph{Right:} the six blocks of a variable $z_v\,$, positive ($P_i\,$, in class $i$) and negative ($N_j\,$, in class $j\ne i$). Each edge carries an early and a late pair coordinate and imposes $w+w'=1$, so the positive blocks carry $z_v$ and the negative blocks $1-z_v\,$.}
\label{fig-profile}
\end{figure}

\begin{lemma}[Bounded systems as profile equations]\label{lem-bounded}
Let a feasible system on $n\ge1$ variables $z_v\in[0,1]$ consist of $p$ equations $\sum_{j=1}^{\ell}\epsilon_jz_{v_j}=b$ with $1\le\ell\le3$, $\epsilon_j\in\{\pm1\}$ and $b\in\Z$, a variable possibly occurring more than once in an equation. In polynomial time one can construct three staircase simplices in dimension $d=p+18n-2$, whose face $\{x_g=0\}$ is a facet, and a lattice point $h=(x,2n,2n,2n)$ of the cone over this facet, such that for every $k\in\Z_{>0}$ the real vertex expressions of $kh$ correspond bijectively to the solutions of the system with bounds $0\le z_v\le k$ and right-hand sides $kb$, integrality being preserved in both directions.
\end{lemma}
\begin{proof}
See Figure~\ref{fig-profile}. In each of three lists create a positive and a negative block for every variable, in a fixed order. A block followed by $L$ blocks gets the interval $[L,L+1]$ at scale one: between consecutive blocks put a singleton coordinate whose value is their common bound, with top and bottom bounds $2n$ and $0$. So a profile entry in a block is $L+w$ with $0\le w\le1$. For each variable join its positive block in class $i$ to its negative block in class $j$ for each $i\ne j$; these six edges form a cycle on its six blocks. For an edge between blocks with offsets $L,L'$, create an early and a late pair coordinate, both of value $L+L'+1$, and order every block as all its early coordinates, then its equation occurrences (possibly none), then all its late coordinates. Equal sums at the early and late coordinates of nonincreasing profiles force both profiles to be constant between them; the intervals overlap and cover the block, so its normalized value $w$ is constant, and each edge imposes $w+w'=1$. By connectivity the positive blocks of a variable carry a common $z_v$ and the negative blocks $1-z_v\,$. For an equation introduce one coordinate, assign its terms to distinct classes, put the coordinate into the block of $v_j$ of sign $\epsilon_j$ in the class of term $j$, and give it the value $b+\sum_jL_j+|\{j:\epsilon_j=-1\}|$; its equation is the prescribed one. Conversely a solution fills every block and singleton uniquely and monotonically. At scale $k$ the offsets and negative values scale accordingly, and integer profiles correspond to integer $z$. Finally append a guard $g$ to every list, followed by the missing coordinates; at $x_g=0$ exactly the original prefixes remain. Every coordinate other than $g$ lies in some list, so the consecutive prefix differences and the three empty prefixes span $\{x_g=0\}$, which is therefore a facet. There are $12n$ pair coordinates, $3(2n-1)$ singletons, $p$ equation coordinates and the guard.
\end{proof}

\section{Cyclic arithmetic, few indecomposable holes, and the two/three threshold}\label{sec-cyclic}
This section proves result~\hyperlink{res-B}{(B)} and the two/three threshold; it studies how much arithmetic a single face can isolate. Subsections~\ref{sec-cyclic-construction} and~\ref{sec-binary} build simplex faces with cyclic residue groups from addition chains, including the binary family. Subsections~\ref{sec-onehole}--\ref{sec-degq} determine what one, or as few as possible, indecomposable holes force on the dimension, the group and the degrees. Subsection~\ref{sec-threshold} proves the two/three threshold for the number of indecomposable holes and for least clearing integers. How IDP fails in the whole hull, beyond the designated faces, is the subject of Section~\ref{sec-failures}.

\subsection{The cyclic construction}\label{sec-cyclic-construction}
An \emph{addition chain} for $q\ge2$ is a sequence
\[
 1=a_0<a_1<\cdots<a_\ell=q,\qquad
 a_j=a_{u(j)}+a_{v(j)}\,,\quad u(j),v(j)<j.
\]
It is \emph{irredundant} if every nonterminal node is used in a later addition. Let $\ell(q)$ be the shortest chain length. Binary chains give $\ell(q)=O(\log q)$, and a shortest chain is irredundant.

\begin{theorem}[The cyclic construction]\label{thm-cyclic}
For every irredundant addition chain of length $\ell$ for $q\ge2$, three staircase simplices in spatial dimension at most $5\ell+1$ have a simplex face $F$ with $\dim F=3\ell+2$ and $G_F\cong\Z/q\Z$. The face has exactly $q-1$ indecomposable holes, all in the ambient Hilbert basis. For each nonzero residue $g$,
$\mathcal C_K(H_g)=\ord(g)\Z_{>0}\,.$
\end{theorem}
\begin{proof}
For $0\le j<\ell$ introduce variables $X_j\,,Y_j\,,Z_j\,$, assigned respectively to three column classes. Form a square $0/1$ matrix $B$ of order $3\ell$ with rows
\[
 R_j=X_j+Y_j\,,\qquad E_j=Y_j+Z_j\quad(0\le j<\ell),
\]
\[
 G_j=X_{u(j)}+Z_{v(j)}+Y_j\quad(1\le j<\ell),\qquad
 G_\ell=X_{u(\ell)}+Z_{v(\ell)}\,.
\]
The notation specifies the nonzero columns in a row. Even when operands coincide the $X$- and $Z$-columns are distinct. The equations over $\T$ give
\[
 X_j=Z_j=a_j\theta,\qquad Y_j=-a_j\theta,\qquad q\theta=0.
\]
Over $\R$ only the zero solution is possible. Thus $B$ is nonsingular and $\Gamma_B\cong\Z/q\Z$.

For $j>0$, count one operand deficit if $X_j$ occurs in no arithmetic row, and one if $Z_j$ occurs in no arithmetic row. Irredundancy permits at most one deficit per node, so the total $b$ satisfies $0\le b\le\ell-1$. Select all $\ell$ arithmetic rows, $R_0\,,E_0\,$, and $R_j$ or $E_j$ for a deficit of $X_j$ or $Z_j\,$, respectively. These $\ell+2+b$ rows cover every column: $Y_j$ occurs in $G_j$ for $j>0$, and every remaining column is used later or covered by its deficit row. The initial columns are covered by $R_0\,,E_0\,$.

Append duplicates of these rows, forming $A$ with $p=4\ell+2+b$ rows. Assign weight $+1$ to selected original rows, $-1$ to their duplicates, and zero elsewhere. Every column support has both signs and total weight zero. Supports within a class are disjoint. The positive rows $E_0\,,R_0\,,G_\ell$ omit respectively classes $X,Z,Y$, so every class has an outside positive coordinate.

Order the supports in each class by node index and apply Lemma~\ref{lem-growth}. To fix permutations completely, list the original rows as $R_0\,,E_0\,,\ldots,R_{\ell-1}\,,E_{\ell-1}\,,G_1\,,\ldots,G_\ell\,$, followed by the selected duplicates in that order, and preserve this row order within each part of every growth order. Lemma~\ref{lem-shared} gives a simplex face with $3\ell+3$ rays and residue group $\Z/q\Z$ in dimension $p\le5\ell+1$.

It remains to check indecomposability, because joining endpoints can change box points. Let $f_{i,j}$ be the fractional column value in class $i$ at node $j$. The first ray coefficient in class $i$ is $c_{i,0}=\{-f_{i,0}\}$. At a nonzero residue, $f_{0,0}=\theta$ and $f_{1,0}=\{-\theta\}$, so
$c_{0,0}+c_{1,0}=1.$
If a nonzero box point decomposed, its summands would themselves be nonzero box points, because all their ray coefficients are below one. The sum of the two displayed coefficients would be two in the sum and one in the proposed total, a contradiction. Flooring ray coefficients proves completeness of this face basis. Lemma~\ref{lem-face} gives ambient survival, and Lemma~\ref{lem-shared} gives the clearing sets.
\end{proof}

Let $D_3^{\rm den}(q)$ be the least spatial dimension of three staircase simplices whose ambient cone has a Hilbert basis element $h$ with $\mathcal C_K(h)=q\Z_{>0}\,$. A generator for a residue of order $q$ in Theorem~\ref{thm-cyclic} and the bound $q\le3^{(d-1)/3}$ of Proposition~\ref{prop-bounds} give
$\left\lceil1+3\log_3q\right\rceil\le D_3^{\rm den}(q)\le5\ell(q)+1,$
so $D_3^{\rm den}(q)=\Theta(\log q)$; for $q=2^k$ the two bounds are $\lceil1+3k\log_32\rceil\approx\lceil1+1.893k\rceil$ and $4k+2$.

The support orders in Lemma~\ref{lem-shared} are free, and a suitable order brings any prescribed residue down to the minimal degree.
\begin{proposition}[A prescribed residue in degree three]\label{prop-sort}
In Lemma~\ref{lem-shared}, let $f\in\Gamma_A$ be nonzero, with values $f_{i,b}\in[0,1)$, and order the columns of every class by nonincreasing values of $f$. Then the box point $H_f$ has multiplicities $(1,1,1)$, is an ambient Hilbert basis element of degree three, and has clearing set $\ord(f)\Z_{>0}\,$.
\end{proposition}
\begin{proof}
By Proposition~\ref{prop-grouprank}, $f$ is not identically zero on any class. In a class with values $f_1\ge\cdots\ge f_q$ and $f_1>0$, the coefficients~\eqref{eq-shared-coeff} are $1-f_1\,,\ f_1-f_2\,,\ \ldots,\ f_{q-1}-f_q\,,\ f_q\,$, which lie in $[0,1)$ and sum to one. So $H_f$ has multiplicity one in every class and degree three; it is not a vertex sum, by exposure and its nonintegral coefficients, and every lattice point of degree at most two involves at most two constituents and is a vertex sum (Proposition~\ref{prop-two}), so it is a Hilbert basis element. The clearing set is given by Lemma~\ref{lem-shared}.
\end{proof}
\begin{corollary}[Every clearing order in degree three]\label{cor-sort}
For every $q\ge2$ there are three staircase simplices in spatial dimension at most $5\ell(q)+1$, and $4k+2$ if $q=2^k$, whose cone has a Hilbert basis element $h$ with $\lambda(h)=(1,1,1)$, $\deg h=3$ and $\mathcal C_K(h)=q\Z_{>0}\,$.
\end{corollary}
\begin{proof}
Apply Proposition~\ref{prop-sort} to the matrix of the proof of Theorem~\ref{thm-cyclic} and the residue $\theta=1/q$; the weights and dimension do not depend on the order of the supports.
\end{proof}
Thus degree three, the least possible degree of an indecomposable hole, already carries arbitrary clearing orders, in contrast with the word disjunctions of Section~\ref{sec-words}, where every indecomposable hole has degree three and clears at two.

\subsection{The binary family}\label{sec-binary}
For the doubling chain the box points of the cyclic face can be listed completely, which gives exponentially many indecomposable holes in linear dimension.
For a lattice simplex $F$, write $h_F^*(\tau)$ for the degree enumerator of its box points; equivalently,
$\sum_{j\ge0}|jF\cap\Z^{d+2}|\tau^j=h_F^*(\tau)/(1-\tau)^{\dim F+1}$.

\begin{theorem}[The binary family]\label{thm-binary}
For every $k\ge1$ there are three staircase simplices in $\R^{4k+2}$ with a simplex face $F_k$ of dimension $3k+2$ and group $\Z/2^k\Z$. It has exactly $2^k-1$ indecomposable holes, all in $\Hilb(K)$.
A residue of exact order $2^r$ is uniquely $\theta=a/2^r$ with $a$ odd, $1\le a<2^r$. If $w(a)$ counts its binary ones, then
\begin{equation}\label{eq-binary-multiplicity}
 (z_1\,,z_2\,,t)(H_\theta)=\bigl(w(a),r+1-w(a),2r+2-w(a)\bigr),
 \qquad \mathcal C_K(H_\theta)=2^r\Z_{>0}\,.
\end{equation}
The largest degree of an indecomposable hole of this face is $2k+1$.
\end{theorem}
\begin{proof}
Use the doubling chain $a_j=2^j$ in Theorem~\ref{thm-cyclic}. Every nonterminal node is used as both operands of the next addition, so $b=0$ and $d=4k+2$.

For a residue $\theta=a/2^r$ put $f_j=\{2^j\theta\}$. These are nonzero exactly for $0\le j<r$, and $f_{r-1}=1/2$. In classes $X,Z$ the nonzero ray coefficients from \eqref{eq-shared-coeff} are $1-\theta$, the values $1-f_j$ for $0\le j<r-1$, and the final $1/2$. Their sum is
$r+1-\theta-\sum_{j=0}^{r-1}f_j=r+1-w(a),$
because $\sum_{j=0}^{r-1}\{2^ja/2^r\}=w(a)-a/2^r$. This last identity follows by summing the contributions of the binary digits of $a$. In class $Y$ the corresponding nonzero coefficients are the complements of those in class $X$, at the same $r+1$ positions; their sum is $w(a)$. This proves \eqref{eq-binary-multiplicity}.

The maximum degree occurs at $r=k,w=1$. Completeness, survival and clearing sets follow from Theorem~\ref{thm-cyclic}.
\end{proof}

The atom structure of the binary face depends on the support order. With the order of Proposition~\ref{prop-sort} for the generator, every box point is a sum of minimal-degree atoms.
\begin{proposition}[Binary normal form]\label{prop-binary-normal}
For every $k\ge1$, the matrix of Theorem~\ref{thm-binary}, with the $X$ and $Z$ supports ordered by decreasing and the $Y$ supports by increasing node index, gives three staircase simplices in $\R^{4k+2}$ with a simplex face $F'_k$ of residue group $\Z/2^k\Z$ whose face Hilbert basis consists of its rays and exactly $k$ further elements $H_{2^j}\,$, $0\le j<k$, all of degree three, with $\mathcal C_K(H_{2^j})=2^{k-j}\Z_{>0}\,$. For $a=\sum_j\epsilon_j2^j$,
\[
 H_a=\sum_{j=0}^{k-1}\epsilon_jH_{2^j}\,,\qquad \lambda(H_a)=(\operatorname{wt}(a),\operatorname{wt}(a),\operatorname{wt}(a)),\qquad h^*_{F'_k}(\tau)=(1+\tau^3)^k .
\]
\end{proposition}
\begin{proof}
For the generator $\theta=2^{-k}$ the values on the $X$ and $Z$ columns are $2^{j-k}$ at node $j$, and $1-2^{j-k}$ on the $Y$ columns, so the stated orders are nonincreasing and~\eqref{eq-shared-coeff} gives the coefficient lists $(\frac12,\frac14,\ldots,2^{-k},2^{-k})$ in classes $X,Z$ and the reverse list in class $Y$. The coefficients of the residue $a\theta$ are the fractional parts of $a$ times these, that is, $\{a/2^r\}$ for $1\le r\le k$ and a second $\{a/2^k\}$; because $\{a/2^r\}=\sum_{j<r}\epsilon_j2^{j-r}<1$, the coefficients of the right-hand side add without carry, which proves the identity for $H_a\,$. For $a=2^j$ each class has coefficient sum $\sum_{r=j+1}^k2^{j-r}+2^{j-k}=1$, so $H_{2^j}$ has degree three and is indecomposable by Proposition~\ref{prop-two}, while every other nonzero box point decomposes into the $H_{2^j}$ for its nonzero binary digits. Flooring coefficients gives the face Hilbert basis, the residue orders give the clearing sets, and summing $\tau^{3\operatorname{wt}(a)}$ gives $h^*$. Exposure gives ambient indecomposability.
\end{proof}
The same matrix with supports in node order has $2^k-1$ indecomposable box points (Theorem~\ref{thm-binary}); here, with $q=2^k$, it has only $k=\log_2q$ indecomposable holes, which by Theorem~\ref{thm-atoms} is optimal up to a constant factor in dimension $O(\log q)$.

\subsection{Faces with a single indecomposable hole}\label{sec-onehole}
\begin{theorem}[Indecomposable holes versus residue order]\label{thm-holes-residue}
Let $F$ be a simplex face of three staircase simplices in $\R^d$ with $|G_F|>1$ and exactly $r$ indecomposable holes. Then
\[
 |G_F|\le\binom{\lfloor(d-1)/3\rfloor+r}{r}.
\]
If $r=1$, then $d\ge3|G_F|-2$, $G_F$ is cyclic, the indecomposable hole $h$ has $\mathcal C_K(h)=|G_F|\Z_{>0}\,$, and all its nonzero box coefficients equal $1/|G_F|$. For every $q\ge3$, three staircase simplices in $\R^{3q-2}$ have a simplex facet $F_q$ whose only indecomposable hole is $h_q=\frac1q\sum(\text{rays})$, of degree three, with $\mathcal C_K(h_q)=q\Z_{>0}$ and all proper faces unimodular.
\end{theorem}
\par\medskip\noindent\emph{Idea of the proof.} A box point decomposes into indecomposable box points with no ray summands. A box point has $\lambda_i\le s_i-1$ ($s_i$ rays from constituent $i$) and $\sum_is_i\le d+2$, so degree at most $d-1$; each indecomposable box point has degree at least three (Proposition~\ref{prop-bounds}), so a decomposition has at most $\lfloor(d-1)/3\rfloor$ summands, and box points inject into vectors in $\Z_{\ge0}^r$ with coordinate sum at most $\lfloor(d-1)/3\rfloor$. If $r=1$, every nonzero box point is a multiple of the indecomposable hole, so $G_F$ is cyclic, and $h$, with all nonzero coefficients $1/|G_F|$, uses at least $3|G_F|$ rays. The facets $F_q$ have $3q$ rays whose relations over $\T$ have a common coefficient $\theta$ with $q\theta=0$.

\begin{proof}
The bound and the case $r=1$ are Theorem~\ref{thm-atoms} and Corollary~\ref{cor-onehole-rigid}; the facets $F_q$ are constructed in Theorem~\ref{thm-onehole} below.
\end{proof}

We now construct the facets $F_q$ of Theorem~\ref{thm-holes-residue}; Theorems~\ref{thm-khole} and~\ref{thm-degq} reuse them. Put $a:=\lfloor q/2\rfloor$, $c:=q-a$ and $m:=q-1-a$, and use the $3q-2$ coordinates $P_j\,,Q_j\,,R_j$ ($1\le j<q$) and $g$ with weights
\[
 \begin{gathered}\alpha(P_j)=-1+\mathbf 1_{j>a}+\mathbf 1_{j\ge c}\,,\qquad \alpha(R_j)=-\alpha(P_j)-\mathbf 1_{j=a}\,,\\ \alpha(Q_j)=\mathbf 1_{j=1}-\mathbf 1_{j=c}\,,\qquad \alpha(g)=1,\end{gathered}
\]
which sum to zero. Each growth order consists of an initial block, $q-1$ internal blocks and a final block; internal and final blocks are sorted by nondecreasing weight:
\begin{itemize}
\item[$A$:] initial $\{Q_1\}$; internal blocks $\{P_j\,,R_j\}$ for $j=q-1,\ldots,1$, with $g$ added for $j=a$; final block the remaining $Q_j$;
\item[$B$:] initial $\{g\}$; internal blocks $\{P_j\,,Q_j\}$ for $j=1,\ldots,q-1$, with $R_{j-1}$ added for $j\ge2$; final block $\{R_{q-1}\}$;
\item[$C$:] initial $\{P_j:j>a\}$; internal blocks $\{Q_j\}$ for $j=q-1,\ldots,1$, with $g$ added for $j=c$ and all $R_\ell$ added for $j=1$; final block the remaining $P_j\,$.
\end{itemize}
Let $r_{X,j}\,$, $0\le j\le q-1$, be the lifted prefix of list $X$ ending after its initial block ($j=0$) or after its $j$th internal block.
\begin{theorem}[Optimal faces with one indecomposable hole]\label{thm-onehole}
For every $q\ge3$, these three staircase simplices in $\R^{3q-2}$ have a simplex facet $F_q$ with the $3q$ rays $r_{X,j}\,$, residue group $\Z/q\Z$ and normalized volume $q$, whose face Hilbert basis consists of the rays and $h_q=\frac1q\sum_{X,j}r_{X,j}\,$. The point $h_q$ has multiplicities $(1,1,1)$, degree three, and $\mathcal C_K(h_q)=q\Z_{>0}$; every proper face of $F_q$ is unimodular, and $h^*_{F_q}(\tau)=1+\tau^3+\cdots+\tau^{3(q-1)}$. By Theorem~\ref{thm-atoms} the dimension $3q-2$ is the least possible for a simplex face with a single indecomposable hole of clearing order $q$.
\end{theorem}
\begin{proof}
Put $\ell:=\lambda_A+\lambda_B+m\lambda_C-\alpha^Tx$; the three coefficients are the totals of the initial blocks, because the $P_j$ with $j>a$ have total weight $m$. Every internal block has total weight zero and is either a zero-weight singleton or contains both signs, so its nonempty proper partial sums in the prescribed order stay strictly below its final total; every final block starts with a negative weight and has no positive weight. So, as in Lemma~\ref{lem-growth}, $\ell$ is nonnegative on all lifted vertices and vanishes exactly at the $3q$ prefixes $r_{X,j}$; it is primitive because $\alpha(g)=1$.

For a linear relation over $\R$ or $\T$ with coefficients $u_{X,j}$ on these rays, the multiplicity rows give $\sum_ju_{X,j}=0$. Let $A_j$ and $C_j$ be the sums of the coefficients from the endpoint of descending internal block $j$ onward in classes $A$ and $C$, let $B_j$ be the sum from the endpoint of ascending internal block $j$ onward in class $B$, and put $B_q=0$; entries in an initial block contribute a class total, which vanishes. The rows of $P_j\,$, $Q_j\,$, $R_j$ and $g$ read
\[
 A_j+B_j=0,\qquad C_j+B_j=0,\qquad A_j+C_1+B_{j+1}=0,\qquad A_a+C_c=0 .
\]
Hence $A_j=C_j=j\theta$ and $B_j=-j\theta$ with $\theta=A_1\,$, the last $R$ row gives $q\theta=0$, and the $g$ row repeats it because $a+c=q$. Over $\R$ the relation is zero, so the $3q$ rays are independent and, the cone having dimension $3q+1$, span a simplex facet. Over $\T$, consecutive differences show that every coefficient equals $\theta$ with $q\theta=0$, so the induced coefficient lattice is $\Z^{3q}+\Z\frac1q(1,\ldots,1)$. The box points are the $jh_q$ with $0\le j<q$, each a sum of $j$ copies of $h_q$ without carry; flooring gives the face Hilbert basis, $h_q$ has degree three and is indecomposable, and its clearing set follows from independence and exposure. Deleting any ray kills the torus kernel, so every proper face is unimodular. Summing $\tau^{3j}$ over the box points gives $h^*$.
\end{proof}

Unlike the form of Lemma~\ref{univ-incidence}, the exposing form $\ell$ here need not take only the values zero and one on $V$ (it reaches $2$ for $q=4,5$ and $3$ for $q=6,7$).

\begin{example}[$q=3$]\label{ex-q3}
On the coordinates $P_1\,,P_2\,,Q_1\,,Q_2\,,R_1\,,R_2\,,g$ with weights
\[\alpha=(-1,1,1,-1,0,-1,1),\]
the construction of Theorem~\ref{thm-onehole} gives the growth orders
\[
\begin{gathered}
 (Q_1\mid R_2\,,P_2\mid P_1\,,R_1\,,g\mid Q_2),\qquad(g\mid P_1\,,Q_1\mid Q_2\,,R_1\,,P_2\mid R_2),\\
 (P_2\mid Q_2\,,g\mid R_2\,,R_1\,,Q_1\mid P_1).
\end{gathered}
\]
Each consists of an initial block of weight $1$, two blocks of weight $0$ and a final block. The form $t-\alpha^Tx$ takes only the values $0$ and $1$ on $V$ and vanishes exactly at the nine prefixes ending at the bars. These are independent, so they are the vertices of a simplex facet $F_3$ of $\mathcal D$, whose cone $K_{F_3}$ is a facet of the $10$-dimensional cone, and their sum is $3h_3$ with $h_3=(1,2,2,1,1,1,2;\,1,1,1)$, of degree three. Neither $h_3$ nor $2h_3$ splits, $3h_3$ does, and $h_3$ is the only indecomposable hole of the whole hull; here $d=7=3q-2$.
\end{example}
 The facets $F_q$ are faces failing IDP all of whose proper faces have IDP, in unbounded dimension, in contrast with the five-dimensional minimal non-IDP faces of word disjunctions (Theorem~\ref{thm-words-minface}); this does not exclude other, smaller non-IDP faces of the same hull.

\subsection{Faces with few indecomposable holes}\label{sec-onehole-any}
Subsection~\ref{sec-onehole} concerned simplex faces whose only indecomposable hole has a clearing order. We first drop the simplex assumption, and then consider simplex faces with as few indecomposable holes as the residue group permits, namely $\mu(G)$; this determines the least dimension exactly and leads to a hardness result and to a trade-off between the residue order and the number of indecomposable holes.
\begin{theorem}[One indecomposable hole on an arbitrary face]\label{thm-onehole-any}
Suppose the cone $K_F$ of a face $F$, in its induced lattice, has exactly one indecomposable hole $h$, with least clearing integer $q$. Then
\[
 \operatorname{rank}K_F\ge(q-1)\deg h+3\ge3q .
\]
In particular $d\ge3q-2$ if $F$ is proper, and $d\ge3q-3$ if $F=\mathcal D$.
\end{theorem}
\begin{proof}
Let $h=\sum_v\alpha_vv$ be a nonnegative real expression on the rays of $F$. If $\alpha_v\ge1/(q-1)$ for some $v$, then $(q-1)h-v$ is a lattice point of $K_F\,$, so it equals $s+rh$ with $s$ in the vertex monoid and $r\in\Z_{\ge0}\,$. Then $(q-1-r)h=s+v$ splits, and comparing degrees gives $1\le q-1-r<q$, contradicting the minimality of $q$. So all coefficients are below $1/(q-1)$. In an expression on independent rays, the coefficients from constituent $i$ sum to $\lambda_i(h)\ge1$, so at least $(q-1)\lambda_i(h)+1$ rays come from constituent $i$. Summing and using $\deg h\ge3$ gives the bound; a proper face has rank at most $d+2$ and $K$ has rank $d+3$.
\end{proof}

\begin{lemma}[Rigidity with the fewest indecomposable holes]\label{lem-rigid}
Let $L\supseteq\Z^n$ be a lattice, $G=L/\Z^n$ nontrivial with $r=\mu(G)$ generators, and call the nonzero elements of $L\cap[0,1)^n$ that are not sums of two such elements the indecomposable holes. There are at least $r$ indecomposable holes. If there are exactly $r$, say $h_1\,,\ldots,h_r\,$, then $h_i=\frac1{q_i}\mathbf 1_{S_i}$ for pairwise disjoint nonempty sets $S_i$ and integers $q_i\ge2$, the box points are exactly the sums $\sum_ia_ih_i$ with $0\le a_i<q_i\,$, and $G\cong\bigoplus_i\Z/q_i\Z$.
\end{lemma}
\begin{proof}
A box point is a sum of indecomposable holes, because a unit vector cannot be subtracted from it within $\R^n_{\ge0}$; so the residues of the indecomposable holes generate $G$, and there are at least $r$. Suppose there are exactly $r$, with matrix $H=(h_1\cdots h_r)$. Choose a prime $p$ with $\dim_{\F_p}G/pG=r$. The residues of the indecomposable holes form a basis of $G/pG$, so $L_0:=\{a\in\Z^r:Ha\in\Z^n\}\subseteq p\Z^r$. The columns of $H$ are independent, because a primitive integral relation would lie in $L_0\,$. Let $A:=\{a\in\Z^r_{\ge0}:Ha\in[0,1)^n\}$, a down-set; $a\mapsto[Ha]$ maps $A$ onto $G$ by the first remark, and injectively, because two box points with the same residue coincide.

Let $a\notin A$ be minimal, so $a-\mathbf e_i\in A$ for $i\in\supp a$, and let $b\in A$ have the same residue. If $a_i\,,b_i>0$, then $a-\mathbf e_i$ and $b-\mathbf e_i$ are elements of $A$ with the same residue, so $a=b\in A$; hence $a,b$ have disjoint supports. As $a-b\in L_0\subseteq p\Z^r$, both lie in $p\Z^r$. Because $Ha=H(a-\mathbf e_i)+h_i\in[0,2)^n$ and $Hb\in[0,1)^n$, the integral vector $z:=Ha-Hb$ lies in $\{0,1\}^n$. Then $z/p=H\bigl((a-b)/p\bigr)$ is a box point, equal to $Hc$ with $c\in A$, and independence gives $a-b=pc\ge0$; disjointness forces $b=0$. So $Ha\in\{0,1\}^n$ for every minimal $a\notin A$.

For each $i$ let $q_i$ be least with $q_i\mathbf e_i\notin A$. Then $q_ih_i\in\{0,1\}^n$, so $h_i=\frac1{q_i}\mathbf 1_{S_i}\,$, and $q_i\ge2$ because $h_i\ne0$ lies in the box. If $S_i\cap S_j\ne\varnothing$ for $i\ne j$, then $(q_i-1)\mathbf e_i+(q_j-1)\mathbf e_j\notin A$, and a minimal $a\notin A$ below it uses both $i$ and $j$. By independence $S_i\ne S_j\,$, and at a coordinate of $S_i\triangle S_j$ the vector $Ha$ takes the value $a_i/q_i$ or $a_j/q_j\,$, strictly between $0$ and $1$, a contradiction. With disjoint supports, $A=\prod_i\{0,\ldots,q_i-1\}$, so $|G|=\prod_iq_i\,$, and the residues $[h_i]$ of orders $q_i$ generate $G$, whence $G\cong\bigoplus_i\Z/q_i\Z$.
\end{proof}

For a nontrivial finite abelian group $G$ with $r=\mu(G)$, put
\[
 \Sigma(G):=\min\Bigl\{\sum_{j=1}^rq_j:\ G\cong\bigoplus_{j=1}^r\Z/q_j\Z,\ q_j\ge2\Bigr\},
\]
and let $d_{\min}(G)$ be the least spatial dimension of three staircase simplices with a simplex face whose residue group is $G$ and which has exactly $\mu(G)$ indecomposable holes.

\begin{figure}[t]
\centering
\begin{tikzpicture}[x=1.0cm,y=0.9cm,font=\small]
 \node[circle,fill,inner sep=1.8pt,label={[font=\scriptsize]left:$P\ (+1)$}] (P) at (0,0) {};
 \node[circle,fill,inner sep=1.8pt,label={[font=\scriptsize]right:$N\ (-1)$}] (N) at (6.8,0) {};
\node[circle,fill,inner sep=1.6pt,label={[font=\scriptsize]above:$x_0$}] (x0) at (2.2,1.70) {};
\node[circle,fill,inner sep=1.6pt,label={[font=\scriptsize]above:$y_0$}] (y0) at (4.6,1.70) {};
\node[circle,fill,inner sep=1.6pt,label={[font=\scriptsize]above:$x_1$}] (x1) at (2.2,0.00) {};
\node[circle,fill,inner sep=1.6pt,label={[font=\scriptsize]above:$y_1$}] (y1) at (4.6,0.00) {};
\node[circle,fill,inner sep=1.6pt,label={[font=\scriptsize]above:$x_2$}] (x2) at (2.2,-1.70) {};
\node[circle,fill,inner sep=1.6pt,label={[font=\scriptsize]above:$y_2$}] (y2) at (4.6,-1.70) {};
\draw[thick,blue!70!black] (P)--node[font=\scriptsize,fill=white,inner sep=1pt]{$\theta_0$}(x0);
\draw[thick,dotted,green!50!black] (x0)--node[font=\scriptsize,fill=white,inner sep=1pt]{$-\theta_0$}(y0);
\draw[thick,dashed,red!70!black] (y0)--node[font=\scriptsize,fill=white,inner sep=1pt]{$\theta_0$}(N);
\draw[thick,dashed,red!70!black] (P)--node[font=\scriptsize,fill=white,inner sep=1pt]{$\theta_1$}(x1);
\draw[thick,blue!70!black] (x1)--node[font=\scriptsize,fill=white,inner sep=1pt]{$-\theta_1$}(y1);
\draw[thick,dotted,green!50!black] (y1)--node[font=\scriptsize,fill=white,inner sep=1pt]{$\theta_1$}(N);
\draw[thick,dotted,green!50!black] (P)--node[font=\scriptsize,fill=white,inner sep=1pt]{$\theta_2$}(x2);
\draw[thick,dashed,red!70!black] (x2)--node[font=\scriptsize,fill=white,inner sep=1pt]{$-\theta_2$}(y2);
\draw[thick,blue!70!black] (y2)--node[font=\scriptsize,fill=white,inner sep=1pt]{$\theta_2$}(N);
 \node[font=\scriptsize,gray] at (1.4,1.55) {$+z_0$};
 \node[font=\scriptsize,gray] at (5.4,1.55) {$+z_0$};
 \node[font=\scriptsize,gray] at (1.1,-0.32) {$+z_1$};
 \node[font=\scriptsize,gray] at (5.7,-0.32) {$+z_1$};
 \begin{scope}[xshift=8.6cm,yshift=0.6cm]
  \draw[thick,blue!70!black] (0,0.6)--(0.7,0.6); \node[anchor=west,font=\scriptsize] at (0.8,0.6) {color $0$};
  \draw[thick,dashed,red!70!black] (0,0)--(0.7,0); \node[anchor=west,font=\scriptsize] at (0.8,0) {color $1$};
  \draw[thick,dotted,green!50!black] (0,-0.6)--(0.7,-0.6); \node[anchor=west,font=\scriptsize] at (0.8,-0.6) {color $2$};
  \node[anchor=west,font=\scriptsize] at (0,-1.3) {$x_j\,(-1)$, $y_j\,(+1)$};
 \end{scope}
\end{tikzpicture}
\caption{The order-two block for $(\Z/2\Z)^2$ in the proof of Theorem~\ref{thm-khole}. Vertices are coordinates, with their weights; each edge is an internal block of the indicated color, so each color class has disjoint supports, and each color misses one positive and one negative vertex. The edges of path $j$ carry the values $\theta_j\,,-\theta_j\,,\theta_j$ in a torus relation. The zero-weight coordinate $z_0$ lies in the first and last blocks of path $0$, and $z_1$ in those of path $1$. The rows of $P$ and $N$ give $\theta_0+\theta_1+\theta_2=0$, and those of $z_0\,,z_1$ give $2\theta_0=2\theta_1=0$, so the residue group is $(\Z/2\Z)^2$.}
\label{fig-ordertwo}
\end{figure}
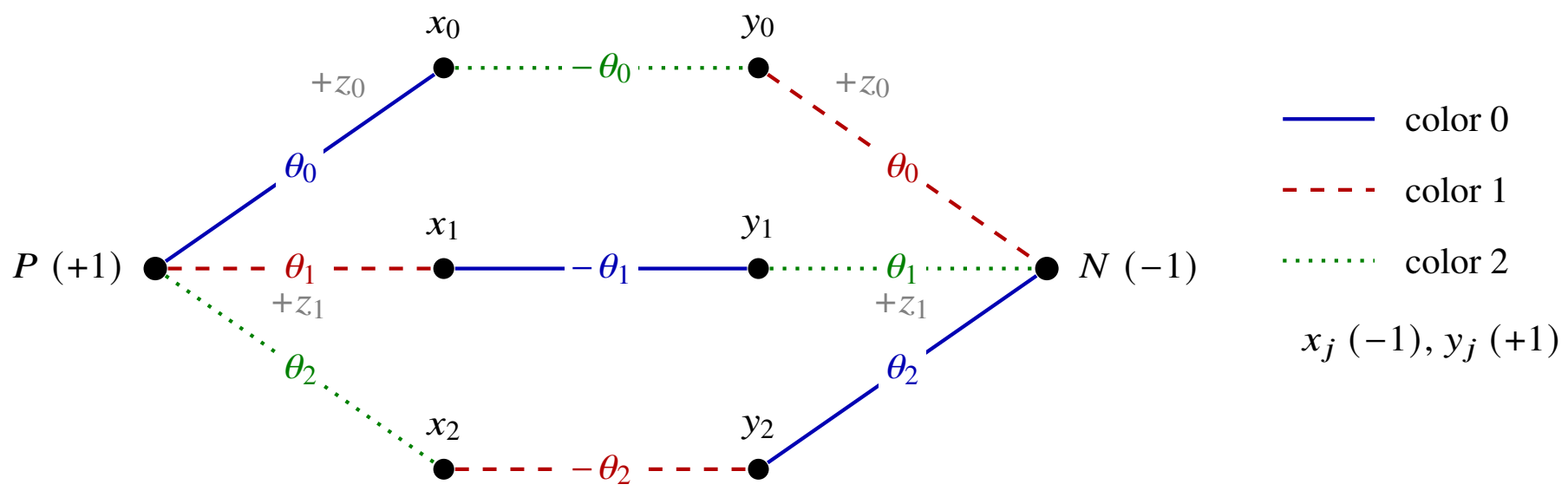

\begin{theorem}[Minimum dimension with the fewest indecomposable holes]\label{thm-khole}
Let $G$ be a nontrivial finite abelian group and $r=\mu(G)$. A simplex face of three staircase simplices in $\R^d$ with residue group $G$ has at least $r$ indecomposable holes. If it has exactly $r$, then they have pairwise disjoint ray supports $S_j\,$, $h_j=\frac1{q_j}\sum_{v\in S_j}v$, $G\cong\bigoplus_j\Z/q_j\Z$, $\mathcal C_K(h_j)=q_j\Z_{>0}\,$, and $d\ge\sum_jq_j\deg h_j-2\ge3\Sigma(G)-2$. For every list $q_1\,,\ldots,q_s\ge2$ other than the single entry $2$, there is a simplex facet in dimension $3\sum_jq_j-2$ with residue group $\bigoplus_j\Z/q_j\Z$ and exactly $s$ indecomposable holes, each of degree three. Consequently
\[
 d_{\min}(G)=3\Sigma(G)-2\quad(G\not\cong\Z/2\Z),\qquad d_{\min}(\Z/2\Z)=6 .
\]
\end{theorem}
\begin{proof}
\emph{Rigidity and the lower bound.} In ray coordinates the face lattice is a lattice $L\supseteq\Z^n$ with $L/\Z^n=G$, and the indecomposable holes are the indecomposable holes of Lemma~\ref{lem-rigid}. This gives the count, the disjoint supports, the coefficients $1/q_j$ and the group. An indecomposable hole has positive multiplicity in every constituent, because a lattice point using at most two simplices splits by Proposition~\ref{prop-two}; so $\deg h_j\ge3$ and $|S_j|=q_j\deg h_j\,$. The face is proper, because $\mathcal D$ has $3(d+1)>d+3$ vertices, so it has at most $d+2$ rays, which gives $d\ge\sum_jq_j\deg h_j-2$; and $(q_j)$ is a decomposition of $G$ into $r$ cyclic factors. Independence and exposure give the clearing sets (Lemma~\ref{lem-face}). For $\Z/2\Z$, Proposition~\ref{prop-dim5} excludes $d\le5$.

\emph{Reading a construction.} For spatial weights $\alpha$ with total zero, let $M_X$ be the maximal prefix weight in class $X$; the form $\sum_XM_X\lambda_X-\alpha^Tx$ is nonnegative on $V$ and vanishes exactly at the maximizing prefixes. Call the segments between consecutive maximizing prefixes the internal blocks. In a relation among these rays over $\R$ or $\T$ the class totals vanish, and the coordinate rows are linear equations in the tail sums over the internal blocks, which determine the coefficients by differences. A block of total weight zero that contains both signs, read with negatives first, then zeros, then positives, has no intermediate maximizing prefix; nor has a zero-weight singleton.

\emph{Cyclic blocks, $q\ge3$.} Use the facet $F_q$ of Theorem~\ref{thm-onehole} on $3q-2$ coordinates. Its maxima are $(1,1,m)$ with $m=q-1-\lfloor q/2\rfloor$, attained exactly at the $q$ block endpoints $r_{X,j}$ per class, and in a relation over $\T$ all $3q$ endpoint coefficients equal one $\theta$ with $q\theta=0$.

\emph{Order-two blocks.} For positive integers $m_X$ use the six coordinates $u^\pm_X$ of weights $\pm m_X\,$, $X\in\{A,B,C\}$. Class $X$ reads $u^+_X\,$, then one internal block of the four coordinates of the other two labels, negatives first, then $u^-_X\,$. The maxima are $m_X$ with two endpoints per class, and either row of label $X$ says that the internal tail sums of the other two classes add to zero; so all six endpoint coefficients equal one $\theta$ with $2\theta=0$, and the difference of the two rows of label $X$ is the class-$X$ total of the block. For two factors of order two (Figure~\ref{fig-ordertwo}) take three paths $P\,x_j\,y_j\,N$ ($j=0,1,2$) with weights $+1,-1,+1,-1$ along each path, so eight coordinates. Color the three edges of path $j$ by $j,j+2,j+1$ modulo $3$; each edge is an internal block of its color containing its two endpoints. Add a zero-weight coordinate $z_0$ to the first and last edge of path $0$ and $z_1$ to the first and last edge of path $1$. Each color misses one positive and one negative vertex; list the missing zero coordinate and the missing positive vertex first, then the blocks of paths $0,2,1$, then the missing negative vertex. The maxima are $(1,1,1)$. The internal rows of path $j$ force its block values $\theta_j\,,-\theta_j\,,\theta_j\,$, the rows of $P,N$ give $\theta_0+\theta_1+\theta_2=0$, and those of $z_0\,,z_1$ give $2\theta_0=2\theta_1=0$. So the real kernel is zero, and the twelve endpoint coefficients are $\theta_0$ on six rays and $\theta_1$ on the other six: two disjoint supports of order two.

\emph{Concatenation.} If some $q_j\ge3$, order the cyclic blocks by nonincreasing $m_t\,$, and prepend the order-two blocks with maxima $(1,1,m_1)$ on disjoint coordinates. Between consecutive cyclic blocks $t,t+1$ add $s^\pm_t$ of weights $\pm(m_t-m_{t+1})$: in $B$ put $s^-_t\,,s^+_t$ at the seam, in $C$ put $s^+_t$ at the seam and every $s^-_t$ after all blocks, and in $A$ put every $s^-_t$ and then every $s^+_t$ after all blocks. Then the offset of cyclic block $t$ in $C$ is $m_1-m_t\,$, every local maximizing prefix is globally maximizing, and the seam and terminal entries stay below the maxima. In a relation, let $T_{t,X}$ be the total coefficient of the class-$X$ rays of block $t$. Each order-two block has $T_{t,X}=0$ by the row differences above. The row of $s^-_t$ says that the class-$B$ total after seam $t$ vanishes, and then the row of $s^+_t$ says the same for class $C$; with the class totals, $T_{t,B}=T_{t,C}=0$ for every cyclic block. The weighted sum of the rows of cyclic block $t$ gives $T_{t,A}+T_{t,B}+m_tT_{t,C}=0$, because earlier prefixes miss the block and later ones contain it with total weight zero; so $T_{t,A}=0$. The remaining equations are the local systems, and conversely local kernel vectors extend. If all $q_j=2$ and $s\ge2$, concatenate $s-2$ six-coordinate blocks and one ten-coordinate block, all with maxima $(1,1,1)$; the row differences and the class totals again decouple the blocks. In both cases there are $3\sum_jq_j$ independent rays in dimension $3\sum_jq_j-2$, so the face is a simplex facet; its box coefficients are constant on the disjoint supports, with independent parameters, so the indecomposable holes are the $s$ block generators, each with $q_j$ rays and coefficient $1/q_j$ per class, hence of degree three. Taking a decomposition attaining $\Sigma(G)$ gives $d_{\min}(G)\le3\Sigma(G)-2$, except for $\Z/2\Z$, where the facet of the introductory example in $\R^6$, exposed by $x_1+x_2-x_4-x_5+x_6\,$, has residue group $\Z/2\Z$ and one indecomposable hole.
\end{proof}

\begin{corollary}[Box points with the fewest indecomposable holes]\label{cor-khole-box}
If a simplex face $F$ has exactly $\mu(G_F)$ indecomposable holes $h_i\,$, of clearing orders $q_i\,$, then its box points are the sums $\sum_ia_ih_i$ with $0\le a_i<q_i\,$, so $h^*_F(\tau)=\prod_i\bigl(1+\tau^{\deg h_i}+\cdots+\tau^{(q_i-1)\deg h_i}\bigr)$, and the minimal faces of $F$ failing IDP are exactly the simplices spanned by the supports $S_i\,$.
\end{corollary}
\begin{proof}
The box description is Lemma~\ref{lem-rigid}. A face of $F$ contains the box point $\sum_ia_ih_i$ exactly when it contains $S_i$ for every $i$ with $a_i>0$, and it fails IDP exactly when it contains a nonzero box point.
\end{proof}

\begin{theorem}[Computing the minimum dimension]\label{thm-dmin-hard}
Given the primary cyclic factors of a finite abelian group $G$ and an integer $D_0$ in binary, deciding whether $d_{\min}(G)\le D_0$ is NP-complete, already for $G=(\Z/3\Z)^{m-1}\oplus\Z/3N\Z$ with $\gcd(N,3)=1$ and the factorization of $N$ given. Unless $\mathrm P=\mathrm{NP}$, $d_{\min}(G)$ cannot be approximated in polynomial time within any factor polynomial in the input length.
\end{theorem}
\begin{proof}
A decomposition of $G$ into $\mu(G)$ cyclic factors is the same as an assignment of the primary cyclic factors to $\mu(G)$ slots, each nonempty and with at most one power of each prime; the slot products are the $q_j\,$. Such an assignment is a polynomial certificate for Theorem~\ref{thm-khole}, and $\Z/2\Z$ is handled directly.

Reduce from $3$-PARTITION, which is strongly NP-complete~\cite{GJ}: integers $B/4<a_i<B/2$, $1\le i\le3m$, polynomial in $m$, with $\sum_ia_i=mB$. Choose distinct primes $p_i\ge5$ of polynomial size, put $P_0=\max_ip_i\,$, take $L$ with $2^L>mP_0^3\,$, and let $p_i^{e_i}$ be the least power with $p_i^{e_i}\ge2^{La_i}\,$, so $p_i^{e_i}<p_i2^{La_i}$. Put $N=\prod_ip_i^{e_i}$. The group has $3$-rank $m$ and every other rank one, so $\mu(G)=m$; each slot receives one factor $3$, and $q_j=3\prod_{i\in I_j}p_i^{e_i}$ for a partition $I_1\,,\ldots,I_m$ of the $3m$ indices. A $3$-partition gives $\sum_jq_j<3mP_0^32^{LB}$. Otherwise every partition has a slot of weight sum at least $B+1$, because a partition with all sums $B$ would have three items per slot; that slot has $q_j\ge3\cdot2^{L(B+1)}>3mP_0^32^{LB}$. With $D_0=9mP_0^32^{LB}-2$, Theorem~\ref{thm-khole} decides the instance, and all numbers have polynomial bit length. Taking $L=m^2+1+\lceil\log_2(mP_0^3)\rceil$ makes the ratio between the two cases exceed $2^{m^2}$, which is larger than any polynomial in the input length for large $m$.
\end{proof}
For instance, $\Z/2\Z\oplus\Z/30\Z$ has $B=6+10$, from the factors $\Z/6\Z\oplus\Z/10\Z$, so $d_{\min}=46$, while its invariant factors would give $94$.

\begin{corollary}[Residue size versus number of indecomposable holes]\label{cor-tradeoff}
For $d\ge7$ and $r\ge1$, let $M(d,r)$ be the largest order of $G_F$ over simplex faces $F$ with at most $r$ indecomposable holes of three staircase simplices in $\R^d$. There are absolute constants $c,C>0$ with
\[
 c\,r\log(1+d/r)\le\log M(d,r)\le C\,r\log(1+d/r).
\]
\end{corollary}
\begin{proof}
\emph{Upper bound.} Put $m=\lfloor(d-1)/3\rfloor$. By Theorem~\ref{thm-atoms}, a simplex face with $\nu\le r$ indecomposable holes has $|G_F|\le\binom{m+\nu}{\nu}\le\binom{m+r}{r}\le\bigl(e(m+r)/r\bigr)^r$. If $r\le d$, then $\log(1+d/r)\ge\log2$, so $\log M(d,r)\le r\bigl(1+\log(1+d/r)\bigr)\le(1+1/\log2)\,r\log(1+d/r)$. If $r>d$, Proposition~\ref{prop-det} gives $|G_F|\le3^{(d-1)/3}$, so $\log M(d,r)\le d\log3$, while $\log(1+x)\ge x/2$ for $0\le x\le1$ gives $r\log(1+d/r)\ge d/2$.

\emph{Lower bound.} Put $s=\min\{r,\lfloor(d+2)/9\rfloor\}\ge1$ and $q=\lfloor(d+2)/(3s)\rfloor\ge3$. Theorem~\ref{thm-khole} gives a simplex facet with group $(\Z/q\Z)^s$ and $s$ indecomposable holes in dimension $3sq-2\le d$. Appending $d-(3sq-2)$ coordinates at the end of all three growth orders embeds that cone as the face where these coordinates vanish, so it is a simplex face in $\R^d$ with the same group and indecomposable holes (Lemma~\ref{lem-face}). If $s=r$, then $x:=(d+2)/(3r)\ge3$, $q\ge\max\{3,2x/3\}$, so $\log q$ is at least an absolute constant times $\log x$, and $1+d/r\le3x+1\le x^3$, so $\log M(d,r)\ge r\log q\ge c\,r\log(1+d/r)$. If $s<r$, then $s\ge(d+2)/18$, so $\log M(d,r)\ge s\log3\ge(d+2)\log3/18$, while $r\log(1+d/r)\le d$.
\end{proof}

\subsection{The largest degree at a fixed clearing order}\label{sec-degq}
The single-hole facets of Theorem~\ref{thm-onehole} have indecomposable holes of degree three. We now ask how large the degree of an indecomposable hole of a given clearing order can be.
On a simplex face $F$, the clearing set of a lattice point $h$ is $q\Z_{>0}$ with $q$ the order of $[h]$ in $G_F\,$, by Lemma~\ref{lem-face} and independence of the rays.
\begin{theorem}[Sharp degree at a fixed clearing order]\label{thm-degq}
Let $h$ be an indecomposable hole on a simplex face $F$ with $s$ rays of three staircase simplices in $\R^d$, of clearing order $q\ge2$. Then
\[
 \deg h\le\Bigl\lfloor\frac{(q-1)s-q+2}{q}\Bigr\rfloor\le\Bigl\lfloor\frac{(q-1)(d+1)+1}{q}\Bigr\rfloor.
\]
The second bound is attained by an indecomposable hole of clearing order $q$ on a simplex facet for every $q\ge3$ and $d\ge3q-1$, and for $q=2$ and $d\ge7$.
\end{theorem}
\begin{proof}[Upper bound]
Write $h=\sum_j\frac{p_j}{q}r_j$ over the rays $r_j$ of $F$, with integers $0\le p_j<q$; the coefficients lie in $\frac1q\Z$ because $q[h]=0$, and below one because a coefficient at least one would let us subtract a ray. Every ray has degree one, so $q\deg h=\sum_jp_j\,$. For $2\le k\le q-1$ let $a_k$ and $b_k$ be the box points of $kh$ and $(q+1-k)h$; both are nonzero, because $1<k,q+1-k<q$. Their $j$th coefficients are $(kp_j\bmod q)/q$ and $((q+1-k)p_j\bmod q)/q$, whose sum is $p_j/q$ or $p_j/q+1$, the latter exactly when $kp_j\bmod q>p_j\,$. If no index carries, then $h=a_k+b_k$ decomposes; so for each such $k$ some $j$ has $kp_j\bmod q>p_j\,$. For fixed $j$, the residues $kp_j\bmod q$, $0\le k<q$, run $\gcd(p_j\,,q)$ times through the multiples of $\gcd(p_j\,,q)$ (with $\gcd(0,q)=q$), and $k=0,1$ do not exceed $p_j$; so at most $q-\gcd(p_j\,,q)-p_j\le q-1-p_j$ values $k\in[2,q-1]$ carry at $j$. Covering all $q-2$ values gives $\sum_j(q-1-p_j)\ge q-2$, which is the first bound. A simplex face is proper, because $\mathcal D$ has $3(d+1)>d+3$ vertices, so $s\le d+2$, which gives the second.
\end{proof}

\par\medskip\noindent\emph{Cloning a ray.} Let $F$ be a simplex facet with exposing form $\ell=c^T\lambda-\alpha^Tx$ and a maximizing prefix $r$ of constituent $X$ among its rays. Insert new coordinates $y_1\,,\ldots,y_m$ with $\alpha(y_t)=0$ in all three growth orders: in $X$ consecutively, immediately after the last coordinate of $r$; in each other constituent anywhere at a \emph{slot}, a position after a prefix with $\ell>0$. Then $\ell\ge0$ on the new lifted vertices, and it vanishes exactly on the old rays other than $r$ and on the $m+1$ \emph{clones} $c_0=r,c_1\,,\ldots,c_m\,$, where $c_t$ also contains $y_1\,,\ldots,y_t\,$. Every old coordinate lies in all or none of the clones, so for a relation over $\R$ or $\T$ the old coordinate and multiplicity rows involve the clone coefficients only through their sum and reproduce the old relations, with that sum in the place of $r$. The row of $y_t$ reads $T_t+\sigma_t=0$, where $T_t$ is the sum of the clone coefficients from $c_t$ on and $\sigma_t$ is the sum of the coefficients of the old rays after the clones in $X$ and after the slots of $y_t$ in the other constituents. So the real kernel stays zero, the $d+2$ rays span a simplex facet, its residue group is the old one, and each clone coefficient of a residue is determined by the old residue and the slots. If, for an old residue $u$, the slots of each $y_t$ can be chosen to give every value of $\sigma_t$ in a set containing $-T_t$ for all $T_t$ we want, then any clone coefficients in $[0,1)$ summing to $u_r$ modulo one can be prescribed.

\begin{proof}[Attainment for $q\ge3$]
Let $m=d-(3q-2)\ge1$ and clone the ray $r_{A,0}$ of the facet $F_q$ of Theorem~\ref{thm-onehole}. Put all $y_t$ at the end of $C$, after its final block, which starts with a negative weight; no ray of $C$ contains them. In $B$, a slot inside the internal block $j$, after its negative entries, or after the final block lies before exactly $n$ rays of $B$, for every $n\in\{0,\ldots,q-1\}$. In the residue of $F_q$ with all coefficients $\theta=\frac{q-1}q$, we get $\sigma_t=(q-1+n)\theta$, which takes every value in $\frac1q\Z/\Z$. Prescribe the clone numerators: one equal to $1$, one equal to $q-1-\kappa$ and the other $m-1$ equal to $q-1$, where $\kappa\in\{0,\ldots,q-1\}$ and $\kappa\equiv1+(m-1)(q-1)\pmod q$; their sum is $\equiv q-1$, as required. Let $h$ be the box point of this residue. The coefficient at the clone with numerator $1$ is a character of the cyclic residue group sending a generator to $\frac1q$, so it is positive on every nonzero residue. A decomposition of $h$ into two nonzero lattice points of the facet would consist of two nonzero box points, with coefficients there summing to at least $\frac2q$; so $h$ is indecomposable, and by Lemma~\ref{lem-face} it is an indecomposable hole of $K$ with clearing set $q\Z_{>0}\,$. Its $3q-1$ old rays and the clones give
\[
 q\deg h=(q-1)(3q-1)+1+m(q-1)-\kappa=(q-1)(d+1)+1-\kappa,
\]
the largest multiple of $q$ at most $(q-1)(d+1)+1$.
\end{proof}

\begin{proof}[Attainment for $q=2$]
Let $\alpha:=(-1,-1,0,1,1,-1,0,0)$ and $\ell:=-\alpha^Tx$, truncated to the dimension at hand. Begin with the triple $(1,2,4,5,6,3)$, $(1,3,4,6,5,2)$, $(2,3,5,6,4,1)$ in $\R^6$. Each order has prefix maximum zero, attained at lengths $0,4$, $0,3,5$ and $0,3,5$, so $\ell$ exposes a facet with these $8$ rays. Write $a_p\,,b_p\,,c_p$ for the coefficients of a relation on the rays of length $p$ in the three orders. The rows of coordinates $1$ and $4$ give $c_5=0$, the row of $6$ gives $b_5=0$, the rows of $1,2,3$ give $b_3=c_3=-a_4$ and $2a_4=0$, and the multiplicity rows give $a_0=-a_4\,$, $b_0=-b_3\,$, $c_0=-c_3$; the row of $5$ then holds. So the real kernel is zero and the torus kernel has order two; its nonzero element has coefficient $\frac12$ on six rays and $0$ on the two length-five rays. Clone the length-five ray of the second order by coordinate $7$, with slots after coordinate $1$ in the first order and at the end of the third; then clone the length-five ray of the third order by coordinate $8$, with slots after coordinate $1$ in the first order and at the end of the second. This gives the facets of
\[
 \begin{gathered}(1,7,2,4,5,6,3),\ (1,3,4,6,5,7,2),\ (2,3,5,6,4,1,7)\quad\text{and}\\(1,7,8,2,4,5,6,3),\ (1,3,4,6,5,7,2,8),\ (2,3,5,6,4,8,1,7),\end{gathered}
\]
with $9$ and $10$ rays. The row of the new coordinate forces both clone coefficients of the nonzero residue to be $\frac12$, so its box point has coefficients $\frac12$ on all rays except one in $\R^7$, and on all ten rays in $\R^8$. Being the only nonzero box point, it is indecomposable, has clearing set $2\Z_{>0}\,$, and has degree $4$, respectively $5$. For $d\ge9$, clone the length-six ray of the second order of the $\R^8$ triple. Place each $y_t$ at the end of the third order, where no ray contains it, and in the first order either after coordinate $1$, below that order's ray of length six with coefficient $\frac12$, or at the end; this gives $\sigma_t=\frac12$ or $0$. Hence any clone coefficients in $\{0,\frac12\}$ with sum $\frac12$ modulo one can be prescribed. Taking the largest odd number of halves among the $d-7$ clones gives degree $\lfloor(d+2)/2\rfloor$.
\end{proof}

\subsection{The two/three threshold}\label{sec-threshold}
The binary family gives exponentially many indecomposable holes and exponentially large clearing integers. Together with the upper bounds of Section~\ref{sec-foundations}, it determines the growth rate for every fixed number of constituents.
\begin{theorem}[The two/three threshold]\label{thm-threshold}
For $n\le1$, the lifted vertex matrix is totally unimodular, $H_n(d)=0$, and $\Delta_n(d)=1$. For every fixed $n\ge2$,
\[
 H_n(d)=2^{\Theta_n(d)},\qquad \Delta_n(d)=2^{\Theta_n(d)}.
\]
For $d\ge6$, $H_n(d)\ge2^{\lfloor(d-2)/4\rfloor}-1$ and $\Delta_n(d)\ge2^{\lfloor(d-2)/4\rfloor}$. These lower bounds for $\Delta_n$ are attained by Hilbert basis elements of degree three, and every indecomposable hole has degree between $3$ and $d$.
\end{theorem}
\begin{proof}
Proposition~\ref{prop-two} proves the assertion for at most two inputs. For fixed $n\ge2$, Proposition~\ref{prop-bounds} gives exponential upper bounds. Set $k=\lfloor(d-2)/4\rfloor$ and use Theorem~\ref{thm-binary}. If necessary, append coordinates at the end of all three growth orders and expose the face on which these coordinates vanish. This embeds the old cone with its induced lattice. To increase the number of constituents, append arbitrary staircase simplices and set their multiplicities to zero. Exposure preserves indecomposability and every clearing set. Hence $H_n(d)\ge2^k-1$ and $\Delta_n(d)\ge2^k$. Degree three is attained by Proposition~\ref{prop-sort} and Corollary~\ref{cor-sort}. The degree range is Proposition~\ref{prop-bounds}.
\end{proof}

\section{Arithmetic universality}\label{sec-universality}
This section proves result~\hyperlink{res-A}{(A)}, arithmetic universality (Theorem~\ref{thm-main}). A short scalar presentation encodes the group; an encoder and a connector make the supports of any two nonzero solutions intersect, which the incidence construction of Section~\ref{sec-realization} turns into indecomposable box points. Subsection~\ref{sec-universality-overview} outlines the argument and Subsections~\ref{sec-support-circuit}--\ref{sec-universality-thm} carry it out. Subsections~\ref{sec-minfaces} and~\ref{sec-simul} give two refinements: minimal non-IDP facets in logarithmic dimension, and one triple for all small groups.

\subsection{Overview}\label{sec-universality-overview}
The proof of Theorem~\ref{thm-main} has four steps, carried out in the following subsections.

\emph{Encode the group.} Write $G=\bigoplus_i\Z/q_i\Z$ and present each factor by an addition chain for $q_i\,$. Internal additions define new scalars; the terminal addition is constrained to zero. The $M=O(\log|G|)$ initial scalars have solutions exactly $G$ over $\T$ and only zero over $\R$.

\emph{Force overlap without changing the group.} Apply a linear-size additive encoder to these scalars (Lemma~\ref{univ-encoder}). Every nonzero input has a fixed positive fraction of nonzero outputs. Give each output position $u$ a fresh node, with value $c_u\,$, and add a node $c_u+c_v$ for each edge $uv$ of a sparse graph that joins any two disjoint vertex sets of linear size (Lemma~\ref{univ-connector}). Call a node \emph{active} for a solution if its value is nonzero. If the output supports of two nonzero solutions meet, they share an active node; if not, a connector edge between them has nonzero sum for both, because exactly one endpoint is active. The new scalars are determined by the old ones, so neither kernel changes. There are $Q=O(M)$ nodes and equations (Proposition~\ref{univ-presentation}).

\begin{figure}[t]
\centering
\begin{tikzpicture}[x=0.8cm,y=0.8cm,scale=1.764]
\draw[thick] (-1.126,0.650)--(0.000,1.300);
\node[font=\scriptsize] at (-0.773,1.339) {$A{:}\,\theta_{j}$};
\draw[thick,dashed] (0.000,1.300)--(1.126,0.650);
\node[font=\scriptsize] at (0.773,1.339) {$B{:}\,{-}\theta_{j}$};
\draw[thick] (1.126,0.650)--(1.126,-0.650);
\node[font=\scriptsize] at (1.746,0.000) {$C{:}\,\theta_{j}$};
\draw[thick,dashed] (1.126,-0.650)--(0.000,-1.300);
\node[font=\scriptsize] at (0.773,-1.339) {$A{:}\,{-}\theta_{j}$};
\draw[thick] (0.000,-1.300)--(-1.126,-0.650);
\node[font=\scriptsize] at (-0.773,-1.339) {$B{:}\,\theta_{j}$};
\draw[thick,dashed] (-1.126,-0.650)--(-1.126,0.650);
\node[font=\scriptsize] at (-1.746,0.000) {$C{:}\,{-}\theta_{j}$};
\fill (-1.126,0.650) circle (1.6pt);
\node[font=\scriptsize] at (-0.814,0.470) {$p_0$};
\fill (0.000,1.300) circle (1.6pt);
\node[font=\scriptsize] at (0.000,0.940) {$n_0$};
\fill (1.126,0.650) circle (1.6pt);
\node[font=\scriptsize] at (0.814,0.470) {$p_1$};
\fill (1.126,-0.650) circle (1.6pt);
\node[font=\scriptsize] at (0.814,-0.470) {$n_1$};
\fill (0.000,-1.300) circle (1.6pt);
\node[font=\scriptsize] at (0.000,-0.940) {$p_2$};
\fill (-1.126,-0.650) circle (1.6pt);
\node[font=\scriptsize] at (-0.814,-0.470) {$n_2$};
\draw[thick] (4.874,0.650)--(6.000,1.300);
\node[font=\scriptsize] at (5.227,1.339) {$A{:}\,\theta_{j+1}$};
\draw[thick,dashed] (6.000,1.300)--(7.126,0.650);
\node[font=\scriptsize] at (6.773,1.339) {$B{:}\,{-}\theta_{j+1}$};
\draw[thick] (7.126,0.650)--(7.126,-0.650);
\node[font=\scriptsize] at (7.746,0.000) {$C{:}\,\theta_{j+1}$};
\draw[thick,dashed] (7.126,-0.650)--(6.000,-1.300);
\node[font=\scriptsize] at (6.773,-1.339) {$A{:}\,{-}\theta_{j+1}$};
\draw[thick] (6.000,-1.300)--(4.874,-0.650);
\node[font=\scriptsize] at (5.227,-1.339) {$B{:}\,\theta_{j+1}$};
\draw[thick,dashed] (4.874,-0.650)--(4.874,0.650);
\node[font=\scriptsize] at (4.254,0.000) {$C{:}\,{-}\theta_{j+1}$};
\fill (4.874,0.650) circle (1.6pt);
\node[font=\scriptsize] at (5.186,0.470) {$p_0$};
\fill (6.000,1.300) circle (1.6pt);
\node[font=\scriptsize] at (6.000,0.940) {$n_0$};
\fill (7.126,0.650) circle (1.6pt);
\node[font=\scriptsize] at (6.814,0.470) {$p_1$};
\fill (7.126,-0.650) circle (1.6pt);
\node[font=\scriptsize] at (6.814,-0.470) {$n_1$};
\fill (6.000,-1.300) circle (1.6pt);
\node[font=\scriptsize] at (6.000,-0.940) {$p_2$};
\fill (4.874,-0.650) circle (1.6pt);
\node[font=\scriptsize] at (5.186,-0.470) {$n_2$};
\draw[gray,dotted,very thick] (1.126,0.650)--(4.874,-0.650);
\node[font=\scriptsize,gray] at (3.250,0.320) {$A{:}\,0$};
\node[font=\small] at (0,-2.1) {node $j$}; \node[font=\small] at (6.0,-2.1) {node $j+1$};
\end{tikzpicture}
\caption{The gadgets of two scalars in the proof of Theorem~\ref{thm-main} (Lemma~\ref{univ-incidence}). Each edge is a column of $A$, labeled by its class and its value in every torus solution (solid: $\theta_j\,$, dashed: $-\theta_j$); values cancel at each vertex. The dotted bridge, like the completion edges (not shown), has value $0$. A gate $\theta_v=\theta_u+\theta_w$ is a row containing $p_0n_0$ of node $u$, $p_1n_1$ of node $w$ and $p_1n_0$ of node $v$.}
\label{fig-gadget}
\end{figure}
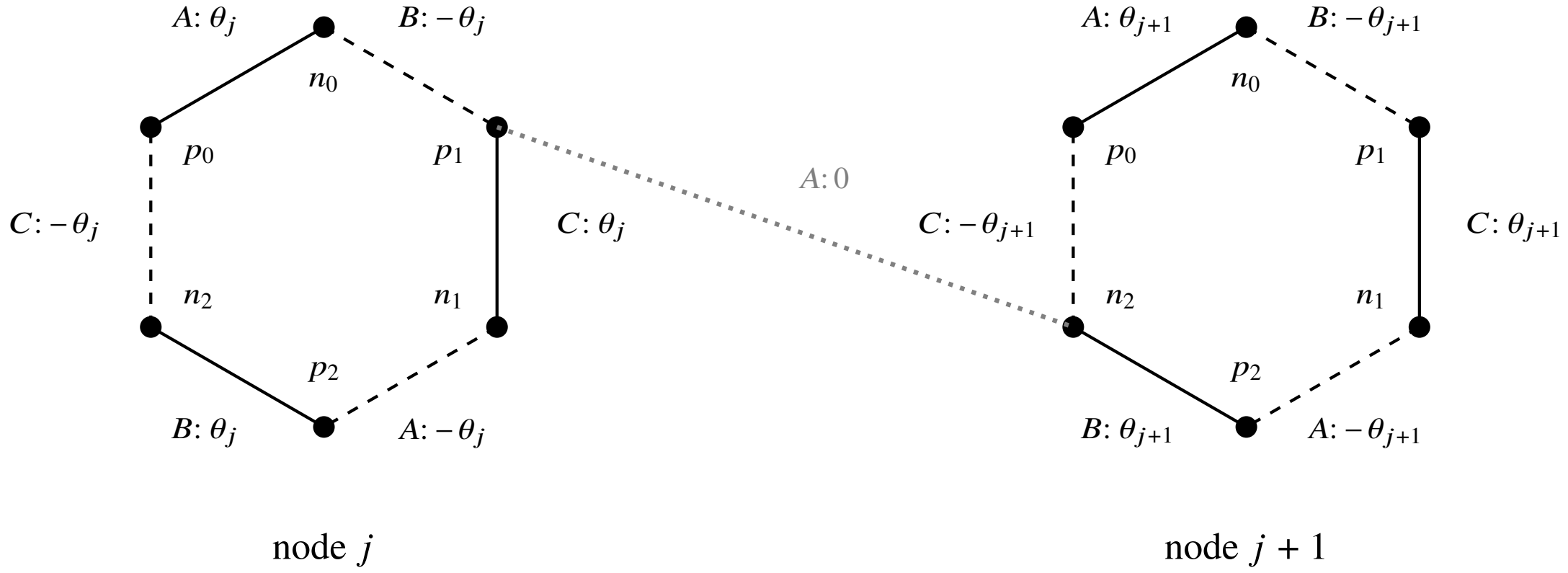

\emph{Realize a facet} (Figure~\ref{fig-gadget}). Replace each scalar $\theta_j$ by a six-cycle on positive vertices $p_0\,,p_1\,,p_2$ and negative vertices $n_0\,,n_1\,,n_2\,$, properly colored by the three classes $A,B,C$, whose edges carry $\pm\theta_j$ in every torus solution. A gate $\theta_v=\theta_u+\theta_w$ becomes the row containing the $\theta_u$-edge of color $A$ of node $u$, the $\theta_w$-edge of color $C$ of node $w$ and the $(-\theta_v)$-edge of color $B$ of node $v$; it meets each class once. Lemma~\ref{univ-incidence} gives a $0/1$ matrix with $9Q+1$ rows, $9Q$ independent columns, disjoint supports within colors, and torus kernel $G$; Lemmas~\ref{lem-growth} and~\ref{lem-shared} make it a simplex facet in dimension $9Q+1$ with $0/1$ slack.

\emph{Identify the arithmetic.} At each node, order its two signed columns and one zero-valued column as $(\theta_j\,,-\theta_j\,,0)$ in two colors and $(-\theta_j\,,\theta_j\,,0)$ in the third. Two fixed coefficient positions sum to one at every active node: their values are $\{-\theta_j\}$ and $\{\theta_j\}$. Two box points active at the same node cannot sum to another box point: those two positions would total two, forcing a carry. In a decomposition of a box point, every summand is itself a box point, so intersecting supports make every nonzero $H_g$ indecomposable (Proposition~\ref{univ-facet}). Independence gives the exact clearing set $\ord(g)\Z_{>0}$; exposure preserves both indecomposability and clearing in the ambient cone (Lemma~\ref{lem-face}). Each active node contributes between three and five to the degree, and every nonzero solution has linearly many active nodes, which gives the degree bounds; Propositions~\ref{prop-det} and~\ref{prop-grouprank} give the lower bounds on $d$.

\subsection{Intersecting additive evaluations}\label{sec-support-circuit}
An addition circuit uses binary additions and a shared zero input; wires may be reused, and a wire can be copied into a fresh scalar node by adding zero. For a vector over an abelian group, $\wt$ denotes the number of nonzero coordinates. Fix
$\alpha_0:=2^{-20},\qquad \delta_0:=2^{-22}.$
The encoder below is built from sparse neighbor-sum maps, as in linear-time codes \cite{Spielman}, followed by a support connector. The auxiliary graphs are supplied by finite probabilistic existence proofs independent of the coefficient group.

\begin{lemma}[Sparse neighbor sums]\label{univ-neighbor}
For every even $M\ge32$, there is a bipartite graph with $M$ left vertices and
$M/2$ right vertices, sixteen distinct neighbors at each left vertex, and
\[
 |\Gamma(S)|>12|S|\qquad(1\le |S|\le\alpha_0M).
\]
Its neighbor-sum map $R_M:B^M\to B^{M/2}$, for any abelian group $B$, can be
computed using at most $16M$ additions and satisfies
\[
 \wt(R_Mv)>8\wt(v)\qquad(0<\wt(v)\le\alpha_0M).
\]
\end{lemma}
\begin{proof}
Choose independently a uniform sixteen-element subset of the right vertices for
each left vertex. If $\alpha_0M<1$, the expansion requirement is empty. Otherwise,
for fixed sets $S,T$ with $|S|=s\le\alpha_0M$ and $|T|=12s$, the probability
that every neighbor of $S$ lies in $T$ is at most $(24s/M)^{16s}$. When
$12s<16$ this probability is zero. Every failure of the required expansion is
covered by such a pair, so the probability of a failure at size $s$ is at most
\[
 \binom Ms\binom{M/2}{12s}(24s/M)^{16s}
 \le \bigl[e^{13}24^4(s/M)^3\bigr]^s
 \le 2^{-20s}.
\]
Here $12s\le M/2$, and the last inequality follows from $e^{13}<2^{21}$,
$24^4<2^{19}$ and $s/M\le2^{-20}$. The sum over $s\ge1$ is less than one.

Let $S=\supp v$, $s=|S|$, and let $u$ count right vertices with exactly one
neighbor in $S$. Because there are $16s$ incident edges,
$16s\ge u+2(|\Gamma(S)|-u),\qquad u>8s.$
At every such right vertex the neighbor sum is one nonzero input, with
coefficient one. This proves the weight assertion over every abelian group,
including groups with torsion. Summing the right neighbor lists uses at most
the total number $16M$ of edges in binary additions.
\end{proof}

\begin{lemma}[A systematic additive encoder]\label{univ-encoder}
For every $r\ge1$, there is an addition circuit with $O(r)$ additions and an
output list of length $L=O(r)$, $L\ge4r$, such that every nonzero input over
any abelian group has at least $\delta_0 L$ nonzero output positions. The circuit
and the implicit constants are independent of the group.
\end{lemma}
\begin{proof}
For $M=32\cdot2^k$, construct a systematic map $C_M:B^M\to B^{4M}$
recursively. For $M=32$, repeat every input coordinate four times. For $M>32$,
set
\[
 a:=R_Mv,\qquad b:=C_{M/2}(a),\qquad
 C_M(v):=(v,b,R_{2M}b).
\]
The three block lengths are $M,2M,M$. We claim that
$\wt(C_Mv)\ge4\delta_0 M\qquad(v\ne0).$
The base case follows from the repetition. If $\wt(v)\ge4\delta_0 M$, the
first block proves the claim. Otherwise $0<\wt(v)<\alpha_0M$, and
Lemma~\ref{univ-neighbor} makes $a$ nonzero. Induction gives
$\wt(b)\ge2\delta_0 M$. If $\wt(b)\ge4\delta_0 M$, the second block suffices.
Otherwise $0<\wt(b)<\alpha_0M\le\alpha_0(2M)$, so the last block has
weight greater than $8\wt(b)\ge16\delta_0 M$.

The addition count satisfies
$T(M)\le T(M/2)+16M+32M<96M$.
Pad the $r$ inputs by zeros to the least $M$ of the prescribed form with
$M\ge r$. Then $M\le64r$, and $L=4M$ has the required properties.
\end{proof}

\begin{lemma}[A support connector]\label{univ-connector}
For every $L\ge1$, there is a simple graph on $L$ vertices with at most
$\Delta L$ edges, where $\Delta=2^{45}$, such that every two disjoint vertex
sets of size at least $\delta_0 L$ have an edge between them.
\end{lemma}
\begin{proof}
Take the undirected union of $\Delta$ independent uniform permutations of
$[L]$. For fixed disjoint $S,T$, one permutation has no directed edge from
$S$ to $T$ with probability at most $(1-|T|/L)^{|S|}$: expose the distinct
images of $S$ successively. Consequently the probability of no edge in the
union is at most
$\exp(-\Delta|S||T|/L)\le\exp(-\Delta\delta_0^2L).$
There are at most $3^L$ ordered disjoint pairs. Because
$\Delta\delta_0^2=2>\log3$, the union bound is less than one. Deleting loops
and repeated edges preserves the property and leaves at most $\Delta L$ edges.
\end{proof}

\begin{proposition}[Intersecting additive evaluations]\label{univ-circuit}
For every $r\ge1$, there is an addition circuit with $O(r)$ scalar nodes, including a shared zero node, such that any two nonzero input evaluations over an arbitrary abelian group have a scalar node nonzero in both. Moreover, every nonzero input evaluation has at least $\delta_0 L$ nonzero scalar nodes, where $L\ge4r$.
\end{proposition}
\begin{proof}
Apply Lemma~\ref{univ-encoder}, and give every output position its own fresh
node by adding its output wire to zero. This is necessary because the output
list can repeat wires. For every edge $uv$ of the graph in
Lemma~\ref{univ-connector}, append the scalar node $c_u+c_v\,$.

The fresh-output supports of two nonzero inputs have sizes at least
$\delta_0 L$. If they meet, a common output node proves the assertion. If they
are disjoint, choose an edge $uv$ between them. On the first evaluation,
$c_u+c_v=c_u\ne0$, and on the second, $c_u+c_v=c_v\ne0$. The addition and
node counts are linear in $r$, with absolute constants. The fresh output nodes
also give the individual support bound.
\end{proof}
\subsection{Presenting a prescribed finite abelian group}
Proposition~\ref{univ-circuit} supplies the circuit; we now attach it to a presentation of a given group by addition chains.
\begin{proposition}[A short scalar presentation]\label{univ-presentation}
Every nontrivial finite abelian group $G$ has a scalar presentation with $Q$
nodes and $Q$ equations whose real kernel is zero and whose kernel over $\T$
is isomorphic to $G$. One can choose $Q=O(\log|G|)$ so that every two nonzero
torus solutions have intersecting scalar supports. More precisely, there is an
integer $M$ with
$\log_2|G|\le M\le2\log_2|G|,\qquad Q\le C_0M,$
for an absolute constant $C_0\,$, and every nonzero torus solution has at least
$4\delta_0 M$ nonzero scalar nodes.
\end{proposition}
\begin{proof}
Write $G=\bigoplus_{i=1}^r\Z/q_i\Z$, $q_i\ge2$. Choose binary addition
chains
\[
 1=a_{i,0}<a_{i,1}<\cdots<a_{i,\ell_i}=q_i\,,
 \qquad \log_2q_i\le\ell_i\le2\log_2q_i\,.
\]
For each factor, make a node for each nonterminal multiple
$a_{i,j}\theta_i\,$, $0\le j<\ell_i\,$. Its internal additions are gates, while
the terminal addition is imposed equal to zero. Over $\R$ these equations force
all $\theta_i=0$; over $\T$ their solutions are precisely
$\theta_i\in(1/q_i)\Z/\Z$. Put $M:=\sum_i\ell_i\,$.

Apply Proposition~\ref{univ-circuit} to the list of all $M$ nonterminal
nodes, including its dependent nodes. Append the circuit gates and a shared
zero node, and constrain that zero node to vanish. The appended equations
determine new scalars uniquely from the old ones, so neither kernel changes.
Every nonzero solution gives a nonzero $M$-tuple. Its scalar support therefore
has size at least $\delta_0 L\ge4\delta_0 M$, and the supports of any two nonzero
solutions meet.

There are $r$ unconstrained starting nodes, one zero node, and $Q-r-1$ gate
nodes. The equations consist of their $Q-r-1$ gate equations, $r$ terminal
constraints and the zero constraint, for a total of $Q$. The circuit size gives
$Q\le C_0M$, and summing the addition-chain bounds proves the stated estimates.
\end{proof}
\subsection{From scalar supports to facet generators}
Next we order the columns of the incidence construction of Lemma~\ref{univ-incidence} so that two box points that are active at a common node cannot add without a carry.
\begin{proposition}[Node-separated transfer]\label{univ-facet}
Under the hypotheses of Lemma~\ref{univ-incidence}, order the support columns at each scalar node by
\[
 A,C:\ (\theta_j\,,-\theta_j\,,0_j),\qquad
 B:\ (-\theta_j\,,\theta_j\,,0_j),
\]
where the designated zero columns are bijected with the nodes. Write $s(g)=\{j:\theta_j(g)\ne0\}$. Box addition has no ray carry if and only if the scalar supports are disjoint. Consequently
\begin{multline*}
 \Hilb(K_F)=\{\text{facet rays}\}\\ \cup\
 \{H_g:g\ne0\text{ has no nontrivial disjoint-support decomposition}\}.
\end{multline*}
Every displayed generator is ambient indecomposable. In particular, if the supports of every two nonzero torus solutions meet, the indecomposable holes of the facet are all $|G|-1$ nonzero representatives. The facet has dimension $d+1$, normalized volume $|G|$, primitive $0/1$ slack, and exact clearing sets as in Lemma~\ref{univ-incidence}.
\end{proposition}
\begin{proof}
Formula~\eqref{eq-shared-coeff} gives the node coefficient blocks
\[
 A,C:\ (\{-\theta_j\},\{2\theta_j\},\{-\theta_j\}),\qquad
 B:\ (\{\theta_j\},\{-2\theta_j\},\{\theta_j\}),
\]
with the final coefficient in each class zero. Distinct nodes use disjoint coefficient blocks. If two solutions have disjoint scalar supports, their coefficient blocks add without carries. If both are active at a node $j$, the fixed positions with values $\{-\theta_j\}$ in color $A$ and $\{\theta_j\}$ in color $B$ sum to one for each solution. The two coefficient sums therefore total two, so at least one is at least one and a carry occurs.

In any decomposition of a box point into nonzero lattice points of the facet cone, every summand has nonnegative ray coefficients below one, and hence is another box point. The no-carry criterion identifies such decompositions with disjoint-support decompositions. Flooring arbitrary facet points now proves the Hilbert-basis assertion; exposure preserves these indecomposables in the ambient cone by Lemma~\ref{lem-face}. If all nonzero supports intersect, every nonzero representative is indecomposable. The remaining conclusions are supplied by Lemma~\ref{univ-incidence}.
\end{proof}
\subsection{The universality theorem}\label{sec-universality-thm}
Combining the presentation of Proposition~\ref{univ-presentation} with the transfer of Proposition~\ref{univ-facet} proves the main theorem; the converse bounds come from Section~\ref{sec-foundations}.
\begin{theorem}[Arithmetic universality]\label{thm-main}
There are absolute constants $a,C>0$ such that every nontrivial finite abelian group $G$ is realized by three staircase simplices in spatial dimension $d\le C\log_2|G|$, with a simplex facet $F$ satisfying:
\begin{enumerate}[label=\textup{(\roman*)}]
\item $G_F\cong G$ and $\Vol(F)=|G|$;
\item $\Hilb(K_F)$ consists of its rays and all $|G|-1$ nonzero box points $H_g\,$, each also in $\Hilb(K)$;
\item $ad\le\deg H_g\le d$ and $\mathcal C_K(H_g)=\ord_G(g)\Z_{>0}$ for $g\ne0$;
\item a primitive integral exposing form takes only the values zero and one on the lifted vertices.
\end{enumerate}
Conversely, every simplex face with residue group $G$ for three staircase simplices in $\R^d$ satisfies
\begin{equation}\label{eq-dim-lower}
 d\ge\max\{\lceil1+3\log_3|G|\rceil,\,3\mu(G)+1\},
\end{equation}
where $\mu(G)$ is the minimal number of generators. Thus the least dimension for these simultaneous properties is $\Theta(\log|G|)$.
\end{theorem}
\begin{proof}
Apply Proposition~\ref{univ-presentation}, followed by
Proposition~\ref{univ-facet}. This gives $d=9Q+1=O(\log|G|)$ and all the
claimed arithmetic properties. Every nonzero solution has at least
$4\delta_0 M$ active scalars. At each active node the coefficient sums in colors $A,C$ are $2,1,1$ as $\theta$ lies in $(0,1/2)$, equals $1/2$, or lies in $(1/2,1)$; those in color $B$ are $1,1,2$. Thus an active node contributes respectively $5,3,4$ to the degree, and
\[
 12\delta_0 M\le\deg H_g\le5Q<\frac59d.
\]
Because $Q\le C_0M$ and $M\ge1$, we have $d\le(9C_0+1)M$; hence the lower
degree bound holds with the absolute constant $a=12\delta_0/(9C_0+1)$.
Proposition~\ref{prop-det} and Lemma~\ref{lem-face} give
$|G|\le3^{(d-1)/3}$, and Proposition~\ref{prop-grouprank} gives the rank bound. This proves the matching lower bound on dimension. The explicit values of $C$ and $a$ are traced in Remark~\ref{univ-scope}.
\end{proof}
\begin{remark}[Constants]\label{univ-scope}
The finite probabilistic proofs allow a fixed first admissible graph, so no Hilbert basis computation is required. The encoder pads $M$ inputs to at most $64M$, uses fewer than $96\cdot64M$ gates, and has output length at most $256M$. Fresh output nodes and the connector give $Q\le2^{54}M$, hence $d\le2^{59}\log_2|G|$. The lower degree constant can be taken as $12\delta_0/(9\cdot2^{54}+1)\ge2^{-76}$. These are coarse constants; the next construction trades asymptotic dimension for smaller explicit coefficients.
\end{remark}

\begin{remark}[Explicit constants for elementary abelian groups]\label{univ-elementary}
Let $\Gamma$ be a generator matrix of a linear intersecting code of dimension $r$ and length $N$ over $\F_p\,$, that is, any two nonzero codewords have intersecting supports. Present each factor of $(\Z/p\Z)^r$ by a binary addition chain for $p$, adjoin a zero node, and discard the zero columns of $\Gamma$, which keeps the code intersecting, and compute every remaining code coordinate $\sum_i\Gamma_{ik}\theta_i$ by successive additions into a fresh node. Two nonzero torus solutions are active at a common output node because the code is intersecting, so Proposition~\ref{univ-facet} gives properties (i), (ii), (iv) of Theorem~\ref{thm-main} with exact clearing sets $p\Z_{>0}$ and $Q\le r\ell(p)+\sum_{i,k}\Gamma_{ik}+1$, with the entries of $\Gamma$ taken as integers in $\{0,\ldots,p-1\}$. Binary intersecting codes of length $N\le\lceil4.82\,r\rceil$ exist by the probabilistic bound \cite[Theorem~3.12]{BSS}, because $2/\log_2(4/3)<4.82$ (the codes are not explicit, only the constants are); this gives dimension at most $44r^2+18r+10$ for $(\Z/2\Z)^r$, with much smaller constants than Remark~\ref{univ-scope}.
\end{remark}

\subsection{Minimal non-IDP faces}\label{sec-minfaces}
A non-IDP simplex face all of whose proper faces have IDP has a cyclic residue group, yet, unlike the bound $d\ge3q-2$ for a single indecomposable hole, every order $q$ occurs on such a facet in dimension $O(\log q)$, with all $q-1$ box points indecomposable and of the same degree (Theorem~\ref{thm-minfacet}). For $G=(\Z/2\Z)^r$, Theorem~\ref{thm-main} gives $2^{\Theta(d)}$ inclusion-minimal non-IDP faces (Corollary~\ref{cor-manymin}).

For a simplex face $F$ and a ray $v$ of $F$, write $\chi_v\colon G_F\to\T$ for the coefficient of $v$ in box coordinates. Extending coefficient vectors by zero identifies the residue group of the subface spanned by a set $I$ of rays with $\bigcap_{v\notin I}\ker\chi_v\subseteq G_F\,$, and a simplex face has IDP exactly when its residue group is trivial, because a nonzero residue has a box point with a fractional coefficient.

\begin{corollary}[Many minimal non-IDP faces]\label{cor-manymin}
The maximum number of inclusion-minimal non-IDP faces of three staircase simplices in $\R^d$ is $2^{\Theta(d)}$. The lower bound is attained by simplex faces of normalized volume two whose pairwise intersections are unimodular.
\end{corollary}
\begin{proof}
Apply Theorem~\ref{thm-main} to $G=(\Z/2\Z)^r$, in dimension $O(r)$. By Proposition~\ref{univ-facet} every ray coefficient of a box point is $0$ or $\frac12$, and every nonzero box point is indecomposable. If $\supp H_g\subseteq\supp H_h$ for nonzero $g\ne h$, then $H_h-H_g$ is a lattice point with coefficients in $\{0,\frac12\}$, hence equal to $H_{g+h}\ne0$, and $H_h=H_g+H_{g+h}$ would decompose; so the supports are pairwise incomparable. The face spanned by $\supp H_g$ therefore has residue group $\{0,g\}$, while each of its proper faces misses a ray of $\supp H_g$ and has trivial residue group. This gives $2^r-1$ distinct minimal non-IDP faces, of normalized volume two, with unimodular pairwise intersections. Appending coordinates at the end of all three growth orders keeps the hull as the face where they vanish, so the bound holds in every large dimension. The hull has $3(d+1)$ vertices and each face is determined by its vertices, which gives the upper bound.
\end{proof}

\begin{theorem}[Minimal non-IDP facets in logarithmic dimension]\label{thm-minfacet}
A non-IDP simplex face is inclusion-minimal among non-IDP faces if and only if every $\chi_v$ is injective; in particular its residue group is cyclic. Conversely, for every $q\ge2$, three staircase simplices in an even spatial dimension $d=O(\log q)$ have a simplex facet $F$ such that
\begin{enumerate}
\item $G_F\cong\Z/q\Z$ and every proper face of $F$ is unimodular;
\item every nonzero box point lies in the relative interior of the cone over $F$, has multiplicities $\bigl((d+2)/6,(d+2)/6,(d+2)/6\bigr)$, and lies in $\Hilb(K)$ with clearing order equal to its residue order;
\item $\Hilb(K_F)$ consists of the rays and the $q-1$ nonzero box points, and $h^*_F(\tau)=1+(q-1)\tau^{(d+2)/2}$.
\end{enumerate}
By Proposition~\ref{prop-det}, the dimension is optimal up to an absolute factor.
\end{theorem}
\begin{proof}
\emph{Characterization.} Deleting a ray $v$ leaves a subface with residue group $\ker\chi_v\,$. All proper faces have IDP exactly when all these kernels are trivial, because every proper face lies in such a subface. A finite subgroup of $\T$ is cyclic.

\emph{The scalar presentation.} We use Lemma~\ref{univ-incidence}: a square scalar presentation on $Q$ nodes with zero real kernel gives a simplex facet in dimension $9Q+1$; each class has a column of each sign for every node and $Q$ further zero columns, in any order, and by~\eqref{eq-shared-coeff} the ray characters of a class with ordered column values $f_1\,,\ldots,f_{3Q}$ are $-f_1\,,\ f_1-f_2\,,\ldots,f_{3Q-1}-f_{3Q}\,,\ f_{3Q}\,$, the negatives of the successive increments of the closed sequence $0,f_1\,,\ldots,f_{3Q}\,,0$. We arrange all increments to be units times a generator $\theta$.

Take a node $\theta$ and a zero-constrained node. Follow the binary addition chain for $q$: the values $1,2,4,\ldots,2^{\lfloor\log_2q\rfloor}$, then the remaining binary powers of $q$, added successively in descending order, with a node of value $a\theta$ for every intermediate value $a<q$, defined by its gate, and the last addition imposed as the terminal equation $q\theta=0$. Read these nodes as a \emph{walk} from the virtual value $0$ to the virtual value $q\theta=0$; its increments are powers of two. If $q$ is odd, they are units. If $q$ is even, let $m$ be its odd part. If $m>1$, compute $m\theta$ by the binary chain for $m$, using fresh copies (sums with the zero node) of the power nodes, and retrace from $m\theta$ to $0$ with fresh copies of the intermediate nodes; this is a second walk from $0$ to $0$ with power-of-two steps. For each $j\ge1$ that occurs below, add the node $w_j=(m+2^j)\theta$ by one gate, as its own walk from $0$ to $0$; its coefficient is a unit. In the first walk and, when present, the second walk, replace each step of size two, between labels $a$ and $a+2$, by two unit steps in the same direction through a fresh node $(a+1)\theta$, and each step of size $2^k$ with $k\ge2$, between $a$ and $a+2^k$, by a two-step detour through a fresh node $a\theta+w_{k-1}\,$, using the zero node when $a=0$. The two new increments, smaller endpoint first, are $2^{k-1}+m$ and $2^{k-1}-m$, the second possibly negative: both are odd and both are congruent to $2^{k-1}$ modulo every odd prime dividing $q$, hence units modulo $q$. A scheduled node may itself be zero modulo $q$; it stays in its walk, because only the increments matter. Now there are $c$ walks from $0$ to $0$ containing $P=O(\log q)$ \emph{scheduled} nodes, each once, and all increments, including those at the virtual endpoints, are units times $\theta$. Every scheduled node other than $\theta$ has one defining gate and $\theta$ has the terminal equation, so the presentation stays square, with zero real kernel and torus kernel $\{\theta:q\theta=0\}\cong\Z/q\Z$.

\emph{Using every column.} Take $Z=\max\{1,\lceil(2c-1-P)/3\rceil\}$ zero-constrained nodes in all, and add $E=P+3Z-(2c-1)\ge0$ fresh copies of $\theta$, each a walk by itself. With $P'=P+E$ scheduled nodes, $c'=c+E$ walks and $Q=P'+Z$ nodes, each class has $Q+2Z=P'+3Z=2c'-1$ zero columns besides the columns of the scheduled nodes. In every class list the $c'$ walks by the positive columns of their nodes, then the $c'$ walks by the negative columns in the same order, with one zero column between consecutive walks and none at the ends. All increments of the closed sequence are units times $\theta$, and those of the second half are the negatives of those of the first. Note that $Q=2P+4Z-2c+1$ is odd, so $d=9Q+1$ is even, and all counts are $O(\log q)$.

\emph{Conclusion.} Every ray character is injective, so every proper face is unimodular by the characterization. For a nonzero residue every coefficient is positive, and the two coefficients of each opposite pair sum to one. Each class has $3Q+1=(d+2)/3$ rays in opposite pairs, so each multiplicity is $(d+2)/6$ and the degree is $(d+2)/2$. A decomposition of a box point into two nonzero lattice points of the facet cone would consist of two nonzero box points, whose coefficients on an opposite pair would total two; so every nonzero box point is indecomposable, flooring gives $\Hilb(K_F)$, and Lemma~\ref{lem-face} gives ambient indecomposability and the clearing orders. The $h^*$-polynomial follows from the box degrees. Finally $q\le3^{(d-1)/3}$ by Proposition~\ref{prop-det}.
\end{proof}

\begin{corollary}[Least degree of an indecomposable hole of a simplex face]\label{cor-halfrays}
A non-IDP simplex face with $s$ rays has an indecomposable hole of degree at most $\lfloor s/2\rfloor$, and the facets of Theorem~\ref{thm-minfacet} attain this bound for every $q\ge2$.
\end{corollary}
\begin{proof}
For a nonzero residue $g$, the box points of $g$ and $-g$ have complementary coefficients $c$ and $1-c$ on the support of $H_g\,$, so their degrees sum to $|\supp H_g|\le s$, and one of them has degree at most $\lfloor s/2\rfloor$. Both fail to split, because their unique ray expressions are fractional and $F$ is a face (Lemma~\ref{lem-face}). Hence the Hilbert-basis decomposition of a box point among them of degree at most $\lfloor s/2\rfloor$ contains a summand that is not a ray, and this summand is an indecomposable hole of degree at most $\lfloor s/2\rfloor$. In Theorem~\ref{thm-minfacet} every nonzero box point has degree $s/2$, and every indecomposable hole has a nonzero box point below it after flooring.
\end{proof}

\subsection{One triple for all small groups}\label{sec-simul}
A single fixed triple in dimension $O(m)$ has, for every $G$ with $1<|G|\le2^m$, its own simplex face $F_G$ satisfying (i)--(iii) of Theorem~\ref{thm-main} with residue group $G$ (Theorem~\ref{thm-simul}). Hence the least positive integer $k$ with $k\overline S_K\subseteq S_K$ can be $2^{2^{\Theta(d)}}$, although $\delta_K(h)\le3^{(d-1)/3}$ for every single $h$ (Corollary~\ref{cor-uniform}).

In Theorem~\ref{thm-main} the growth orders depend on $G$. Here they are fixed, and $G$ is selected by the target point, through its minimal face. We use the profile encoding of Lemma~\ref{lem-bounded}, whose proof applies verbatim to bounds $0\le z_v\le B_v$ with positive integers $B_v\,$: a block of $z_v$ then occupies an interval $[L,L+B_v]$, a pair coordinate has value $L+L'+B_v\,$, and an equation coordinate has value $b+\sum_jL_j+\sum_{j:\epsilon_j=-1}B_{v_j}\,$. The growth orders depend only on the bounds and on the signed occurrences of the variables in the equations; the right-hand sides enter only the target point.

\begin{lemma}[Minimal faces of bounded systems]\label{lem-simul-face}
Suppose that the system of Lemma~\ref{lem-bounded}, with bounds $B_v\,$, has a unique real solution $z^*$, and let $h$ be the corresponding target point and $F$ the minimal face of $\mathcal D$ containing $h/\deg h$. Then $F$ is a simplex face in $\{x_g=0\}$, and $G_F$ is isomorphic to the group of $\bar z\in\T^n$ satisfying the homogeneous equations $\sum_j\epsilon_j\bar z_{v_j}=0$ together with $\bar z_v=0$ whenever $z^*_v\in\{0,B_v\}$. For each $v$ with $0<z_v^*<B_v$ there are two rays of $F$, distinct for distinct $v$, at which the box point of $\bar z$ has coefficients $\{-\bar z_v\}$ and $\{\bar z_v\}$.
\end{lemma}
\begin{proof}
By Lemma~\ref{lem-bounded}, $h$ has a unique real vertex expression. The point $h$ lies in the relative interior of $K_F\,$, so it has an expression positive on every ray of $F$; by uniqueness this is the expression, and the rays of $F$ are independent, because otherwise a dependence would perturb it. Hence $F$ is a simplex, its rays are the prefixes with positive coefficient, and it lies in $\{x_g=0\}$. In each class the profile takes the value $\lambda_i$ before the first block, a constant value in each block, the integral value of each singleton, and zero from $g$ on. A block of $z_v$ on $[L,L+B_v]$ with value $L+w$ thus has an entering prefix with coefficient $B_v-w$ and a leaving prefix with coefficient $w$, and all other drops vanish. So the rays are the entering and leaving prefixes of the blocks with $w<B_v\,$, respectively $w>0$.

By independence, $G_F$ is the kernel of the ray matrix over $\T$. Tail sums identify such a kernel vector with a $\T$-valued profile in each class that changes only across rays, vanishes at the start (multiplicity rows) and from $g$ on, and whose three values add to zero at every coordinate. A singleton lies in one list before $g$, so the profile vanishes there; a block with only one ray adjacent therefore has value zero, which is the boundary condition. The pair coordinates give values $\bar z_v$ on positive and $-\bar z_v$ on negative blocks of $v$, and the equation coordinates give the homogeneous equations. Conversely, every such $\bar z$ defines these profiles, and their differences are supported on rays. For an interior $v$, the entering and leaving rays of its positive block in class $0$ have coefficients $\{-\bar z_v\}$ and $\{\bar z_v\}$.
\end{proof}

\begin{theorem}[Simultaneous arithmetic universality]\label{thm-simul}
There are absolute constants $C,c>0$ such that, for every integer $m\ge1$, one triple of staircase simplices in dimension $d\le Cm$ has, for every finite abelian group $G$ with $1<|G|\le2^m$, a simplex face $F_G\subseteq\{x_g=0\}$ with $G_{F_G}\cong G$. Every nonzero box point $H_a$ of $F_G$ lies in $\Hilb(K)$, has $cm\le\deg H_a\le d$, and satisfies $\mathcal C_K(H_a)=\ord(a)\Z_{>0}\,$. Any triple with this property has $d\ge(m-4)/3$.
\end{theorem}
\begin{proof}
\emph{Switches.} For $v\in[0,1]$ and $\epsilon\in\{0,1\}$, the rows $w+w'-v=0$, $w+a=\epsilon$, $w'+b=1-\epsilon$ with $w,w',a,b\in[0,1]$ have the unique solution $w=\epsilon v$, $w'=(1-\epsilon)v$, $a=\epsilon(1-v)$, $b=(1-\epsilon)(1-v)$. Write $S_\epsilon(v)$ for $w$. The flag occurs only in right-hand sides.

\emph{The fixed system.} Put $L=2m$. Write $G=\bigoplus_{i=1}^r\Z/q_i\Z$ with $q_i\ge2$; the binary lengths of the $q_i$ add up to at most $2\log_2|G|\le L$. Write these expansions, most significant bit first, in consecutive blocks of positions $1,\ldots,L$, padding the first block with leading zeros. Let $\beta_j$ be the bit at position $j$, and let $s_j\,,e_j\in\{0,1\}$ indicate whether $j$ starts or ends its block. Use $z_0=0$ and variables $\theta_j\,,z_j\,,u_j\,,y_j\,,v_j\in[0,1]$ for $1\le j\le L$, with fresh switches and the rows
\[
 \begin{gathered}S_{1-s_j}(\theta_j)=S_{1-s_j}(\theta_{j-1})\ (j\ge2),\quad u_j=S_{1-s_j}(z_{j-1}),\quad y_j=u_j+u_j\,,\\ v_j=S_{\beta_j}(\theta_j),\quad z_j=y_j+v_j\,,\quad S_{e_j}(z_j)=e_j\,.\end{gathered}
\]
Inside a block the $\theta_j$ agree, $z$ restarts at the beginning, and $z_j$ is $\theta_j$ times the binary value of the bits read so far in the block. At the end of block $i$ the last row forces $q_i\theta_j=1$. So the unique real solution has $\theta_j=1/q_i$ on block $i$, and every variable lies in $[0,1]$ because the partial binary values are at most $q_i\,$. Then attach the circuit of Proposition~\ref{univ-circuit} to $\theta_1\,,\ldots,\theta_L\,$, with its shared zero a variable fixed to $0$. Give each addition node $c=a+b$ a variable $c\in[0,2]$ and the row $a+b-c=\lfloor a^*+b^*\rfloor-1$, where $a^*,b^*$ are the real values already determined. Then $c^*=1+\{a^*+b^*\}$ lies strictly between the bounds. All rows have at most three signed occurrences, there are $O(L)$ variables and rows, and the bounds and occurrences do not depend on $G$. Lemma~\ref{lem-bounded} gives one triple in dimension $d=O(m)$ for all $G$, and the target point $h_G$ depends on $G$.

\emph{The residue group.} Let $F_G$ be the minimal face of $h_G\,$. By Lemma~\ref{lem-simul-face}, $G_{F_G}$ consists of the homogeneous solutions with the boundary conditions. For a switch, $\epsilon=0$ forces $\bar w=0$ because $w^*=0$, and $\epsilon=1$ forces $\bar w'=0$, hence $\bar w=\bar v$. So the homogeneous system follows the same rules: $\bar\theta_j=\eta_i$ on block $i$, all other entries are determined by the $\eta_i\,$, and the last row of block $i$ gives $q_i\eta_i=0$. Conversely, let $q_i\eta_i=0$. Replacing $1/q_i$ by a formal variable $\tau_i\,$, every variable before the circuit is an integer-affine function of one $\tau_i\,$, and all rows except the final rows of the blocks hold identically. Substituting $\tau_i=\eta_i$ and dropping constants gives a homogeneous solution. A variable at a bound has a real value $a+b/q_i\in\Z$, so $q_i\mid b$ and its entry $b\eta_i$ vanishes. The circuit nodes are interior and only add $\bar c=\bar a+\bar b$. Hence $G_{F_G}\cong\bigoplus_i\Z/q_i\Z=G$.

\emph{Box points.} A nonzero residue has some $\eta_i\ne0$, hence a nonzero input vector $(\bar\theta_j)$. By Proposition~\ref{univ-circuit}, two nonzero residues share a nonzero circuit node $c$. Its two rays from Lemma~\ref{lem-simul-face} carry coefficients summing to one in each box point, so a sum of two nonzero box points is not a box point. In a decomposition of $H_a$ into two nonzero lattice points of $K$, exposure puts both in $K_{F_G}\,$, and nonnegativity of the unique coefficients makes both box points. So $H_a$ is indecomposable, hence in $\Hilb(K)$ by Lemma~\ref{lem-face}, and independence gives $\mathcal C_K(H_a)=\ord(a)\Z_{>0}\,$. Each nonzero node contributes one to the degree through its own pair of rays, and there are at least $4\delta_0L=8\delta_0m$ of them by Proposition~\ref{univ-circuit}; the upper bound is Proposition~\ref{prop-bounds}.

\emph{Lower bound on $d$.} The $2^m-1$ groups $\Z/q\Z$, $2\le q\le2^m$, need distinct faces, and $\mathcal D$ has at most $2^{3(d+1)}$ faces.
\end{proof}

For a triple, let $u(\mathcal D)$ be the least $k\ge1$ with $k\overline S_K\subseteq S_K\,$, and let $U(d)$ be its maximum over triples in $\R^d$.
\begin{corollary}[A uniform clearing multiplier]\label{cor-uniform}
$\log_2\log_2U(d)=\Theta(d)$. In contrast, $\delta_K(h)\le3^{(d-1)/3}$ for every single lattice point $h$ (Proposition~\ref{prop-det}).
\end{corollary}
\begin{proof}
For the upper bound, every $h\in\overline S_K$ lies in the cone of a basis $B\subseteq V$, and Cramer's rule gives $|\det B|\,h\in S_K\,$. So the least common multiple of all $|\det B|$ is a uniform multiplier. There are at most $2^{3(d+1)}$ bases, each with $|\det B|\le3^{(d-1)/3}$ by Proposition~\ref{prop-det}, so $\log_2u(\mathcal D)\le2^{3(d+1)}d$. For the lower bound, take $m=\lfloor d/C\rfloor\ge2$ in Theorem~\ref{thm-simul} and append the missing coordinates at the end of all three growth orders. The old cone becomes the face where they vanish, and Lemma~\ref{lem-face} keeps the clearing sets. Every $q\le2^m$ is then a clearing order, so $u(\mathcal D)$ is divisible by $\operatorname{lcm}(1,\ldots,2^m)$. With $R=2^{m-1}$, every prime power dividing $\binom{2R}{R}$ is at most $2R$, so $\binom{2R}{R}\mid\operatorname{lcm}(1,\ldots,2R)$ and $\log_2u(\mathcal D)\ge2R-\log_2(2R+1)\ge2^{m-2}$.
\end{proof}

\section{Integer decomposition of $0/1$ polytopes}\label{sec-idp}
This section proves result~\hyperlink{res-C}{(C)}: the integer decomposition property of a $0/1$ polytope, together with its indecomposable holes and their clearing sets, is realized on a face of three staircase simplices (Theorem~\ref{thm-idp-transfer}). The series lift is a generalized Lawrence construction \cite{PTV1,PTV2}; the realization identifies its induced lattice and preserves indecomposable holes and clearing sets.

Let $\mathcal Q=\conv(v_1\,,\ldots,v_N)\subseteq[0,1]^r$ be a $0/1$ polytope with $N\ge1$ distinct vertices, $b_j=(v_j\,,1)$, and
\[
 B:=(b_1\cdots b_N),\qquad C_B:=\cone B,\qquad \Lambda_B:=\spanR B\cap\Z^{r+1},\qquad S_B:=\Z_{\ge0}B .
\]
All statements use the induced lattice $\Lambda_B\,$. The lattice points of $\mathcal Q$ are its vertices, because a lattice point of the cube is a cube vertex and is not a nontrivial convex combination of other cube vertices; hence $\mathcal Q$ has IDP exactly when $S_B=C_B\cap\Lambda_B\,$, and its \emph{indecomposable holes} are the elements of $\mathcal H_B:=\Hilb(C_B\,,\Lambda_B)\setminus\{b_1\,,\ldots,b_N\}$. For $h\in C_B\cap\Lambda_B$ put $\mathcal C_B(h):=\{q\in\Z_{>0}:qh\in S_B\}$ and let $t(h)$ be its last coordinate.

\subsection{Overview}\label{sec-idp-overview}
We outline the proof of Theorem~\ref{thm-idp-transfer}; the details occupy the rest of this section.

Integer matrices of equal width and equal torus kernels (Subsection~\ref{sec-incidence}) have isomorphic cones, lattices and monoids (Lemma~\ref{lem-toruskernel}).

\emph{Encode $\ker_\T B$ by chains.} A global chain $0=u_0\,,u_1\,,\ldots,u_N=0$ has increments $g_j:=u_j-u_{j-1}\,$. In each spatial row, the indices with entry $0$ and those with entry $1$ get chains of their own from $0$ to $0$, tied to the global chain by gates $u_j=s_{\ell,0}+s_{\ell,1}\,$, where $s_{\ell,0}\,,s_{\ell,1}$ are the current values of the two chains of row $\ell$. So $g_j$ is the increment of its list in every row, and $g_j=u_j-u_{j-1}$ identifies the scalar kernel with $\ker B$ over $\R$ and $\T$.

\emph{Order the supports so that differences repeat characters.} In the incidence construction of Subsection~\ref{sec-incidence} applied to these chains (Proposition~\ref{prop-extended-incidence}), whose rays may be dependent, list each chain with internal nodes $w_1\,,\ldots,w_{m-1}$ in each class as $0,w_{m-1}\,,\ldots,w_1\,,0$ and $0,-w_1\,,\ldots,-w_{m-1}\,,0$. Endpoint coefficients are backward differences, which here are the increments, so $g_j$ gets two copies per class from the global chain and from one chain per spatial row: $2c$ per class with $c:=r+1$, $6c$ in all, and every other coefficient vanishes on the kernel. Grouping positions by original column, the ray matrix $R$ of the face $F_{\mathcal Q}$ satisfies
\[
 \ker_\T R=\{(\operatorname{repeat}_{6c}(g),0):g\in\ker_\T B\},
\]
where $\operatorname{repeat}_k(g)$ lists each entry of $g$ $k$ times.

\emph{Identify with the series lift} (Figure~\ref{fig-series}). The series lift $C^{(6c)}$ of $B$ (Lemma~\ref{lem-series}) replaces each column $b_j$ by a chain of $6c$ generators summing to $(b_j\,,0)$; its indecomposable holes are exactly the points $(h,0)$, $h\in\mathcal H_B\,$, with the same clearing sets, and its degrees are multiplied by $6c$. Its generators, padded with free unit rays to the same number of columns, have the same torus kernel, and free unimodular rays change neither indecomposable holes nor clearing sets. This gives parts~(1) and~(2) of Theorem~\ref{thm-idp-transfer}, with $\widehat h$ the image of $(h,0)$; exposure (Lemma~\ref{lem-face}) puts $\widehat h$ in $\Hilb(K)$ with clearing set $\mathcal C_B(h)$. If $Bp=h$, the coefficient $p_j$ on the $2c$ positions per class of column $j$ gives part~(3).

\begin{figure}[t]
\centering
\begin{tikzpicture}[x=0.9cm,y=0.9cm,scale=1.450]
\fill (-1.250,0.000) circle (1.6pt);
\node[font=\small] at (-1.775,0.000) {$0$};
\fill (-0.386,1.189) circle (1.6pt);
\node[font=\small] at (-0.549,1.688) {$\mathbf e_{j,1}$};
\fill (1.011,0.735) circle (1.6pt);
\node[font=\small] at (1.436,1.043) {$\mathbf e_{j,2}$};
\node at (1.011,-0.735) {$\cdots$};
\fill (-0.386,-1.189) circle (1.6pt);
\node[font=\small] at (-0.549,-1.688) {$\mathbf e_{j,k-1}$};
\draw[->,very thick,blue!70!black] (169.0:1.25) arc[start angle=169.0,end angle=119.0,radius=1.25];
\node[font=\scriptsize] at (-0.627,0.456) {$a_{j,1}$};
\draw[->,thick] (97.0:1.25) arc[start angle=97.0,end angle=47.0,radius=1.25];
\node[font=\scriptsize] at (0.239,0.737) {$a_{j,2}$};
\draw[->,thick] (25.0:1.25) arc[start angle=25.0,end angle=-25.0,radius=1.25];
\node[font=\scriptsize] at (0.775,0.000) {$a_{j,3}$};
\draw[->,thick] (-47.0:1.25) arc[start angle=-47.0,end angle=-97.0,radius=1.25];
\node[font=\scriptsize] at (0.191,-0.590) {$a_{j,k-1}$};
\draw[->,thick] (-119.0:1.25) arc[start angle=-119.0,end angle=-169.0,radius=1.25];
\node[font=\scriptsize] at (-0.627,-0.456) {$a_{j,k}$};
\node[anchor=west,font=\small,align=left] at (2.6,0.45) {$a_{j,1}=(b_j\,,\mathbf e_{j,1})$, the only stage with a spatial part,\\[1pt] $a_{j,l}=(0,\mathbf e_{j,l}-\mathbf e_{j,l-1})$ for $1<l<k$,\\[1pt] $a_{j,k}=(0,-\mathbf e_{j,k-1})$;};
\node[anchor=west,font=\small] at (2.6,-0.75) {$\textstyle\sum_{l=1}^k a_{j,l}=(b_j\,,0)$.};
\end{tikzpicture}
\caption{The series lift of one column $b_j$ of $B$ (Lemma~\ref{lem-series}). In the private coordinates of column $j$ the $k$ stages form a closed walk summing to $(b_j\,,0)$, so a point with zero private part uses all stages of a column equally. Hence the indecomposable holes of the lift are the points $(h,0)$, $h\in\mathcal H_B\,$, with the same clearing sets and degrees multiplied by $k$ ($k=6c$ in Theorem~\ref{thm-idp-transfer}).}
\label{fig-series}
\end{figure}

\subsection{A series lift with no new indecomposable holes}
We first show that the series lift of Figure~\ref{fig-series} changes neither the indecomposable holes nor their clearing sets, and that cones with the same torus kernel are isomorphic.
Fix $k\ge2$. Give each column $j$ a private block $\Z^{k-1}$ with unit vectors $\mathbf e_{j,1}\,,\ldots,\mathbf e_{j,k-1}\,$, and in $\Lambda_B\oplus\Z^{N(k-1)}$ define
\[
 a_{j,1}=(b_j\,,\mathbf e_{j,1}),\qquad a_{j,l}=(0,-\mathbf e_{j,l-1}+\mathbf e_{j,l})\ (2\le l<k),\qquad a_{j,k}=(0,-\mathbf e_{j,k-1}).
\]
Let $C^{(k)}$ be their cone and $S^{(k)}$ their nonnegative integer monoid; the induced lattice of $C^{(k)}$ is $\Lambda_B\oplus\Z^{N(k-1)}$. Write points as $(x,z)$ with $z_{j,a}$ the private coordinates.

\begin{lemma}[Series lift]\label{lem-series}
For $(x,z)$ put $s_{j,1}:=0$, $s_{j,l}:=\sum_{a<l}z_{j,a}\,$, $m_j:=\max_{1\le l\le k}s_{j,l}$ and $r_{j,l}:=m_j-s_{j,l}\,$. Then
\begin{gather*}
 (x,z)\in C^{(k)}\iff x-Bm\in C_B\,,\\ (x,z)\in S^{(k)}\iff x-Bm\in S_B\ \text{for integer }(x,z),
\end{gather*}
and $\Hilb(C^{(k)})=\{a_{j,l}\}\cup\{(h,0):h\in\mathcal H_B\}$. The grading $\tau(x,z)=k\,t(x)-\sum_{j,a}(k-a)z_{j,a}$ is one on every $a_{j,l}$ and $k\,t(h)$ on $(h,0)$, and $\{q\in\Z_{>0}:q(h,0)\in S^{(k)}\}=\mathcal C_B(h)$. In particular $C^{(k)}$ has the integer decomposition property in its induced lattice if and only if $\mathcal Q$ does.
\end{lemma}
\begin{proof}
In a real expression $(x,z)=\sum_{j,l}p_{j,l}a_{j,l}$ with $p\ge0$, the private coordinates give $z_{j,a}=p_{j,a}-p_{j,a+1}\,$, hence $p_{j,l}=p_{j,1}-s_{j,l}\,$, and nonnegativity of block $j$ is equivalent to $p_{j,1}\ge m_j\,$. Then $x=\sum_jp_{j,1}b_j$ and $x-Bm=B(p_{\cdot,1}-m)\in C_B\,$. Conversely, if $x-Bm=\sum_jq_jb_j$ with $q\ge0$, then $p_{j,1}=m_j+q_j$ and $p_{j,l}=p_{j,1}-s_{j,l}$ are nonnegative and reproduce $(x,z)$. For integer $(x,z)$ the $m_j$ and $r_{j,l}$ are integers and the same argument gives the integer statement. Every lattice point of the cone therefore has the decomposition
$(x,z)=\sum_{j,l}r_{j,l}a_{j,l}+(x-Bm,0),$
where $x-Bm\in C_B\cap\Lambda_B$ decomposes into vertices $b_j$ and indecomposable holes, and $(b_j\,,0)=\sum_la_{j,l}\,$; so the displayed set generates. Each $a_{j,l}$ has degree one and is indecomposable. If $(h,0)$ with $h\in\mathcal H_B$ were a sum of other generators, the private coordinates of the sum vanish, so within each block all stage coefficients are equal; if some $a_{j,1}$ occurred, then $b_j$ would occur in the $x$-part, so $h-b_j\in C_B\cap\Lambda_B$ would be nonzero, contradicting indecomposability of $h$; hence no stage generator occurs, $h$ is a sum of indecomposable holes, and indecomposability forces it to be a single one. The grading values are read off from the generators, and setting $z=0$ in the membership criterion gives $m=0$ and the clearing statement. Finally, $C^{(k)}$ is IDP exactly when $\mathcal H_B$ is empty.
\end{proof}

\begin{lemma}[Torus kernels determine induced lattices]\label{lem-toruskernel}
Let $R$ and $S$ be integer matrices with the same number $m$ of columns and with $\ker_\T R=\ker_\T S$, where $\ker_\T R=\{f\in\T^m:Rf=0\}$. Then $Ru\mapsto Su$ is a well-defined linear isomorphism between the real column spans which maps the induced lattice, the column cone and the column monoid of $R$ onto those of $S$. If all columns of both matrices have degree one, it preserves degree.
\end{lemma}
\begin{proof}
If $Rw=0$ for a real $w$, then every real multiple of $w$ reduces modulo $\Z^m$ to an element of $\ker_\T R=\ker_\T S$, so $aSw$ is integral for every real $a$ and $Sw=0$; by symmetry the real kernels agree and the map is a well-defined linear isomorphism. Moreover $Ru$ is integral exactly when the class of $u$ lies in $\ker_\T R$, which is exactly when $Su$ is integral, so the induced lattices correspond. Cones, monoids and degrees correspond by construction.
\end{proof}
For the series matrix of Lemma~\ref{lem-series}, the private equations over $\T$ force all stage coefficients in block $j$ to a common value $g_j\,$, so its torus kernel is $\{\operatorname{repeat}_k(g):g\in\ker_\T B\}$; appending free unit rays in new coordinates appends zero positions.

\subsection{Transfer to a face}
We now realize a series lift of $B$ on a face of three staircase simplices, by the chain presentation and column orders outlined in Subsection~\ref{sec-idp-overview}.
\begin{theorem}[IDP transfer to a face]\label{thm-idp-transfer}
Put $c:=r+1$. There are three staircase simplices in spatial dimension
\[
 d\le(9r+8)N+9r+11
\]
whose disjunctive hull has a face $F_{\mathcal Q}$ with the following properties, where $K_{F_{\mathcal Q}}\,$, $\Lambda_{F_{\mathcal Q}}=\spanR K_{F_{\mathcal Q}}\cap\Z^{d+3}$ and $S_{F_{\mathcal Q}}=\Z_{\ge0}(V\cap K_{F_{\mathcal Q}})$ are its cone, induced lattice and vertex monoid.
\begin{enumerate}
\item $K_{F_{\mathcal Q}}\cap\Lambda_{F_{\mathcal Q}}=S_{F_{\mathcal Q}}$ if and only if $\mathcal Q$ has the integer decomposition property.
\item There is an injective map $h\mapsto\widehat h$ from $\mathcal H_B$ into $\Lambda_{F_{\mathcal Q}}$ such that $\Hilb(K_{F_{\mathcal Q}}\,,\Lambda_{F_{\mathcal Q}})$ consists of the lifted vertices on $F_{\mathcal Q}$ together with the $\widehat h$; every $\widehat h$ belongs to $\Hilb(K)$, and $\mathcal C_K(\widehat h)=\mathcal C_B(h)$.
\item $\widehat h$ has multiplicities $(\lambda_0\,,\lambda_1\,,\lambda_2)=(2c\,t(h),2c\,t(h),2c\,t(h))$ and degree $6c\,t(h)$.
\item A primitive integral form exposes $F_{\mathcal Q}$ and takes only the values zero and one on $V$.
\item If $N\ge2$ and $\mathcal Q$ is full dimensional in $\R^r$, then $F_{\mathcal Q}$ is a facet and $d=(9r+8)N+r+11$.
\end{enumerate}
In particular, if the disjunctive hull has the integer decomposition property then so does $\mathcal Q$.
\end{theorem}
\begin{proof}
If $N=1$, take three copies of $[0,1]$ and their unimodular facet $x=0$; there are no indecomposable holes and every assertion follows. Hence assume $N\ge2$.

\emph{A chain presentation of $\ker_\T B$.} Let $\theta_\star$ be a zero node. Introduce the nodes $u_1\,,\ldots,u_{N-1}$ of a \emph{global chain}, and put $u_0=u_N:=\theta_\star\,$. For each spatial row $\ell\le r$ and $b\in\{0,1\}$, list $J_{\ell,b}:=\{j:v_{j,\ell}=b\}$ increasingly, put $n_{\ell,b}:=|J_{\ell,b}|$ and $m_{\ell,b}(j):=|J_{\ell,b}\cap\{1,\ldots,j\}|$, and introduce the nodes $s_{\ell,b,1}\,,\ldots,s_{\ell,b,n_{\ell,b}-1}$ of a chain for $J_{\ell,b}\,$, with $s_{\ell,b,0}=s_{\ell,b,n_{\ell,b}}:=\theta_\star$ (if $J_{\ell,b}$ is empty, there is no chain). Impose the gates
\[
 u_j=s_{\ell,0,m_{\ell,0}(j)}+s_{\ell,1,m_{\ell,1}(j)}\qquad(1\le\ell\le r,\ 1\le j<N),
\]
and add $z:=2r+2$ zero nodes, $\theta_\star$ among them, with their zero constraints. The same identities hold for $j=0$ and $j=N$, where both sides vanish. Subtracting consecutive identities shows that $g_j:=u_j-u_{j-1}$ is, for every $\ell$, the increment of the chain of the list $J_{\ell,b}$ containing $j$ at the position of $j$. Every chain starts and ends at zero, so $\sum_jg_j=0$ and $\sum_{j\in J_{\ell,1}}g_j=0$ for every $\ell$, that is, $Bg=0$. Conversely, every $g$ with $Bg=0$ determines all chain values by partial sums, and these satisfy the gates. Hence, over $\R$ and over $\T$, the scalar kernel is $\{g:Bg=0\}$, extended by the determined values. With $a_\ell\in\{1,2\}$ the number of nonempty lists $J_{\ell,b}\,$, the numbers of nodes and equations are
\[
 Q=N-1+\sum_{\ell=1}^r(N-a_\ell)+z,\qquad E=r(N-1)+z .
\]

\emph{The one-sign ordering.} Apply the extended incidence construction of Subsection~\ref{sec-extended} with the following column order in each class: for each chain with internal nodes $w_1\,,\ldots,w_{m-1}\,$, where $m\ge1$ is its number of increments, use the two lists $0,w_{m-1}\,,\ldots,w_1\,,0$ and $0,-w_1\,,\ldots,-w_{m-1}\,,0$, whose backward differences are its increments in reverse and in forward order. Every signed column of every internal node is used once. There are at most $4r+2$ lists, which need at most $4r+3$ zero separators, consecutive lists sharing a boundary zero, and the zero pool has $Q+2z\ge4r+4$ columns; concatenate with zero separators and append the unused zero columns. Every index $j$ lies in the global chain and in exactly one chain per spatial row, and each of these occurrences gives two positions per class whose coefficient character is $g_j$ on the kernel; for a chain with a single increment this character vanishes identically, and so does $g_j\,$. Every original character therefore occupies $2c$ positions per class and $6c$ in all, and all remaining coefficient characters vanish. Consequently, after grouping the coefficient positions by original column, the torus kernel of the endpoint matrix $R$ is
\[
 \ker_\T R=\{(\operatorname{repeat}_{6c}(g),0):g\in\ker_\T B\},
\]
the zeros occupying the $f=(9Q+3)-6cN$ remaining positions; here $f>0$, because $Q\ge(r+1)N+1$.

\emph{Identification and conclusion.} The series matrix of Lemma~\ref{lem-series} with $k=6c$, augmented by $f$ free unit rays in new coordinates, has the same number $9Q+3$ of columns and the same torus kernel. By Lemma~\ref{lem-toruskernel} the endpoint cone $K_{F_{\mathcal Q}}$ with its induced lattice and vertex monoid is isomorphic, degree-preservingly, to $C^{(k)}\times\R_{\ge0}^f$ with its induced lattice and the monoid generated by the stage generators and the free rays. Free unimodular rays change neither IDP nor the indecomposable holes, and only adjoin themselves to the Hilbert basis, so Lemma~\ref{lem-series} gives~(1), the Hilbert basis in~(2) with $\widehat h$ the image of $(h,0)$, and the clearing sets within $S_{F_{\mathcal Q}}\,$. By Lemma~\ref{lem-face}, exposure confines every ambient decomposition and every expression by $V$ of a point of $F_{\mathcal Q}$ to the face, so the $\widehat h$ are ambient Hilbert basis elements and $\mathcal C_K(\widehat h)=\mathcal C_B(h)$, which completes~(2). For~(3), $(h,0)$ has a representation with equal coefficient $p_j$ at all $k$ positions of block $j$ and zero on the free rays, where $Bp=h$; the class-$i$ multiplicity of $\widehat h$ is the sum of its coefficients on the $2c$ positions of each block in class $i$, namely $2c\sum_jp_j=2c\,t(h)$. Part~(4) is the exposing form of the incidence construction. Because $a_\ell\ge1$, we have $Q\le(r+1)N+r+1$ and $E=rN+r+2$, so $d=8Q+E+1\le(9r+8)N+9r+11$.

For~(5), no spatial row is constant, so $a_\ell=2$ for every $\ell$, $Q=(r+1)N+1$ and $d=(9r+8)N+r+11$. The real scalar kernel is $\ker_\R B$, of dimension $N-r-1$, so the endpoint matrix, with its $9Q+3$ columns, has rank $9Q+3-(N-r-1)=d+2$, using $E=Q-N+r+1$; the face it spans has codimension one in the $(d+3)$-dimensional ambient cone. Fixing the orders of ties and of the chains and the zero pool makes the three permutations an explicit function of $\mathcal Q$.
\end{proof}

\begin{example}[A tetrahedron]\label{ex-tetra}
The tetrahedron $\mathcal Q=\conv\{0,\mathbf e_1+\mathbf e_2\,,\mathbf e_1+\mathbf e_3\,,\mathbf e_2+\mathbf e_3\}\subseteq\R^3$ fails IDP: $h=(1,1,1;2)$, half the sum of the four columns $b_j\,$, is a lattice point of $\cone B$ that is not a sum of two columns, while $2h$ is the sum of all four, and $h$ is the only indecomposable hole. With $N=4$ and $r=3$, Theorem~\ref{thm-idp-transfer} transfers the IDP obstruction of $\mathcal Q$ to a facet $F_{\mathcal Q}$ of three staircase simplices in $\R^{154}$. The image $\widehat h$ is the only indecomposable hole of $F_{\mathcal Q}\,$; it lies in the ambient Hilbert basis, has multiplicities $(16,16,16)$ and degree $48$, and has clearing order $2$, like $h$.
\end{example}

IDP of the full hull implies IDP of the designated face.
Corollary~\ref{cor-words-unimodular} shows that a unimodular facet can coexist with indecomposable holes of the hull.

\section{Initial faces, the schedule criterion, and the complexity of splitting}\label{sec-ambient}
This section studies initial faces, the faces obtained by setting some coordinates to zero. The initial facet of a word disjunction (Section~\ref{sec-words}) is such a face, and so are the faces in the hardness constructions below. Subsections~\ref{sec-dichotomy} and~\ref{sec-twosat} give criteria for initial faces of coordinate availability at most two. Subsections~\ref{sec-splitting}--\ref{sec-holedeg} prove result~\hyperlink{res-E}{(E)}: the complexity of splitting a given lattice point, of recognizing IDP, and of the least degree of an indecomposable hole of a \emph{designated} facet, that is, a facet specified by an exposing form given as part of the input. Recall that the growth order $\gamma_i$ is the reverse of $\sigma_i\,$; we use multiplicity coordinates $(x,\lambda)$.

\subsection{Prefix chains and initial faces}\label{sec-dichotomy}
A \emph{prefix chain} on an ordered list $L=(l_1\,,\ldots,l_\ell)$ of distinct coordinates consists of the indicators of its $\ell+1$ initial segments, including the empty one. For lists $L_0\,,L_1\,,L_2$ let
\[
 K(L_0\,,L_1\,,L_2)=\cone\{(\mathbf e_{\{l_1,\ldots,l_b\}}\,,\eps_i):\ 0\le b\le|L_i|,\ i=0,1,2\}
\]
be the cone over the three chains, with its induced lattice, its vertex monoid $S(L)$ and its saturation $\overline S(L)$; for lists $L_0\,,\ldots,L_n$ define $K(L_0\,,\ldots,L_n)$ in the same way, with the unit vectors $\eps_0\,,\ldots,\eps_n\,$. The \emph{availability set} of a coordinate $k$ is the set of indices $i$ with $k\in L_i\,$. If $U\subseteq[d]$ and $L_i$ is the longest initial segment of $\gamma_i$ disjoint from $U$ (all of $\gamma_i$ if $U=\varnothing$), then the lifted vertices of the triple with $x_U=0$ are exactly the prefix indicators of the $L_i\,$, so $K(L_0\,,L_1\,,L_2)$ is the face $\{x_U=0\}$ of $K$, which we call an \emph{initial face}; by Lemma~\ref{lem-face} its Hilbert basis elements are ambient Hilbert basis elements with the same clearing sets. Two-chain subconfigurations are totally unimodular, by the proof of Proposition~\ref{prop-two}, so every lattice point of degree at most two in $K(L)$ lies in $S(L)$, and a degree-three point outside $S(L)$ is a Hilbert basis element.

\begin{theorem}[Clearing dichotomy for coordinates of availability at most two]\label{thm-dichotomy}
Let $L_0\,,\ldots,L_n$ be lists in which every coordinate has availability of size at most two, and let $S\subseteq\overline S$ be the vertex monoid and the saturation of the cone over their prefix chains. Then $2\overline S\subseteq S$, and for $y\in\overline S$,
\[
 \{q\in\Z_{>0}:qy\in S\}=\begin{cases}\Z_{>0}\,,&y\in S,\\ 2\Z_{>0}\,,&y\in\overline S\setminus S.\end{cases}
\]
\end{theorem}
\begin{proof}
Fix a lattice point $(x,\lambda)$ of the cone. A real expression by the chain vertices of list $i$ is the same as a nonincreasing \emph{cumulative profile} $\lambda_i\ge c_{i,1}\ge\cdots\ge c_{i,\ell_i}\ge0$, where $c_{i,r}$ is the total coefficient of the prefixes of length at least $r$, the coefficient of the prefix of length $r$ being $c_{i,r}-c_{i,r+1}$ (with $c_{i,0}=\lambda_i$ and $c_{i,\ell_i+1}=0$). The spatial equations are $c_{i,r}=x_k$ or $c_{i,r}+c_{j,s}=x_k\,$, one for each coordinate $k$, according as $k$ has availability one or two. The set of profiles is a nonempty polytope; choose a vertex of it, at which a full-rank set of constraints is active. The active constraints are spatial equations, consecutive equalities $c_{i,r}=c_{i,r+1}\,$, and the bounds $c_{i,1}=\lambda_i\,$, $c_{i,\ell_i}=0$, all with integer right-hand sides.

Form the graph whose nodes are the variables $c_{i,r}\,$, with an edge for each active equation of the form $u\pm v=b$, and call a node pinned if an active equation gives it an integer value. In a connected component, every variable equals $\pm u_0+b$ for a fixed variable $u_0$ and integers $b$. A pinned component is integral. An unpinned component must contain a cycle whose signs force $u_0=-u_0+b$, because otherwise it would retain a free real parameter, contradicting full rank; hence $u_0\in\frac12\Z$. Thus every variable at this vertex is half-integral. Doubling gives integral nonincreasing profiles for $(2x,2\lambda)$, whose consecutive differences are a nonnegative integer vertex expression. This proves $2\overline S\subseteq S$.

Conversely, suppose $q$ is odd and $qy\in S$, with integer profiles $c$ for $qy$. Then $c/q$ is a real profile for $y$ with all denominators odd, hence no strict half-integers. Round every entry to the nearest integer. Rounding is monotone, so the profiles stay nonincreasing, and the integer bounds are preserved. A singleton spatial equation is already integral. For $u+v=x_k\in\Z$ the fractional parts of $u,v$ are complementary and neither is $1/2$, so nearest-integer rounding preserves the sum. The rounded profiles are an integer vertex expression of $y$. Hence an odd multiple of $y$ lies in $S$ only when $y$ does.
\end{proof}

\begin{corollary}[Twofold clearing on initial faces]\label{cor-initial-two}
Let $F=\{x_U=0\}$ be an initial face of a triple of staircase simplices on which every coordinate has availability at most two. Then twice every lattice point of $K_F$ lies in the vertex monoid, and every Hilbert basis element of $K_F$ that is not a lifted vertex has ambient clearing set exactly $2\Z_{>0}\,$.
\end{corollary}
\begin{proof}
On an initial face every lifted vertex is a prefix of one of the lists $L_i\,$, so Theorem~\ref{thm-dichotomy} applies to $K_F$ with its induced lattice, and Lemma~\ref{lem-face} transfers Hilbert basis membership and clearing sets to $K$.
\end{proof}

\subsection{The schedule criterion}\label{sec-twosat}
Let $L_0\,,\ldots,L_n$ be lists in which every coordinate has availability of size at most two. Call a coordinate of availability $\{i\}$ a \emph{singleton of class $i$} and one of availability $\{i,j\}$ a \emph{pair coordinate}. The singletons of class $i$ cut $L_i$ into consecutive \emph{segments}, the maximal runs of pair coordinates between them. A \emph{schedule} assigns a bit $b_i(k)\in\{0,1\}$ to every class $i$ and every pair coordinate $k\in L_i\,$, nonincreasing along each segment of $L_i$; it is \emph{admissible} if $b_i(k)\ne b_j(k)$ at every pair coordinate $k$ of availability $\{i,j\}$.

\begin{theorem}[Schedule criterion]\label{thm-twosat}
The cone $K(L_0\,,\ldots,L_n)$ has the integer decomposition property in its induced lattice if and only if an admissible schedule exists. If none exists, the lattice point $y$ defined below is a lattice point of the cone outside the vertex monoid, and $2y$ lies in the vertex monoid. The existence of an admissible schedule is an instance of $2$-satisfiability, decidable in time linear in the total length of the lists.
\end{theorem}
\begin{proof}
Use cumulative profiles as in the proof of Theorem~\ref{thm-dichotomy}: a lattice point $(x,\lambda)$ of the cone has real profiles $u_i\,$, nonincreasing along $L_i$ with values in $[0,\lambda_i]$, with $u_i(k)=x_k$ at a singleton $k$ of class $i$ and $u_i(k)+u_j(k)=x_k$ at a pair coordinate of availability $\{i,j\}$; integer vertex expressions are exactly the integer profiles.

\emph{Sufficiency.} Given an admissible schedule, replace $u_i(k)$ by its ceiling where $b_i(k)=1$ and by its floor where $b_i(k)=0$, and leave the singleton values, which are integers. Within a segment the bits are nonincreasing, and $\lceil a\rceil\ge\lfloor a'\rfloor$ for $a\ge a'$, so the rounded profile is nonincreasing there. Across a singleton $k$ of class $i$ the profile passes the integer $x_k$: values before it are at least $x_k$ and remain so after rounding, values after it are at most $x_k$ and remain so. The bounds $0$ and $\lambda_i$ are integers and are preserved. At a pair coordinate the two profiles are rounded in opposite directions, and $\lfloor a\rfloor+\lceil a'\rceil=a+a'$ when $a+a'\in\Z$. So the rounded profiles are integral and give an integer vertex expression.

\emph{Necessity.} Let class $i$ have $\sigma_i$ singletons, so $\sigma_i+1$ segments, numbered $0,\ldots,\sigma_i$ in list order (some possibly empty). Define real profiles by $u_i=\sigma_i-t+\frac12$ on segment $t$, $u_i=\sigma_i-t+1$ at the $t$th singleton ($t=1,\ldots,\sigma_i$), and $\lambda_i=\sigma_i+1$. They are nonincreasing with values in $[0,\lambda_i]$, the singleton values are integers, and at a pair coordinate the sum of two half-integers is an integer. So they define a lattice point $y$ of the cone, with $x_k=0$ at every coordinate in no list. On the coordinates used by the lists, the induced lattice is the full lattice generated by the prefix chains, because each of them lies in some list; every other coordinate vanishes on the cone. In an integer profile for $y$, the values on segment $t$ of class $i$ lie between the integer values at its two ends, which differ by one, so they are $\sigma_i-t+b_i(k)$ with bits $b_i(k)$ nonincreasing along the segment; and at a pair coordinate the equation $b_i(k)+b_j(k)=1$ says that the bits differ. Thus integer vertex expressions of $y$ are exactly admissible schedules, and $y$ lies outside the vertex monoid when none exists. By Theorem~\ref{thm-dichotomy}, $2y$ is in the vertex monoid.

\emph{Complexity.} With one Boolean variable per bit, monotonicity along a segment is the set of implications $b_i(k')\Rightarrow b_i(k)$ for consecutive pair coordinates $k$ before $k'$ in a segment, and admissibility $b_i(k)\ne b_j(k)$ is the pair of clauses $b_i(k)\vee b_j(k)$ and $\neg b_i(k)\vee\neg b_j(k)$; all are $2$-clauses, and their number is linear in the total length.
\end{proof}
For three chains in a common order, the schedule condition becomes an interval condition.

\begin{corollary}[Interval criterion]\label{thm-interval}
Let three lists be restrictions of a common order, in which every nonsingleton coordinate has type $01$, $02$ or $12$ and every singleton coordinate comes after all pair coordinates. For a pair type $A$ let $I_A$ be the closed interval of positions between its first and its last occurrence, empty if $A$ does not occur. Then $K(L_0\,,L_1\,,L_2)$ has the integer decomposition property if and only if $I_{01}\cap I_{02}\cap I_{12}=\varnothing$. If the intersection is nonempty, the point $h$ with $x_k=1$ at every pair coordinate, $x_k=0$ at every singleton coordinate and $\lambda=(1,1,1)$ is a Hilbert basis element of degree three with exact clearing set $2\Z_{>0}\,$.
\end{corollary}
\begin{proof}
All pair coordinates of class $i$ precede its singletons, so a schedule is one cut per class, and its bits can be read along the common order as three nonincreasing bit sequences. If the intervals do not meet, two of them are disjoint by the Helly property of intervals on a line, say every occurrence of $B$ precedes every occurrence of $A$; cut between them, give bit $1$ to the two classes of $A$ before the cut and only to the class not in $B$ after it. These schedules are nonincreasing, and before (after) the cut the only pair whose bits agree is $A$ ($B$), which does not occur there. If a type is absent, a constant assignment whose agreeing pair is that type works. Conversely, in an admissible schedule the three-bit state at a pair position has one or two ones, and the pair whose bits agree does not occur at that position. If the state is constant, its agreeing pair never occurs, so its interval is empty and the intervals do not meet. Otherwise monotonicity allows exactly one change, from a two-one state with agreeing pair $A$ to a one-one state with agreeing pair $B\ne A$; then every occurrence of $B$ precedes every occurrence of $A$, so the intervals do not meet. Finally, $h$ is half the sum of the empty prefix and the longest pair prefix in each class, so it is a lattice point of the cone; its integer vertex expressions are binary profiles vanishing at the singletons, that is, admissible schedules. So $h$ is outside the vertex monoid exactly when the intervals meet; in that case it is a Hilbert basis element because it has degree three, and its clearing set is $2\Z_{>0}$ by Theorem~\ref{thm-dichotomy}.
\end{proof}

\subsection{Splitting a specified point}\label{sec-splitting}
On an initial face $K(L_0\,,L_1\,,L_2)$ (Subsection~\ref{sec-dichotomy}), a real expression of $(x,\lambda)$ is a triple of nonincreasing profiles on the lists with values in $[0,\lambda_i]$, whose entries sum to $x_j$ at each coordinate; a coordinate of availability one fixes a profile entry and one of availability two fixes the sum of two entries.

\begin{theorem}[Availability at most two]\label{thm-split-two}
If every coordinate of $L_0\,,L_1\,,L_2$ has availability at most two, whether a given lattice point splits is decided by one linear program followed by a $2$-SAT instance of linear size.
\end{theorem}
\begin{proof}
Find a real profile $f$ by linear programming; if none exists the point is not in the cone. Restrict each entry to $\lfloor f\rfloor$ or $\lceil f\rceil$. If $p$ is any integer profile, clip it coordinatewise to $[\lfloor f\rfloor,\lceil f\rceil]$. The clipped profiles stay nonincreasing, because the profiles and both sequences of interval endpoints are nonincreasing. At a pair coordinate the two intervals are reflections of each other under $u\mapsto x_j-u$, so clipping preserves the sum; singleton entries and the integer bounds are preserved. Hence the restricted choices contain an integer solution whenever one exists. Forbidding, for each consecutive order inequality and each pair equation, the violating choices of its at most two rounding bits gives linearly many $2$-clauses.
\end{proof}

\begin{theorem}[Hardness at availability three]\label{thm-split-hard}
Deciding whether a lattice point splits is strongly NP-complete, even on a facet $\{x_g=0\}$ of three staircase simplices, for $h=(x,m,m,m)$ with $m\le d$ and $0\le x_j\le3m$, and even when integral vertex expressions of $2h$ and $3h$ are given. In the reduction, $\mathcal C_K(h)$ is $\Z_{>0}$ or $\Z_{\ge2}\,$, and $h$ lies in the lattice generated by the rays of its minimal face.
\end{theorem}
\begin{proof}
Reduce from positive one-in-three SAT~\cite{Schaefer}. Use variables $X_i$ for the Boolean variables and $y,u,v,p_2\,,q_2\,,q_3\,,w\in[0,1]$ with
\[
 p_2=u+u,\quad q_2=v+v,\quad q_3=q_2+v,\quad w=p_2+q_3\,,\quad y+w=1,
\]
and, for each clause $(i,j,\ell)$, a fresh $z\in[0,1]$ with $z=X_i+X_j$ and $z+X_\ell=y$. At an integer solution $2u,2v\le1$ force $u=v=0$, hence $y=1$, and the clause rows are the one-in-three conditions; every satisfying assignment extends. Independently of the formula, setting all $X_i\,$, clause variables and $y$ to zero, the choices $(u,v)=(\frac12,0)$ and $(0,\frac13)$ give solutions that are integral after doubling and tripling. Apply Lemma~\ref{lem-bounded}: these give integral expressions of $2h$ and $3h$, and $h$ splits exactly when the formula is satisfiable. All sizes are linear, $m=2n\le d$, and feasibility gives $0\le x_j\le3m$. Membership in NP is witnessed by the integer coefficients of the $3(d+1)$ lifted vertices. Because $2$ and $3$ generate every integer at least two, the clearing set is as stated; both supplied expressions lie in the minimal face of $h$, and their difference puts $h$ in its ray lattice.
\end{proof}

\begin{theorem}[Two parameters]\label{thm-split-param}
For a lattice point $h=(x,\lambda)$ on an initial face of three staircase simplices, let $q$ be the number of coordinates of availability three, $k=\min_i\lambda_i\,$, and $\ell$ the input length. Splitting $h$ takes time $2^{O(q\log(q+2))}\operatorname{poly}(\ell)$ and $(d+1)^k\operatorname{poly}(\ell)$. It is W[1]-hard with parameter $k$. The hardness holds on a face $\{x_g=0\}$ with integral expressions of $2h$ and $3h$ given and all multiplicities and coordinates $O(d)$.
\end{theorem}
\begin{proof}
\emph{Integrality only at triple coordinates.} Delete the $q$ triple coordinates from the lists and test the resulting availability-two instance by Theorem~\ref{thm-split-two}; if it fails, $h$ does not split, and otherwise fix an integer profile $z$ for it. If $f$ is a real profile on the full lists that is integral at the triple coordinates, keep these entries and clip $z$ to $[\lfloor f\rfloor,\lceil f\rceil]$ elsewhere. The reduced constraints hold as before, and an entry before (after) an integral triple entry has floor at least (ceiling at most) that entry, so the full lists stay nonincreasing. Hence $h$ splits if and only if some real profile is integral at the triple coordinates.

\emph{The parameter $q$.} Substitute the fixed singleton entries and eliminate one entry of each pair equation. Apart from the $q$ triple equations, all constraints have at most two variables with coefficients $\pm1$. Eliminate the remaining non-triple variables by Fourier--Motzkin; combining two such constraints again has at most two signed terms (a repeated variable is divided by two), only the tightest bound per sign pattern is kept, and the bit lengths stay polynomial, growing linearly in the number of eliminations. The result is an integer program in $3q$ variables, solved in time $2^{O(q\log(q+2))}\operatorname{poly}(\ell)$~\cite{Kannan}; a solution extends by linear programming and the clipping above.

\emph{The parameter $k$.} Enumerate the at most $(d+1)^k$ multisets of $k$ prefixes of a list of multiplicity $k$, subtract, and test real feasibility of the remaining two lists, which suffices because rounding one profile down and the other up preserves integral pair sums, order and bounds.

\emph{W[1]-hardness.} Reduce from Multicolored Clique~\cite{FHRV} with color classes of sizes $n_1\,,\ldots,n_k\,$, after tripling every vertex. Encode thresholds $1=X_{i,1}\ge\cdots\ge X_{i,n_i+1}=0$, indicators $a_{iu}=X_{i,u}-X_{i,u+1}$ and complements $b_{iu}=1-a_{iu}\,$, and impose $b_{jv}\ge a_{iu}$ for every non-edge between options $(i,u)$ and $(j,v)$, $i<j$; integral thresholds select one option per class, forming a clique exactly when all these inequalities hold. List $A$ has one block of unit width per class, carrying the thresholds $X_{i,u}$ on consecutive marked subintervals; list $C$ has a companion block for each scalar, tied to every subinterval using it by early and late marker coordinates as in Lemma~\ref{lem-bounded}; list $B$ has definition blocks for $a_{iu}$ and $b_{iu}\,$, each containing one coordinate of availability three imposing $a_{iu}+X_{i,u+1}+(1-X_{i,u})=1$ or $b_{iu}+X_{i,u}+(1-X_{i,u+1})=2$, and one block per non-edge with two marked subintervals at the same offset, tied to $b_{jv}$ and $a_{iu}\,$, so that monotonicity imposes $b_{jv}\ge a_{iu}\,$. Then the multiplicity of $A$ is $k$, there are $q=2n$ coordinates of availability three, and $d=O(n^2)$. Giving weight $1/p$ to $p$ distinct options per class, $p=2,3$, satisfies all inequalities and yields expressions of $2h$ and $3h$. Multicolored Clique is W[1]-hard~\cite{FHRV}.
\end{proof}

We collect these results.

\begin{theorem}[Splitting by availability]\label{thm-complexity}
Let $h=(x,\lambda)$ be a lattice point of $K(L)$.
\begin{enumerate}
\item If every availability has size at most two, then $2h$ splits, $h$ splits as soon as an odd multiple does, and every indecomposable hole has clearing set exactly $2\Z_{>0}\,$. Whether $h$ splits is decided by a linear program and a linear-size $2$-SAT instance, and IDP of $K(L)$ is a linear-time $2$-SAT test.
\item In general, deciding whether $h$ splits is strongly NP-complete, even if $K(L)$ is a facet $\{x_g=0\}$, $h=(x,m,m,m)$ with $m\le d$ and $0\le x\le3m$, and integral vertex expressions of $2h$ and $3h$ are given.
\item Splitting is fixed-parameter tractable in the number $q$ of availability-three coordinates; in $k:=\min_i\lambda_i$ it takes time $(d+1)^k\operatorname{poly}(\ell)$, $\ell$ the input length, yet is W[1]-hard.
\end{enumerate}
\end{theorem}
\begin{proof}
Part~(1) combines Theorems~\ref{thm-dichotomy}, \ref{thm-twosat} and~\ref{thm-split-two}; part~(2) is Theorem~\ref{thm-split-hard}; part~(3) is Theorem~\ref{thm-split-param}.
\end{proof}

\subsection{Recognizing the integer decomposition property}\label{sec-recognition}
Certificates of degree below the rank place the recognition of IDP in coNP. Hardness comes from edge polytopes, transferred to a designated facet by Theorem~\ref{thm-idp-transfer}.
\begin{lemma}[Certificates]\label{lem-cert}
Let $C$ be a pointed rational cone, graded by a linear form that is integral on its lattice, whose degree-one lattice points $a_1\,,\ldots,a_N$ are listed explicitly and generate it. If some lattice point of $C$ is not a nonnegative integer combination of the $a_j\,$, there is one, $h$, of degree at most $\operatorname{rank}C-1$ with $h-a_j\notin C$ for all $j$. Such an $h$ is verified by linear programming.
\end{lemma}
\begin{proof}
A nonsplitting point of least degree admits no subtraction $h-a_j\in C$. In an expression of $h$ on independent generators every coefficient is below one, and the coefficients sum to the integer degree, which is therefore at most $\operatorname{rank}C-1$. A point admitting no subtraction does not split.
\end{proof}

\begin{corollary}[Recognizing IDP]\label{cor-recognition}
Deciding whether a $0/1$ polytope given by its vertices has IDP is coNP-complete, already for edge polytopes of connected graphs. Deciding whether the facet of three staircase simplices exposed by a given primitive integral form with values in $\{0,1\}$ on $V$ has IDP is coNP-complete, even for instances in which twice every lattice point of the facet cone splits.
\end{corollary}
\begin{proof}
The lattice points of a $0/1$ polytope are its vertices, and those of a facet of $\mathcal D$ are its lifted vertices, so Lemma~\ref{lem-cert} gives membership in coNP. For hardness, let $H$ be a graph. If $H$ is bipartite, it has no odd cycle and we output a fixed IDP instance. Otherwise join a new vertex by one edge to one vertex of each component; the new edges are bridges, so the odd cycles of $H$ and the edges between them are unchanged, and the resulting graph $G$ is connected and not bipartite. Two vertex-disjoint odd cycles joined by no edge may be shortened along chords to chordless ones on subsets of their vertex sets, so $G$ violates the odd cycle condition of Ohsugi and Hibi exactly when $H$ has two vertex-disjoint odd cycles with no edge between them, which is NP-complete to decide~\cite{GKPT}. By~\cite{OH}, the edge polytope $\mathcal Q_G=\conv\{\mathbf e_u+\mathbf e_v:uv\in E(G)\}$ has IDP exactly when the odd cycle condition holds; because $G$ is connected and not bipartite, the edge vectors generate every lattice point of even coordinate sum, so normality in~\cite{OH} is IDP in $\Z^{V(G)}$. Deleting the coordinate of one vertex $w$ maps the hyperplane $\sum_vx_v=2t$ isomorphically, with its lattice, onto $\R^{V(G)-w}\times\R$, so it yields a full-dimensional $0/1$ polytope with distinct vertices and the same IDP status. A vertex of $\{c\ge0:\sum_{e\ni v}c_e=x_v\}$ has a support whose components contain at most one cycle, which is odd, so it is half-integral; hence twice every lattice point of the cone over $\mathcal Q_G$ splits. Theorem~\ref{thm-idp-transfer}(5) realizes this polytope on a facet in polynomial size, with a primitive $0/1$ exposing form, preserving IDP and all clearing sets.
\end{proof}

\subsection{The least indecomposable hole degree of a designated facet}\label{sec-holedeg}
For a simplex face $F$ let $\nu(F)$ be the least degree of a lattice point of $K_F$ outside its vertex monoid, or $\infty$. A lattice point of $K_F$ is a nonnegative integral ray sum plus a box point, and it splits exactly when the box point is zero; so $\nu(F)$ is the least degree of a nonzero box point.

\begin{theorem}[Hardness of the least indecomposable hole degree]\label{thm-holedeg-hard}
Given three staircase simplices, a primitive integral form with values in $\{0,1\}$ on $V$ exposing a simplex facet $F$, and an integer $T$, deciding whether $\nu(F)\le T$ is strongly NP-complete, even when $G_F\cong(\Z/2\Z)^k$ and every nonzero box point has clearing set $2\Z_{>0}\,$.
\end{theorem}
\begin{proof}
For membership, check that the exposed rays are independent and guess a lattice point of $K_F$ of degree at most $T$; its ray coefficients have polynomial size, and it lies outside the vertex monoid exactly when one of them is not an integer.

We reduce from the minimum distance of binary linear codes, which is NP-complete~\cite{Vardy} (stated there for a parity-check matrix, from which elimination gives a generator matrix): given a $k\times n$ generator matrix $\Gamma$ over $\F_2$ of full row rank and $t_0\le n$, decide whether some nonzero codeword has weight at most $t_0\,$. Build a scalar presentation with a zero node $\zeta=0$, input nodes $a_1\,,\ldots,a_k$ with $a_i+a_i=0$, and addition nodes computing every code coordinate $c_j=\sum_i\Gamma_{ij}a_i\,$, with $c_j=\zeta+\zeta$ for a zero column. Let $S_0\le1+k+nk$ be the number of nodes so far, put $R=S_0+1$, and append $R$ copies $c_j+\zeta$ of every $c_j\,$. The presentation is square, its real kernel is zero, and its torus kernel is $\{a\in(\frac12\Z/\Z)^k\}\cong(\Z/2\Z)^k$, every node value lying in $\{0,\frac12\}$. Apply Lemma~\ref{univ-incidence} and Proposition~\ref{univ-facet}. At a node of value $\frac12$, the coefficient blocks are $(\frac12,0,\frac12)$ in all three colors, so a box point has degree three times its number of active nodes, namely
\[
 \deg H_a=3\bigl(R\,\wt(c(a))+s_0(a)\bigr),\qquad 0\le s_0(a)\le S_0\,.
\]
Put $T=3(Rt_0+S_0)$. If $\wt(c(a))\le t_0$ then $\deg H_a\le T$, and if $\wt(c(a))\ge t_0+1$ then $\deg H_a\ge3R(t_0+1)>T$. Because $\Gamma$ has full row rank, $a\ne0$ exactly when $c(a)\ne0$. So $\nu(F)\le T$ if and only if the code has a nonzero word of weight at most $t_0\,$. All numbers are polynomial in $n$ and $k$, and the clearing sets follow from Lemma~\ref{univ-incidence}.
\end{proof}

\section{Where integer decomposition fails}\label{sec-failures}
This section proves result~\hyperlink{res-F}{(F)}. The designated faces of Sections~\ref{sec-cyclic}--\ref{sec-idp} isolate arithmetic on one face. Here we turn to the whole hull and ask how and where IDP can fail. IDP can first fail in large degree, so that testing lattice points of low degree cannot certify it, and the certificates of degree below the rank in Lemma~\ref{lem-cert} are tight up to a constant factor (Subsection~\ref{sec-firsthole}). A multiple of a lattice point that lies in the ray lattice of its minimal face can fail to split for exponentially long (Subsection~\ref{sec-delay}). Six dimensions are needed for any failure (Subsection~\ref{sec-dim5}), and an indecomposable hole need not lie in any simplex face (Subsection~\ref{sec-offface}).

\subsection{The degree of the first indecomposable hole}\label{sec-firsthole}
The indecomposable holes constructed in Section~\ref{sec-cyclic} have small degree. Here the least degree of a nonsplitting lattice point grows linearly with the dimension.
Write $\nu(\mathcal D)$ for the least degree of a lattice point of $K$ outside the vertex monoid, with $\nu(\mathcal D)=\infty$ under IDP. By Proposition~\ref{prop-two}, $\nu(\mathcal D)\ge3$, and by Proposition~\ref{prop-bounds}, $\nu(\mathcal D)\le d$ when IDP fails.

\begin{theorem}[First indecomposable holes of large degree]\label{thm-firsthole}
For integers $m\ge1$ and $\kappa\in\{0,1\}$ there are three staircase simplices in spatial dimension $d=15m+16+11\kappa$ with $\nu(\mathcal D)=2m+1+\kappa$. An indecomposable hole of this degree can be taken with multiplicities $(1,1,2m+\kappa-1)$ and clearing set $2\Z_{>0}\,$, at height one above a facet that has IDP and from which every off-facet vertex has height one; and twice every lattice point of $K$ splits. In particular, for every $k\ge3$ some triple in dimension at most $8k+12$ has $\nu(\mathcal D)=k$.
\end{theorem}

\par
\begin{figure}[t]
\centering
\begin{tikzpicture}[x=0.95cm,y=0.95cm,font=\small]
 \tikzset{ov/.style={circle,fill,inner sep=1.6pt},pv/.style={rectangle,draw,fill=white,inner sep=2.6pt}}
 \node[ov] (a1) at (0,0) {}; \node[ov] (a2) at (-0.9,0.55) {}; \node[ov] (a3) at (-0.9,-0.55) {};
 \node[ov] (b1) at (1.4,0) {}; \node[ov] (b2) at (2.3,0.55) {}; \node[ov] (b3) at (2.3,-0.55) {};
 \draw[thick] (a1)--(a2)--(a3)--(a1); \draw[thick] (b1)--(b2)--(b3)--(b1);
 \draw[thick,blue!70!black] (a1)--(b1);
 \node at (-0.6,0) {$Z$}; \node at (2.0,0) {$Z'$};
 \node[font=\scriptsize,align=center] at (0.7,-1.25) {(a) alternating path of length one:\\ $b_Z+b_{Z'}$ splits};
 \begin{scope}[xshift=6.2cm]
 \node[ov] (c1) at (0,0) {}; \node[ov] (c2) at (-0.9,0.55) {}; \node[ov] (c3) at (-0.9,-0.55) {};
 \node[ov] (d1) at (2.4,0) {}; \node[ov] (d2) at (3.3,0.55) {}; \node[ov] (d3) at (3.3,-0.55) {};
 \node[pv] (p) at (1.2,0.75) {};
 \draw[thick] (c1)--(c2)--(c3)--(c1); \draw[thick] (d1)--(d2)--(d3)--(d1);
 \draw[thick,dashed] (c1)--(p)--(d1);
 \node[above,font=\scriptsize] at (p.north) {pivot};
 \node at (-0.6,0) {$Z$}; \node at (3.0,0) {$Z'$};
 \node[font=\scriptsize,align=center] at (1.2,-1.25) {(b) the only path passes a pivot:\\ $b_Z+b_{Z'}$ does not split; $2(b_Z+b_{Z'})$ does};
 \end{scope}
\end{tikzpicture}
\caption{Lemma~\ref{lem-signed-odd}. Ordinary vertices (dots) carry negative edges (solid) and pivots (squares) carry positive edges (dashed). \emph{(a)} A single negative edge is an alternating path between the odd cycles $Z$ and $Z'$, so $b_Z+b_{Z'}$ splits. \emph{(b)} The only path from $Z$ to $Z'$ passes through a pivot, whose two edges are both positive, so it is not alternating; $b_Z+b_{Z'}$ is a lattice point of degree three that does not split and splits exactly at its even multiples.}
\label{fig-signed}
\end{figure}
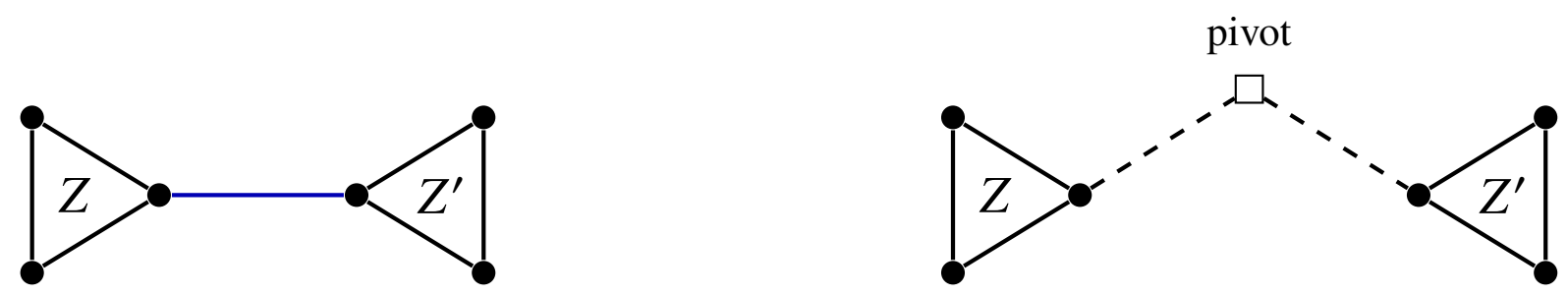

\medskip\noindent\emph{A signed edge model.} Let $\Omega$ (ordinary vertices) and $P$ (pivots) be disjoint, and take a simple graph with \emph{negative} edges inside $\Omega$ and \emph{positive} edges between $\Omega$ and $P$. For an edge $e=uv$ of sign $s_e$ put $r_e=s_e(\mathbf e_u+\mathbf e_v)$, and let $\Lambda=\{y\in\Z^{\Omega\cup P}:\sum_uy_u\in2\Z\}$. All $r_e$ have degree one for $\deg y=\frac12\bigl(3\sum_{p\in P}y_p-\sum_{w\in \Omega}y_w\bigr)$. For a simple odd cycle $Z$ put $b_Z=\frac12\sum_{e\in E(Z)}r_e\,$, an integral vector with odd coordinate sum. An \emph{alternating path} between two disjoint cycles has its ends on the two cycles, internal vertices outside them, and successive edges of opposite signs; a single edge counts (Figure~\ref{fig-signed}).

\begin{lemma}[Signed odd cycles; cf.~\cite{BL}]\label{lem-signed-odd}
The lattice points of $\cone\{r_e\}$ in $\Lambda$ are generated by the $r_e$ and the points $b_Z+b_{Z'}$ for vertex-disjoint odd cycles with no alternating path between them. Such a point has degree $(|Z|+|Z'|)/2$, does not split, and splits exactly at its even multiples; if an alternating path exists, $b_Z+b_{Z'}$ splits. Hence the least degree of a nonsplitting lattice point is the least $(|Z|+|Z'|)/2$ over such pairs, and twice every lattice point splits.
\end{lemma}
\begin{proof}
Take an expression on independent rays. In the unsigned support graph every component has at most as many edges as vertices and no even cycle, so it is a tree or an odd unicyclic graph; leaf elimination makes the tree coefficients integral, and on each odd cycle all coefficients are integral or all lie in $\frac12+\Z$. After removing integer parts, the remainder is a sum of $b_Z$ over vertex-disjoint odd cycles, evenly many because the coordinate sum is even; pair them arbitrarily. This also shows that twice every lattice point splits. On an odd cycle one can move the half coefficients to $0$ or $1$, alternating around the cycle, so that the signed incidence vector changes by $\pm\mathbf e_u$ at a chosen vertex $u$. Along an alternating path the sum of the path rays vanishes at internal vertices, so perturbing the two cycles at the ends of the path expresses $b_Z+b_{Z'}$ integrally. Conversely, suppose an odd multiple of $b_Z+b_{Z'}$ has an integral expression. At each vertex outside the cycles the target is zero, so positive and negative incidences pair off, and the edge copies decompose into alternating walks ending on the cycles and closed walks. The target has odd coordinate sum on $V(Z)$, so some walk joins the two cycles. A pivot has only positive edges and cannot be internal to an alternating walk, so the walk contains an alternating path. The even multiples are split by the displayed half-sums.
\end{proof}

\par\medskip\noindent\emph{The construction.} Put $b=m+\kappa$, $n=4m+4+3\kappa$, $a=n-b$ and $r=b+1$. Use ordinary tracks $A_0\,,\ldots,A_a\,$, $B_0\,,\ldots,B_n\,$, $C_0\,,\ldots,C_n$ and pivots $p_0\,,\ldots,p_{n-1}\,$. For each pair of tracks, take a monotone shuffle starting at index zero in both, and join the current pair of vertices by a negative edge at every visited state, including the initial one; the shuffles are $B^nA^a$ for $AB$, $C^nA^a$ for $AC$, and $B^bT_1\cdots T_{n-b}C^b$ for $BC$, a letter meaning that its track advances, where each $T_j$ is $BC$ except $T_r=CB$ when $\kappa=1$. The positive neighbors of $p_j$ are $B_j\,,B_{j+1}\,,C_j\,,C_{j+1}$ and $A_i$ for $\alpha_j\le i\le\alpha_{j+1}\,$, with $\alpha_j=\max(0,j-b)$; in each track these intervals cover the track and successive ones meet in one vertex.

Index the constituents by $A,B,C$. Use one spatial coordinate per ordinary increment $X_{i-1}\to X_i\,$, one per pivot increment $p_j\to p_{j-1}\,$, and a guard $g$, in total $a+3n=d$. The growth order of constituent $X$ first follows the shuffle of the other two tracks, then reads $g$, and then follows its own track and the pivots in reverse from $X_{\rm last}$ and $p_{n-1}$: at $p_j$ it reads the own-track increments backwards down to the lower end of the interval of $p_j\,$, and then the increment $p_j\to p_{j-1}$ if $j>0$. Map an ordinary increment to $\mathbf e_{X_{i-1}}-\mathbf e_{X_i}\,$, a pivot increment to $\mathbf e_{p_{j-1}}-\mathbf e_{p_j}\,$, the guard to $\mathbf e_{A_a}+\mathbf e_{B_n}+\mathbf e_{C_n}+\mathbf e_{p_{n-1}}\,$, and the selector of $X$ to the negative ray of the edge joining the initial vertices of the other two tracks. An increment of a negative endpoint and a reverse increment of a positive endpoint have the same image, and the guard turns the last negative ray of a constituent into its first positive ray; so the prefixes map bijectively onto the signed rays. The increments generate the kernel of the map summing the coordinates of each track and of the pivots, and the images of the three selectors and the guard there, $(0,-1,-1,0)$, $(-1,0,-1,0)$, $(-1,-1,0,0)$, $(1,1,1,1)$, have determinant $\pm2$ and generate the even-sum lattice. Hence this linear map is a lattice isomorphism from $(K,\Z^{d+3})$ onto $(\cone\{r_e\},\Lambda)$ carrying the vertex monoid to the ray monoid.

\begin{lemma}[Lower bound]\label{lem-firsthole-lower}
Two vertex-disjoint odd cycles with no alternating path between them have total length at least $2(2m+1+\kappa)$.
\end{lemma}
\begin{proof}
Every odd cycle meets $B\cup C$, because the graph on $A\cup P$ is bipartite. A cycle through $A_0$ is joined by a negative edge to any $B$- or $C$-vertex of the other cycle, and a cycle through $p_j$ with $j\le b$ reaches the other cycle through the positive edge $p_jA_0$ and a negative edge from $A_0\,$. So both cycles avoid $A_0\,,p_0\,,\ldots,p_b\,$. They also avoid $B_n\,,C_n\,$: if both meet $\{B_n\,,C_n\}$ the edge $B_nC_n$ joins them; if only one does, the other contains a pivot, because the negative graph without $A_0\,,B_n\,,C_n$ is bipartite, and that pivot has a positive edge to a vertex of $A$, which is joined negatively to both $B_n$ and $C_n\,$. After these deletions $B_0\,,\ldots,B_b$ and $C_0$ lie on no cycle, and every remaining negative edge is $B_iC_j$ with $i-j\in\{b,b+1\}$, except for the edge $B_{r+b-1}C_r$ with difference $b-1$ when $\kappa=1$, which we call the \emph{exceptional edge}. Label $\tau(B_i)=\tau(C_i)=i$, $\tau(p_j)=j+\frac12$, $\tau(A_i)=b+i$. Along a positive edge $\tau$ changes by at most $\frac12$, and along an oriented negative edge by $\varepsilon(b+\frac12)+\eta$ with $\varepsilon=\pm1$ and $|\eta|\le\frac12$, or $|\eta|\le\frac32$ on the exceptional edge. An odd cycle $Z$ has evenly many positive edges, hence oddly many negative ones, and summing the changes around $Z$ gives $0=(b+\frac12)\sum\varepsilon_e+\sum\eta_e$ with $\sum\varepsilon_e$ odd. Thus $|Z|/2\ge b+\frac12-e(Z)$, where $e(Z)\in\{0,1\}$ records the exceptional edge, which at most one of the two cycles uses. Adding gives $(|Z|+|Z'|)/2\ge2b+1-\kappa=2m+1+\kappa$.
\end{proof}

\begin{lemma}[Upper bound]\label{lem-firsthole-upper}
For $(u,q)=(r,m)$ and $(u,q)=(s,b)$ with $s=3m+3+2\kappa$, the cycle $Z_u$ following $C_u\,,p_u\,,C_{u+1}\,,p_{u+1}\,,\ldots,C_{u+q-1}\,,p_{u+q-1}\,,B_{u+q}$ and closing with the negative edge $B_{u+q}C_u$ is odd of length $2q+1$, and there is no alternating path between $Z_r$ and $Z_s\,$.
\end{lemma}
\begin{proof}
The closing edges exist: for $u=r$ and $\kappa=1$ it is the exceptional edge, and otherwise it is visited right after $T_u\,$. The ordinary vertices of each cycle and the ordinary neighbors of its pivots have $B,C$ indices in $[u,u+q]$ and $A$ indices in $[u-b,u+q-b]$; these sets are disjoint for the two cycles, so no positive edge joins a pivot of one cycle to an ordinary vertex of the other. A negative edge between the two sets would be a $BC$ edge, as neither set contains $A_0\,,B_n\,,C_n\,$; but the $BC$ index intervals are separated by $s-(r+m)=b+2$, while every $BC$ edge has index difference at most $b+1$. An alternating path has length at most three, with only its end edges positive, so it would need one of these joins.
\end{proof}

\begin{proof}[Proof of Theorem~\ref{thm-firsthole}]
Lemmas~\ref{lem-signed-odd}--\ref{lem-firsthole-upper} and the lattice isomorphism give $\nu(\mathcal D)=2m+1+\kappa$ and that twice every lattice point splits. Let $h$ be the preimage of $b_{Z_r}+b_{Z_s}\,$. In each cycle the negative edge lies in constituent $A$, the last positive edge in $B$ and the others in $C$, which gives the multiplicities. It is indecomposable because all lattice points of smaller degree split, and Lemma~\ref{lem-signed-odd} gives its clearing set. The form $\ell(y)=\frac12\bigl(\sum_Py_p-\sum_Wy_w\bigr)$ is integral on $\Lambda$, vanishes on positive rays and equals one on negative rays; it corresponds to $t-x_g\,$. The positive graph is connected and bipartite, so its rays span a hyperplane, and $\ell=0$ is a facet; it has IDP because its graph has no odd cycle. Each cycle contributes half of the slack of its negative edge, so $\ell(h)=1$. Finally, $k=2m+1+\kappa$ with $d=15m+16+11\kappa\le8k+12$.
\end{proof}

\begin{proposition}[Doubling caps the first indecomposable hole]\label{prop-double}
If $h$ is an indecomposable hole of three staircase simplices in $\R^d$ and $2h\in S_K\,$, then $\deg h\le\lfloor(d+3)/2\rfloor$. Hence if twice every lattice point of $K$ splits and IDP fails, then $\nu(\mathcal D)\le\lfloor(d+3)/2\rfloor$; with Theorem~\ref{thm-firsthole}, the largest $\nu(\mathcal D)$ over such triples in $\R^d$ is $\Theta(d)$.
\end{proposition}
\begin{proof}
Write $2h$ as a sum of $2\deg h$ lifted vertices. If a vertex $v$ occurred twice, then $h-v$ would be a lattice point of $K$, and $h=v+(h-v)$ would decompose. So $h=\frac12\sum_{v\in U}v$ for a set $U$ of $2\deg h$ distinct vertices. If $U$ were dependent, a real relation $\sum_U\mu_vv=0$ scaled so that $\max_v|\mu_v|=\frac12=\mu_w$ for some $w$ would give $h=\sum_U(\frac12+\mu_v)v$ with nonnegative coefficients and coefficient one at $w$, and again $h-w\in K$. Hence $2\deg h=|U|\le d+3$. A least-degree nonsplitting lattice point is indecomposable, because a decomposition would have summands of smaller degree, which split.
\end{proof}

\subsection{How late a multiple can fail to split}\label{sec-delay}
For a nonzero lattice point $h$ of $K$, let $F$ be its minimal face, $A$ the matrix of the rays of $F$, and $q$ the order of $h$ modulo $\Z A$. Only the multiples of $q$ can clear $h$, and we call them \emph{lattice-compatible}. Let $L(d)$ be the largest lattice-compatible $k$ with $kh\notin S_K\,$, over all lattice points $h$ of all three-staircase cones in $\R^d$.

\begin{theorem}[Clearing delay]\label{thm-delay}
Let $\rho$ be the plastic number, $\rho^3=\rho+1$. There is $c>0$ such that for all sufficiently large $d$
\[
 c\rho^{2d}\le L(d)\le(2d+1)B(d-1)+B(d-1)^2 .
\]
In particular $L(d)=2^{\Theta(d)}$, and the lower bound is attained by indecomposable holes of degree three that are trivial modulo the ray lattice of their minimal face.
\end{theorem}

\begin{lemma}[Conductor bound]\label{lem-conductor}
In the notation above, every lattice-compatible $k\ge(2d+1)B(d-1)+B(d-1)^2$ satisfies $kh\in S_K\,$.
\end{lemma}
\begin{proof}
By exposure (Lemma~\ref{lem-face}), $kh\in S_K$ exactly when $kh$ is a nonnegative integer combination of the columns of $A$. Let $P_h=\{c\ge0:Ac=h\}$ and $Q_h=\{(c,k)\ge0:Ac=kh\}$. The polytope $P_h$ is bounded, because the coordinates of $c$ sum to $\deg h$, and it has a strictly positive point, because $h$ lies in the relative interior of $F$; put $s=1+\dim P_h\,$. The heights $k$ of the lattice points of $Q_h$ form $\{0\}\cup\mathcal C_K(h)$. A vertex $u$ of $P_h$ gives the lattice point $(D_uu,D_u)$ of $Q_h\,$, where $D_u$ is the least common denominator of $u$; let $D_h=\max_uD_u$ and $\delta=\delta_K(h)\le D_h\,$. The support of $u$ is an independent set of lifted vertices, which extends to a basis of $\R^{d+3}$ by lifted vertices, so Cramer's rule and Proposition~\ref{prop-bounds} give $D_h\le B(d-1)$. Also $s-1\le\dim\ker A\le3(d+1)-(d+3)=2d$.

Let $\Lambda$ be the set of lattice points $(c,k)$ with $Ac=kh$; because $P_h$ has a positive point, this is the lattice of the linear span of $Q_h\,$, and its heights are exactly $q\Z$. Let $g$ be the greatest common divisor of the clearing multipliers at most $D_h$ (it is a multiple of $q$). Choose the lattice point of $Q_h$ of height $\delta$ and complete it to a basis of the span of $Q_h$ by $s-1$ of the extreme rays $(D_uu,D_u)$; all basis heights are multiples of $g$, with total at most $\delta+(s-1)D_h\,$. For any $m\in q\Z$, take a point of $\Lambda$ of height $m$ and subtract the integer parts of its coordinates in this basis. The result is a lattice point of $Q_h$ of height at most $\delta+(s-1)D_h$ and congruent to $m$ modulo $g$. The clearing multipliers at most $D_h\,$, divided by $g$, generate a numerical semigroup with least element $\delta/g$ and generators at most $D_h/g$; a shortest path in the residues modulo $\delta/g$ shows that its conductor is at most $(\delta/g-1)(D_h/g-1)$. Adding, every $m\in q\Z$ with $m\ge\delta+(s-1)D_h+\delta D_h/g$ is the height of a lattice point of $Q_h\,$, hence clears. The bounds $\delta\le D_h\le B(d-1)$, $s-1\le2d$ and $g\ge1$ give the claim.
\end{proof}

\begin{proposition}[Exponential delay at degree three]\label{prop-delay}
Let $a_0=\cdots=a_4=1$ and $a_i=a_{i-1}+a_{i-5}$ for $i\ge5$. For every $n\ge8$ with $n\equiv2\pmod6$, three staircase simplices in $\R^{n+6}$ have an indecomposable hole $h$ of degree three, trivial modulo the ray lattice of its minimal face, with $\mathcal C_K(h)=\langle a_n\,,a_n+1\rangle\setminus\{0\}$. Its largest nonclearing multiplier is $a_n^2-a_n-1$.
\end{proposition}
\begin{proof}
Use profile variables $x_0\,,\ldots,x_n\,$, with $x_i$ in class $i\bmod3$, and further variables $D,W$ in class $0$ and $Y,Z$ in class $1$. At scale $k$ impose the rows
\[
 x_{i-1}+x_i=k\ (1\le i\le4),\qquad x_i+x_{i-1}+x_{i-5}=k\ \text{or}\ 2k\ (5\le i\le n,\ i\text{ odd or even}),
\]
\[
 D+x_n=k\ (\text{twice}),\quad Y+x_{n-2}=k,\quad Z+x_{n-3}=k,\quad W+Z=k,\quad x_n+Y+W=2k .
\]
Each row uses at most one variable per class, and there are $n+6$ rows; they are the spatial coordinates. Put $t_i=x_i$ for even $i$ and $t_i=k-x_i$ for odd $i$. The first two families give $t_i=a_is$ with $s=x_0\,$, and then $D=k-a_ns$, $Y=k-a_{n-2}s$, $Z=a_{n-3}s$, $W=k-a_{n-3}s$. With $f_i=t_i-t_{i-2}-t_{i-3}$ one has $t_i-t_{i-1}-t_{i-5}=f_i-f_{i-1}+f_{i-2}\,$, so the gate rows give $f_i=f_{i-1}-f_{i-2}$ for $i\ge5$, while $f_3=f_4=-s$; hence $f_n=0$ for $n\equiv2\pmod6$, and the last row, whose residual is a combination of the residuals of the three rows before it and $f_n\,$, is implied by the others. The resulting linear dependence among the rows supplies weights: $1,-1,0,1$ on the four seed rows, $(-1)^{i+1}c_{n-i}$ on gate row $i$, where $(c_0\,,\ldots,c_5)=(1,1,0,-1,-1,0)$ is extended periodically, $-1$ and $1$ on the two copies of $D+x_n=k$, and $-1,1,-1,1$ on the last four rows. Every variable has total weight zero on its rows, and a nonzero weight on some row, so every nonsingleton support contains both signs; each class has a positive row outside its supports.

In each class order the variables by nonincreasing value at $s=1$, $k=a_n+\frac12$. List the outside rows of weight zero, then the outside positive rows, then the supports of the variables in this order, each with its negative rows first, its zero rows next and its positive rows last, and finally the outside negative rows. As in Lemma~\ref{lem-ordered}, the support endpoints are exactly the maximum-weight prefixes, and with $u_i$ the prefix before the first support of class $i$ and $b$ the vector of row coefficients of $k$, the vertex expressions of $kh$, $h=(\sum_iu_i+b,1,1,1)$, correspond to the solutions of the rows that are nonincreasing in each class with values in $[0,k]$. These constraints force $s>0$ and $a_n\le k/s\le a_n+1$; conversely on this interval all values lie in $[0,k]$ and respect the orders, because every value divided by $s$ is an integer or $k/s$ minus an integer and no integer lies strictly between $a_n$ and $a_n+1$. With $b'=D$ and $a'=s-D$ the solutions are $a',b'\ge0$ with $a_na'+(a_n+1)b'=k$, integral exactly when the profiles are. So $\mathcal C_K(h)=\langle a_n\,,a_n+1\rangle\setminus\{0\}$, $h$ is an indecomposable hole of degree three, hence indecomposable, and because $a_n$ and $a_n+1$ both clear, $h$ is trivial modulo the ray lattice of its minimal face. The largest gap of $\langle a,a+1\rangle$ is $a^2-a-1$.
\end{proof}

\begin{proof}[Proof of Theorem~\ref{thm-delay}]
The upper bound is Lemma~\ref{lem-conductor}. Because $x^5-x^4-1=(x^2-x+1)(x^3-x-1)$, the dominant root of the recurrence is $\rho$ and $a_n=\Theta(\rho^n)$. Proposition~\ref{prop-delay} gives $L(n+6)\ge a_n^2-a_n-1$, and appending at most five coordinates at the end of all three orders embeds the cone as a face with the same indecomposable hole and clearing set, which covers every large $d$.
\end{proof}

\subsection{Dimension at most five}\label{sec-dim5}
The example of Section~\ref{sec-intro}, which is the word disjunction of $ABCAB$ (Theorem~\ref{thm-words}), lives in dimension six. No smaller example exists.

\begin{proposition}[Dimension at most five]\label{prop-dim5}
Three staircase simplices in $\R^d$ with $d\le5$ have IDP. Thus $H_2(d)=0$ for $d\le5$, while $H_2(6)\ge1$: six is the smallest spatial dimension in which three staircase simplices fail IDP.
\end{proposition}

Throughout this subsection we use the multiplicity coordinates~\eqref{eq-multiplicity}, in which the lifted vertices are $(v,\eps_i)$. The \emph{growth order} of $\mathcal P_\sigma$ is $(\sigma(d),\ldots,\sigma(1))$; the vertices of $\mathcal P_\sigma$ are the indicators $\mathbf e_S:=\sum_{p\in S}\mathbf e_p$ of the prefixes $S$ of its growth order, and $\mathbf e:=\mathbf e_{[d]}\,$. A lattice point of $K$ \emph{splits} if it lies in $S_K\,$. We use two elementary facts.

\smallskip\noindent\emph{Level sets.} If $y\in\Z^d$ is nonincreasing along the growth order of $\mathcal P_i$ and $0\le y\le\lambda_i\,$, then $y$ is the sum of the $\lambda_i$ vertices $\mathbf e_{\{y\ge s\}}\,$, $1\le s\le\lambda_i\,$, of $\mathcal P_i\,$.

\begin{lemma}[Complementation]\label{lem-complement}
The map $(x,\lambda)\mapsto\bigl((\textstyle\sum_i\lambda_i)\mathbf e-x,\ \lambda\bigr)$ is a lattice automorphism of $\Z^{d+3}$ that maps the lifted vertices of $\mathcal P_{\sigma_0}\,,\mathcal P_{\sigma_1}\,,\mathcal P_{\sigma_2}$ bijectively onto those of the staircase simplices with reversed growth orders. It therefore maps the cone, its vertex monoid, and the split points of one triple onto those of the other. A chain of selected vertices of $\mathcal P_i$ with lowest vertex $\mathbf e_L\,$, blocks $S_1\,,\ldots,S_\ell$ in chronological order and heights $\beta_1\,,\ldots,\beta_\ell$ (defined in~\eqref{eq-dim5-heights} below) goes to a chain with the same blocks in reversed chronological order and heights $\lambda_i-\beta_j\,$.
\end{lemma}
\begin{proof}
The map is integral and is its own inverse. It sends $(v,\eps_i)$ to $(\mathbf e-v,\eps_i)$, and $\mathbf e-\mathbf e_S$ is the indicator of a prefix of the reversed growth order. Applied to $\lambda_i\mathbf e_L+\sum_j\beta_j\mathbf e_{S_j}\,$, it gives $\lambda_i\mathbf e_{L'}+\sum_j(\lambda_i-\beta_j)\mathbf e_{S_j}\,$, where $L'$ is the complement of $L\cup\bigcup_jS_j\,$.
\end{proof}

\begin{lemma}[Prefix partition]\label{lem-prefix-partition}
Let $S_0\,,S_1\,,S_2$ be linearly ordered finite sets such that every element of $T:=S_0\cup S_1\cup S_2$ lies in exactly two of them, and put $I_{jk}:=S_j\cap S_k\,$. If some $I_{jk}$ is empty or $|T|\le4$, then $T$ is the disjoint union of a prefix of $S_0\,$, a prefix of $S_1$ and a prefix of $S_2\,$.
\end{lemma}
\begin{proof}
The sets $I_{01}\,,I_{02}\,,I_{12}$ partition $T$. If $I_{12}=\varnothing$, then $S_0=T$; take all of $S_0$ and empty prefixes elsewhere, and similarly for the other cases. Otherwise all three are nonempty and $|T|\le4$, so two of them are singletons, and these two share one set. Relabel so that $I_{02}=\{p\}$ and $I_{12}=\{q\}$; then $S_2=\{p,q\}$. If $p$ precedes $q$ in $S_2\,$, take the prefix $\{p\}$ of $S_2$ and all of $S_1=I_{01}\cup\{q\}$; otherwise take $\{q\}$ and all of $S_0=I_{01}\cup\{p\}$.
\end{proof}

\begin{lemma}[Small three-class determinants]\label{lem-det5}
Let $M$ be a nonsingular $k\times k$ $0/1$ matrix with $k\le5$ whose columns are partitioned into three nonempty classes with pairwise disjoint supports inside each class. Then $|\det M|\le3$. If $|\det M|=3$, then $k=5$, the classes have sizes $2,2,1$, and after permuting rows, interchanging the two classes of size two, and permuting columns inside classes, $M$ is, with columns ordered $(X_1\,,X_2\,,Y_1\,,Y_2\,,Z_1)$,
\[
 M_*:=\begin{pmatrix}
 0&0&0&1&1\\
 0&1&0&0&1\\
 0&1&1&0&0\\
 1&0&0&1&0\\
 1&0&1&0&1
 \end{pmatrix}.
\]
\end{lemma}
\begin{proof}
The bound $|\det M|\le3$ is \cite[Lemma~5.5]{Stair}; we need the equality case, so we give a self-contained argument. Because $k\le5$, some class is a single column $w$; call the other classes $C_1$ and $C_2$ and negate the columns of $C_2\,$. Every row has at most one entry in $C_1$ and at most one in $C_2\,$. Let $\Gamma$ be the multigraph on the vertex set $\{0\}\cup C_1\cup C_2$ in which a row with entries in both classes is an edge between its two columns and a row with an entry in only one class is an edge from that column to $0$. Orient every edge from its $C_1$ end to its $C_2$ end, where $0$ counts as the $C_2$ end of an edge to $C_1$ and as the $C_1$ end of an edge to $C_2\,$. The columns $C_1\cup C_2$ of $M$ then form the transpose $N$ of the incidence matrix of the directed graph $\Gamma$ with the row of vertex $0$ deleted, so $N$ is totally unimodular.

If a row meets neither $C_1$ nor $C_2\,$, its only entry is in $w$, and expanding along that row gives $|\det M|=1$. Otherwise all $k$ rows are edges, so $\Gamma$ has $k$ vertices and $k$ edges. As $N$ has rank $k-1$, $\Gamma$ is connected, hence has a unique cycle $Z$, possibly of length two. The cofactors of the column $w$ form a vector in the left kernel of $N$. This kernel is spanned by the signed cycle vector $\varepsilon$ of $Z$, where $\varepsilon_r=\pm1$ according as $Z$ traverses the edge $r$ with or against its orientation and $\varepsilon_r=0$ off $Z$. Because the cofactors are $0$ or $\pm1$ and not all zero, they equal $\pm\varepsilon$, so
\[
 |\det M|=\Bigl|\sum_{r\in Z}\varepsilon_rw_r\Bigr|.
\]
Both edges of $Z$ at a vertex of $C_1$ leave it and both edges at a vertex of $C_2$ enter it, so the signs alternate there. At $0$ they alternate unless $Z$ passes from $C_1$ through $0$ to $C_2\,$. If the signs alternate all around $Z$, its length is even and at most $k\le5$, hence at most four, and the sum is at most $2$ in absolute value. Otherwise $Z$ has odd length $\ell\le5$, and the sum is at most $(\ell+1)/2\le3$ in absolute value. Equality forces $\ell=5$. Then $Z$ passes through all five vertices, alternates between $C_1$ and $C_2$ away from $0$, and so $|C_1|=|C_2|=2$. The three edges of the majority sign are the two edges at $0$ and the edge opposite $0$; equality forces $w=1$ on them and $w=0$ on the remaining two edges. Because $\Gamma=Z$, this determines $M$ up to the stated relabelings. With $C_1=\{X_1\,,X_2\}$, $C_2=\{Y_1\,,Y_2\}$, $w=Z_1$ and $Z=0X_2Y_1X_1Y_20$, it gives $M_*\,$.
\end{proof}

\begin{proof}[of Proposition~\ref{prop-dim5}]
\emph{Reduction.} Suppose some lattice point of $K$ does not split, and let $h=(x,\lambda)$ be one of least degree $\sum_i\lambda_i\,$. In every expression $h=\sum_r\gamma_rr$ with $\gamma_r\ge0$ over lifted vertices, all $\gamma_r<1$: otherwise $h-r$ is a lattice point of $K$ of smaller degree, so it splits, and then $h$ splits. By Carath\'eodory's theorem, choose such an expression with linearly independent rays and all $\gamma_r\in(0,1)$. The rays of any two constituents form a totally unimodular matrix, by the argument of Proposition~\ref{prop-two} in multiplicity coordinates: in a square submatrix, replace the selected columns of each constituent, except the lowest, by consecutive differences, and negate one constituent; then every row, including the multiplicity rows, has at most one $+1$ and at most one $-1$. An independent expression of a lattice point in such rays has integral coefficients, so all three constituents occur. Because $\lambda_i=\sum\gamma_r$ over the rays of $\mathcal P_i$ is a positive integer, each constituent contributes at least two rays.

\emph{Heights.} Let the selected vertices of $\mathcal P_i$ be $\mathbf e_{L_i}=v_{i,0}\,,\ldots,v_{i,\ell_i}\,$, with nested supports, coefficients $\gamma_{i,0}\,,\ldots,\gamma_{i,\ell_i}\,$, and blocks $S_{i,j}:=\supp v_{i,j}\setminus\supp v_{i,j-1}\,$, which are consecutive segments of the growth order. With the heights $\beta_{i,j}:=\sum_{r\ge j}\gamma_{i,r}\,$,
\begin{equation}\label{eq-dim5-heights}
 x=\sum_{i=0}^2\Bigl(\lambda_i\mathbf e_{L_i}+\sum_{j=1}^{\ell_i}\beta_{i,j}\mathbf e_{S_{i,j}}\Bigr),
\end{equation}
and $\gamma_{i,j}=\beta_{i,j}-\beta_{i,j+1}$ (with $\beta_{i,\ell_i+1}:=0$) and $\gamma_{i,0}=\lambda_i-\beta_{i,1}\,$. Conversely, integer heights $\lambda_i\ge m_{i,1}\ge\cdots\ge m_{i,\ell_i}\ge0$ give the sum $\lambda_i\mathbf e_{L_i}+\sum_jm_{i,j}\mathbf e_{S_{i,j}}$ of $\lambda_i$ selected vertices of $\mathcal P_i\,$, by the level-set fact.

\emph{Denominators.} The block indicators $\mathbf e_{S_{i,j}}$ form a $0/1$ matrix $B$ with $k:=\sum_i\ell_i$ columns in three nonempty classes with disjoint supports inside each class. Subtracting the lowest selected ray of each class from the others and eliminating the multiplicity rows turns the ray matrix into a block matrix with diagonal blocks $B$ and $I_3\,$, so $B$ has full column rank and $k\le d\le5$. By~\eqref{eq-dim5-heights}, $B\beta$ is integral. Choose a nonsingular $k\times k$ row submatrix $M$ of $B$. By Cramer's rule all $\beta_{i,j}\,$, and hence all $\gamma_r\,$, lie in $|\det M|^{-1}\Z$, and Lemma~\ref{lem-det5} gives $|\det M|\in\{1,2,3\}$. Determinant one would make every $\gamma_r$ integral. If $|\det M|=2$, every $\gamma_r$ equals $\frac12$, each constituent contributes an even number of rays, and there are $k+3\le8$ rays, so the numbers of rays are $(2,2,2)$ or, up to order, $(4,2,2)$. If $|\det M|=3$, then $k=d=5$, $M=B$, the numbers of rays are $(3,3,2)$ up to order, and $B$ is $M_*$ up to relabeling. We show that $h$ splits in each case, a contradiction.

\emph{Two rays from each constituent.} Here $\lambda=(1,1,1)$, each $\mathcal P_i$ has one block $S_i$ of height $\frac12$, and integrality of $x$ puts every coordinate in zero or two of the $S_i\,$. With $T:=\bigcup_iS_i\,$, $x=\sum_i\mathbf e_{L_i}+\mathbf e_T\,$. Order each $S_i$ by its growth order. If Lemma~\ref{lem-prefix-partition} applies, with prefixes $Q_i\,$, then the vertices $\mathbf e_{L_i\cup Q_i}$ sum to $x$. Otherwise $T=[5]$ and all $I_{jk}$ are nonempty; some $I_{jk}$ is a singleton, say $I_{12}=\{p\}$. Then $S_0=T\setminus\{p\}$, so $L_0=\varnothing$ or $L_0=\{p\}$. In the first case $x=\mathbf e+\mathbf e_{L_1}+\mathbf e_{L_2}$; in the second $x=\mathbf 0+\mathbf e_{L_1\cup S_1}+\mathbf e_{L_2\cup S_2}\,$, because $\mathbf e_{S_1}+\mathbf e_{S_2}=\mathbf e_T+\mathbf e_{\{p\}}\,$. Both are vertex sums.

\emph{Four rays from one constituent.} Call it $\mathcal P_X\,$, with blocks $D_1\,,D_2\,,D_3$ in chronological order, heights $\frac32,1,\frac12$ and $\lambda_X=2$; the other two, $\mathcal P_Y$ and $\mathcal P_Z\,$, have $\lambda=1$ and single blocks $S_Y\,,S_Z$ of height $\frac12$. Put $S_X:=D_1\cup D_3\,$, ordered by the growth order, and $g:=2\mathbf e_{L_X}+\mathbf e_{D_1}+\mathbf e_{D_2}\,$. Then
\[
 x=g+\mathbf e_{L_Y}+\mathbf e_{L_Z}+\tfrac12(\mathbf e_{S_X}+\mathbf e_{S_Y}+\mathbf e_{S_Z}),
\]
so every coordinate lies in zero or two of $S_X\,,S_Y\,,S_Z$; let $T$ be their union. If Lemma~\ref{lem-prefix-partition} applies, with prefixes $Q_X\,,Q_Y\,,Q_Z\,$, then $g+\mathbf e_{Q_X}$ takes the value $2$ on $L_X\,$, then $2$ and $1$ on $D_1\,$, $1$ on $D_2\,$, then $1$ and $0$ on $D_3\,$, and $0$ afterwards, where $Q_X$ meets $D_3$ only if $D_1\subseteq Q_X\,$. By the level-set fact it is the sum of two vertices of $\mathcal P_X\,$, and with $\mathbf e_{L_Y\cup Q_Y}$ and $\mathbf e_{L_Z\cup Q_Z}$ this splits $h$.

Otherwise $T=[5]$, $d=5$, and $B$ is the nonsingular $5\times5$ matrix with columns $D_1\,,D_2\,,D_3\,,S_Y\,,S_Z\,$, so its rows are distinct. Every coordinate lies in exactly two of $S_X\,,S_Y\,,S_Z\,$, so its row type is one of $D_1Y$, $D_1Z$, $D_3Y$, $D_3Z$, $D_2YZ$ and $YZ$, where for instance $D_1Y$ means a coordinate of $D_1$ in $S_Y$ but not in $S_Z\,$. The nonempty block $D_2$ is therefore a single coordinate of type $D_2YZ$. If type $YZ$ were absent, the other four rows would be $D_1Y,D_1Z,D_3Y,D_3Z$; after expanding along the column $D_2$ they form the incidence matrix of a four-cycle, which is singular. So type $YZ$ occurs once, and the remaining three rows contain both types for one of $D_1\,,D_3$ and one type for the other. By Lemma~\ref{lem-complement}, which reverses the chronological order of $D_1\,,D_2\,,D_3$ and keeps the heights $\frac32,1,\frac12$ in chronological order, and by interchanging $Y$ and $Z$, we may name the coordinates $1,\ldots,5$ so that
\[
 D_1=\{1,2\},\quad D_2=\{4\},\quad D_3=\{5\},\qquad S_Y=\{1,3,4\},\quad S_Z=\{2,3,4,5\}.
\]
Coordinate $3$ lies in no block of $\mathcal P_X\,$, so $L_X=\{3\}$ or $L_X=\varnothing$. In both cases
\[
 x=2\mathbf e_{L_X}+2\mathbf e_{\{1,2,4\}}+\mathbf e_{\{3,5\}}+\mathbf e_{L_Y}+\mathbf e_{L_Z}\,.
\]
If $L_X=\varnothing$, this is $\mathbf e+\mathbf e_{\{1,2,4\}}$ from $\mathcal P_X$ plus $\mathbf e_{L_Y}$ and $\mathbf e_{L_Z}\,$. If $L_X=\{3\}$, it is $\mathbf 0+\mathbf e_{\{1,2,3\}}$ from $\mathcal P_X$ plus $\mathbf e_{L_Y\cup S_Y}$ and $\mathbf e_{L_Z\cup S_Z}\,$. In both cases $h$ splits.

\emph{Determinant three.} Call the constituents with three rays $\mathcal P_X\,,\mathcal P_Y$ and the one with two rays $\mathcal P_Z\,$, with blocks $X_1\,,X_2\,$, $Y_1\,,Y_2\,$, $Z_1$; the indices name blocks and need not be chronological. Name the coordinates $1,\ldots,5$ by the rows of $M_*\,$, so that
\[
 X_1=\{4,5\},\quad X_2=\{2,3\},\quad Y_1=\{3,5\},\quad Y_2=\{1,4\},\quad Z_1=\{1,2,5\},
\]
and hence $L_X\subseteq\{1\}$, $L_Y\subseteq\{2\}$, $L_Z\subseteq\{3,4\}$. The rows of $B\beta\equiv0$ modulo $\Z$ give
\[
 (\beta_{X_1}\,,\beta_{X_2}\,,\beta_{Y_1}\,,\beta_{Y_2}\,,\beta_{Z_1})\equiv\theta(1,-1,1,-1,1),\qquad 3\theta\equiv0,
\]
and $\theta\not\equiv0$, because otherwise all $\gamma_r$ would be integral. Lemma~\ref{lem-complement} keeps the blocks and the matrix $B$ and replaces every height $\beta$ by $\lambda_i-\beta$, hence $\theta$ by $-\theta$; so we may assume $\theta=\frac13$. For a constituent with three rays, all three coefficients lie in $\{\frac13,\frac23\}$ and sum to an integer, so they are all $\frac13$ or all $\frac23$. Let $\alpha_X=1$ if $X_1$ precedes $X_2$ and $\alpha_X=0$ otherwise, and define $\alpha_Y$ likewise. The residues then force
\[
 \beta_{X_1}=\tfrac13+\alpha_X\,,\quad\beta_{X_2}=\tfrac23,\quad\beta_{Y_1}=\tfrac13+\alpha_Y\,,\quad\beta_{Y_2}=\tfrac23,\quad\beta_{Z_1}=\tfrac13,
\]
with $\lambda_X=1+\alpha_X\,$, $\lambda_Y=1+\alpha_Y\,$, $\lambda_Z=1$. Hence
\[
 x-\lambda_X\mathbf e_{L_X}-\lambda_Y\mathbf e_{L_Y}-\mathbf e_{L_Z}=(1,\,1,\,1+\alpha_Y\,,\,1+\alpha_X\,,\,1+\alpha_X+\alpha_Y).
\]
For $(u,v)\in\{(1,0),(0,1),(1,1)\}$ take the integer heights
\[
 m_{X_2}=u,\qquad m_{Y_2}=v,\qquad m_{X_1}=1+\alpha_X-v,\qquad m_{Y_1}=1+\alpha_Y-u.
\]
They are admissible: if $\alpha_X=0$, then $X_2$ comes first and $1\ge u\ge1-v\ge0$ because $u+v\ge1$; if $\alpha_X=1$, then $X_1$ comes first and $2\ge2-v\ge u\ge0$ because $u+v\le2$; the same holds for $Y$ with $u$ and $v$ interchanged. These heights give $\lambda_X$ vertices of $\mathcal P_X$ and $\lambda_Y$ vertices of $\mathcal P_Y\,$. Subtracting them leaves $\mathbf e_{L_Z}$ plus the vector with entries $(1-v,\,1-u,\,u+v-1)$ on the coordinates $(1,2,5)$ and $0$ on $3,4$. The three choices of $(u,v)$ make this the indicator of $\{1\}$, $\{2\}$ and $\{5\}$, respectively. Choosing the first coordinate of $Z_1$ in the growth order of $\mathcal P_Z$ leaves a vertex of $\mathcal P_Z\,$, so $h$ splits.

\emph{Sharpness.} The example of Section~\ref{sec-intro} gives $H_2(6)\ge1$, as does Theorem~\ref{thm-threshold}.
\end{proof}

\subsection{Arithmetic off the designated faces}\label{sec-offface}
The faces constructed above control the arithmetic of one face; they do not describe the whole ambient Hilbert basis. In particular, an indecomposable hole need not lie in any simplex face. Recall that a lattice point of degree at most two involves at most two constituents, so it is a vertex sum by Proposition~\ref{prop-two}.

\par\medskip\noindent\emph{An indecomposable hole outside every simplex face.}
Take $q=4$ in Theorem~\ref{thm-onehole}, the spatial coordinates in the order $(P_1\,,P_2\,,P_3\,,Q_1\,,Q_2\,,Q_3\,,R_1\,,R_2\,,R_3\,,g)$, multiplicities $(\lambda_A\,,\lambda_B\,,\lambda_C)$, and the growth orders
\begin{gather*}
 A=(Q_1\,,R_3\,,P_3\,,R_2\,,P_2\,,g,P_1\,,R_1\,,Q_2\,,Q_3),\\ B=(g,P_1\,,Q_1\,,Q_2\,,P_2\,,R_1\,,R_2\,,Q_3\,,P_3\,,R_3),\\
 C=(P_3\,,Q_3\,,Q_2\,,g,R_2\,,R_3\,,Q_1\,,R_1\,,P_1\,,P_2),
\end{gather*}
which are the growth orders of Theorem~\ref{thm-onehole} for $q=4$, with ties broken as displayed; the form $f$ below exposes a different face of the same hull. Write $A_j\,,B_j\,,C_j$ for the lifted prefixes of length $j$.
\begin{proposition}[An indecomposable hole in no simplex face]\label{prop-offface-y}
The point
\[y=(2,2,1,2,2,1,1,2,1,3;\,1,2,1)\]
is a Hilbert basis element of $K$ of degree four with $\mathcal C_K(y)=2\Z_{>0}\,$. Its minimal face $H$ has the ten rays $A_0\,,A_6\,,A_7\,,B_1\,,B_2\,,B_5\,,B_7\,,B_8\,,C_0\,,C_6\,$, rank nine and the single relation $A_6+B_2=A_7+B_1\,$. In particular $y/4$ lies in no simplex face of the hull, and $\Hilb(K)$ contains indecomposable holes other than $h_4\,$.
\end{proposition}
\begin{proof}
First, $2y=A_0+A_7+B_1+B_5+B_7+B_8+C_0+C_6\,$. The form $f=\lambda_B-2x_{P_2}+3x_{P_3}+2x_{Q_1}+x_{R_1}-x_{R_2}-x_{R_3}-x_g$ takes the following values on the prefixes of lengths $0,\ldots,10$:
\begin{gather*}
 A:\ 0,2,1,4,3,1,0,0,1,1,1;\qquad B:\ 1,0,0,2,2,0,1,0,0,3,2;\\ C:\ 0,3,3,3,2,1,0,2,3,3,1.
\end{gather*}
So $f\ge0$ on $K$, its zero face $H$ has the ten listed rays, and $f(y)=0$. In a linear relation among these rays, let $c$ be the coefficient of $C_6$; the rows $R_1\,,Q_3\,,Q_2$ give $b_7+b_8=0$, $b_8=-c$, $b_5=-c$, and the rows $P_2\,,P_3$ give $a_6+a_7=c=-c$. Over $\R$ this forces $c=0$, and the remaining rows leave only $a_6=b_2=-a_7=-b_1\,$. Hence every nonnegative real expression of $y$ has coefficients
\begin{gather*}
 a_0=b_5=b_7=b_8=c_0=c_6=\tfrac12,\\ a_6=b_2=\tau,\qquad a_7=b_1=\tfrac12-\tau,\qquad 0\le\tau\le\tfrac12 .
\end{gather*}
No lifted vertex can be subtracted from $y$ within $K$. A nontrivial decomposition of a point of degree four has a summand of degree at most two, which is a vertex sum; so $y$ is indecomposable. In an integer expression of $py$ the six fixed coefficients equal $p/2$, so $\mathcal C_K(y)=2\Z_{>0}\,$.

At $\tau=\frac14$ all ten coefficients are positive, so $y$ lies in the relative interior of $H$, which is therefore its minimal face. A face of $K$ containing $y$ contains $H$, which has ten rays in rank nine, so it is not simplicial; hence $y/4$ lies in no simplex face of $\mathcal D$.
\end{proof}

\section{Complete ambient Hilbert bases for word disjunctions}\label{sec-words}
This section proves result~\hyperlink{res-D}{(D)}: for word disjunctions we determine the entire ambient Hilbert basis and the clearing sets of all lattice points. Their initial facets are the common-order pair-prefix chains of Theorem~\ref{thm-interval}, and the full staircase simplices complete them without adding indecomposable holes. Subsection~\ref{sec-words-overview} describes word disjunctions informally, states the main theorem and outlines its proof; Subsections~\ref{sec-words-main}--\ref{sec-words-proof} prove it through a unimodular graph representation and the Hilbert basis of chordal graph cones, and Subsections~\ref{sec-words-cons}--\ref{sec-grouped} derive consequences.

\subsection{Overview}\label{sec-words-overview}
Take any string (a \emph{word}) $W=w_1\cdots w_\ell$ over the symbols $A,B,C$. Its positions $1,\ldots,\ell$ are the coordinates, and a separator coordinate $g:=\ell+1$ (the \emph{guard}) is added, so $d=\ell+1$. Each of the three simplices belongs to a pair of letters. It switches on the positions carrying those two letters from left to right, then the guard, then the positions of the third letter from right to left. With $W|_X$ listing the positions of letters in $X$, the \emph{word disjunction} $\mathcal D_W$ has growth orders
\[
 (W|_{AB}\,,g,\operatorname{rev}W|_C),\quad(W|_{AC}\,,g,\operatorname{rev}W|_B),\quad(W|_{BC}\,,g,\operatorname{rev}W|_A).
\]
Before the guard, every coordinate is shared by exactly two simplices, and all three use the same left-to-right order. Call a substring $w_{i+1}\cdots w_j$ ($0\le i<j\le \ell$) \emph{bad} if it contains all three letters and the spans from the first to the last occurrence of $A$, of $B$ and of $C$ have a common point, as in $ABCAB$; we then also call the pair $(i,j)$ bad. Let $h_{ij}$ have $\lambda=(1,1,1)$, $x_g=0$, and $x_p=2,1,0$ for $p\le i$, $i<p\le j$, $p>j$ respectively. Write $K_W\,$, $V_W$ and $S_W$ for the cone, lifted vertices and vertex monoid of $\mathcal D_W\,$. These notions are made precise in Subsection~\ref{sec-words-main}.

\begin{theorem}[Complete ambient classification]\label{thm-words-summary}
For every word $W$, $\Hilb(K_W)=V_W\cup\{h_{ij}:(i,j)\text{ bad}\}$. Every $h_{ij}$ has degree three and clearing set $2\Z_{>0}\,$, and $2y\in S_W$ for every lattice point $y$ of $K_W\,$. Moreover:
\begin{enumerate}
\item $\mathcal D_W$ has IDP if and only if $W$ contains neither $XYZXY$ nor $XYZYX$ as a subsequence, for distinct letters $X,Y,Z$; equivalently, for some two letters, every occurrence of one precedes every occurrence of the other (Proposition~\ref{prop-grouped});
\item the minimal faces of $\mathcal D_W$ failing IDP are exactly five-dimensional simplices $\Pi_{ij}$ of normalized volume two, one for each bad pair;
\item for $d\ge6$ the maximum number of indecomposable holes among word disjunctions is $\binom{d-4}2$, attained by the prefixes of $ABCABC\cdots$.
\end{enumerate}
\end{theorem}
\par\medskip\noindent\emph{Idea of the proof.} By Theorem~\ref{thm-words-graph}, $K_W$ is isomorphic to $C_{G_W}\,$, where $C_G:=\cone(\{(\mathbf e_v\,,0)\}\cup\{(-\mathbf e_u-\mathbf e_v\,,1):uv\in E(G)\})$ for a graph $G$. The graph $G_W$ is a $2$-tree, hence chordal, and the following argument works for chordal graphs (Theorem~\ref{thm-words-chordal}). A lattice point $(y,k)$ of $C_G$ has a minimum fractional edge cover of the deficits $\max(0,-y_v)$ whose nonintegral edges form vertex-disjoint odd cycles with weights in $\frac12+\Z$ (Lemma~\ref{lem-words-cover}); set the integer parts aside. Chords shorten every such cycle to a triangle. Two half-weight triangles joined by an edge become three integral edges, and an unpaired triangle uses the half unit of slack forced by the integrality of $k$. Only pairs of triangles with no joining edge remain, and they contribute exactly $(-\mathbf e_{T\cup T'}\,,3)$. An integral vertex expression of $q(-\mathbf e_{T\cup T'}\,,3)$ has $3q$ edges and, as its $y$-coordinates sum to $-6q$, no unit vectors; so every vertex of $T\cup T'$ has degree $q$, all edges lie in $T$ or $T'$, and each triangle has even degree sum $3q$, so $q$ is even. After reading $w_1\cdots w_i\,$, the latest vertices of $A$, $B$, $C$ form a triangle $T_i\,$; for $i<j$, $T_i$ and $T_j$ are disjoint and unjoined exactly when $(i,j)$ is bad (Lemma~\ref{lem-words-exc}).

\begin{proof}
The Hilbert basis, the degrees, the clearing sets and doubling are Theorem~\ref{thm-words} and Corollary~\ref{cor-words-two}; part~(1) is Corollary~\ref{cor-words-five} with Proposition~\ref{prop-grouped}, part~(2) is Theorem~\ref{thm-words-minface}, and part~(3) is Theorem~\ref{thm-words-count}.
\end{proof}

The facet $\{x_g=0\}$ of $\mathcal D_W$ is an initial face (Section~\ref{sec-ambient}) on which every coordinate can be switched on by at most two of the simplices. On every such face, twice every lattice point splits and no odd multiple of a nonsplitting lattice point does, so every indecomposable hole has clearing order $2$, and IDP is a linear-time $2$-SAT test (Theorems~\ref{thm-dichotomy} and~\ref{thm-twosat}).

\subsection{The word family}\label{sec-words-main}
Let $W=w_1\cdots w_\ell$ be a word on the alphabet $\{A,B,C\}$, $\ell\ge0$. Its positions $1,\ldots,\ell$ are spatial coordinates; add a \emph{guard} coordinate $g:=\ell+1$, so $d=\ell+1$. For $X\subseteq\{A,B,C\}$ let $W|_X$ be the list of positions whose letters lie in $X$, in word order, and let $\operatorname{rev}$ reverse a list. The \emph{word disjunction} $\mathcal D_W$ is the disjunctive hull of the three staircase simplices with growth orders
\begin{gather*}
 \gamma_0=\bigl(W|_{AB}\,,\,g,\,\operatorname{rev}W|_C\bigr),\qquad
 \gamma_1=\bigl(W|_{AC}\,,\,g,\,\operatorname{rev}W|_B\bigr),\\
 \gamma_2=\bigl(W|_{BC}\,,\,g,\,\operatorname{rev}W|_A\bigr).
\end{gather*}
Thus, before the guard, the letters $A,B,C$ have availability $01$, $02$, $12$, and the face $\{x_g=0\}$ is the initial face of Subsection~\ref{sec-dichotomy} for the three common-order lists $W|_{AB}\,,W|_{AC}\,,W|_{BC}\,$. Write $K_W\,$, $V_W\,$, $S_W:=\Z_{\ge0}V_W$ and $\overline S_W:=K_W\cap\Z^{d+3}$ for its cone, lifted vertices, vertex monoid and saturation, and $\mathcal C_W$ for clearing sets.

For $0\le i<j\le \ell$ let $W_{(i,j]}:=w_{i+1}\cdots w_j\,$, and let $I_X(i,j)$ be the interval between the first and the last occurrence of the letter $X$ in this subword, in the original position indices, empty if $X$ does not occur. Call $(i,j)$ \emph{bad} if $I_A(i,j)\cap I_B(i,j)\cap I_C(i,j)\ne\varnothing$, and let $\mathcal B(W)$ be the set of bad pairs. Define $h_{ij}$ in multiplicity coordinates by
\begin{equation}\label{eq-hij}
 \lambda(h_{ij})=(1,1,1),\qquad (h_{ij})_g=0,\qquad (h_{ij})_p=\begin{cases}2,&p\le i,\\1,&i<p\le j,\\0,&p>j,\end{cases}\qquad(1\le p\le \ell);
\end{equation}
in natural coordinates it has $z_1=z_2=1$ and $t=3$.

We prove the complete classification after establishing the graph representation and its chordal Hilbert basis. The completion is specific to these common-order pair-prefix chains; it does not assert full-hull IDP transfer for arbitrary $0/1$ sources.

\subsection{A unimodular graph representation}\label{sec-words-graph}
Let $a,b,c$ be the numbers of occurrences of $A,B,C$ in $W$. Build a graph $G_W$ as follows. Start with a triangle on vertices $A_0\,,B_0\,,C_0\,$. Reading $W$ from left to right, at the $j$th occurrence of a letter $X$ introduce a vertex $X_j\,$, join it to the two currently latest vertices of the other two letters, and make $X_j$ the latest vertex of its letter. Let $T_i$ be the triangle of latest vertices after the first $i$ positions, so that $T_0=\{A_0\,,B_0\,,C_0\}$. Every new vertex is attached to an existing edge, so $G_W$ is a $2$-tree, hence chordal, and its triangles are exactly the \emph{state triangles} $T_0\,,\ldots,T_\ell\,$. For a finite simple graph $G$ put
\[
 C_G:=\cone\bigl(\{(\mathbf e_v\,,0):v\in V(G)\}\cup\{(-\mathbf e_u-\mathbf e_v\,,1):uv\in E(G)\}\bigr)\subseteq\R^{V(G)}\times\R,
\]
in the lattice $\Z^{V(G)}\times\Z$, graded by $\deg(y,k):=\sum_vy_v+3k$, so that all displayed generators have degree one; let $S_G$ be their monoid.

\begin{theorem}[Unimodular graph representation]\label{thm-words-graph}
The integral linear map $\Psi:(x,\lambda)\mapsto(y,k)$ given by $k=t-x_g$ and
\[
 y_{A_0}=x_{A_1}-\lambda_0-\lambda_1\,,\qquad y_{A_j}=x_{A_{j+1}}-x_{A_j}\ (1\le j<a),\qquad y_{A_a}=x_g-x_{A_a}\,,
\]
where $x_{A_j}$ denotes the coordinate of the position of the $j$th occurrence of $A$ (and $y_{A_0}=x_g-\lambda_0-\lambda_1$ if $a=0$), with the same formulas for $B$ subtracting $\lambda_0+\lambda_2$ and for $C$ subtracting $\lambda_1+\lambda_2\,$, is a lattice isomorphism $\Z^{d+3}\to\Z^{V(G_W)}\times\Z$ that maps $V_W$ bijectively onto the generators of $C_{G_W}\,$. Hence $\Psi(K_W)=C_{G_W}\,$, $\Psi(S_W)=S_{G_W}\,$, and $\Psi$ preserves degree.
\end{theorem}
\begin{proof}
The inverse is given by $\lambda_2=k+\sum_jy_{A_j}\,$, $\lambda_1=k+\sum_jy_{B_j}\,$, $\lambda_0=k+\sum_jy_{C_j}\,$, $x_g=2k+\sum_vy_v\,$, and $x_{A_j}=x_g-\sum_{r\ge j}y_{A_r}$ and similarly for $B,C$, as one checks by telescoping; so $\Psi$ is unimodular.

Consider a lifted vertex of constituent $0$, so $\lambda=(1,0,0)$. If its prefix of $\gamma_0$ ends before the guard, it consists of the positions of $W|_{AB}$ up to some point of the word, $x_g=0$ and $k=1$. On the $A$ coordinates the prefix is the indicator of $A_1\,,\ldots,A_m$ for some $m$, so the telescoping formulas give $y_{A_m}=-1$ and $y_{A_r}=0$ for $r\ne m$ (with $y_{A_0}=-1$ when $m=0$); likewise $y=-\mathbf e_{B_{m'}}$ on the $B$ coordinates, where $B_{m'}$ is the latest $B$ vertex at that point, and $y=0$ on the $C$ coordinates because $\lambda_1+\lambda_2=0$ and no $C$ position is in the prefix. The image is $(-\mathbf e_{A_m}-\mathbf e_{B_{m'}}\,,1)$, the edge between the latest $A$ and $B$ vertices at that moment, which lies in a state triangle. The $a+b+1$ prefixes of $W|_{AB}$ give the $a+b+1$ edges of $G_W$ joining an $A$ vertex to a $B$ vertex, in the order of their creation. If the prefix contains the guard, then $x_g=t=1$, $k=0$, every $A$ and $B$ coordinate is one and the corresponding $y$ vanish, and the prefix contains the last $m$ occurrences of $C$; the telescoping formulas give $y=\mathbf e_{C_{c-m}}\,$. These $c+1$ prefixes give the unit generators of the $C$ vertices. Constituents $1$ and $2$ give, in the same way, the $AC$ edges and the units of the $B$ vertices, and the $BC$ edges and the units of the $A$ vertices. Every edge of $G_W$ joins two different letters, so this is a bijection; degrees are preserved because both gradings are one on the generators.
\end{proof}

Call a graph a \emph{triangle strip} if it arises from a triangle $T_0$ by a sequence of steps $T_{i+1}=(T_i\setminus\{v\})\cup\{w\}$, each adding a new vertex $w$ adjacent exactly to the two vertices of $T_i\setminus\{v\}$ for some $v\in T_i$; the state triangles admit this ordering with consecutive triangles sharing an edge, without requiring their adjacency graph to be a path.
\begin{proposition}[The realized graphs]\label{prop-words-strips}
The graphs $G_W\,$, for words $W$ in three letters, are exactly the triangle strips. Consequently, for every triangle strip $G$ with $\ell+3$ vertices, the cone $C_G$ with its lattice and its generator monoid is isomorphic, by a degree-preserving unimodular map, to the full cone of three staircase simplices in $\R^{\ell+1}$ with its lifted-vertex monoid.
\end{proposition}
\begin{proof}
By construction $G_W$ is a triangle strip, the step at position $i$ replacing the latest vertex of the letter $w_i\,$. Conversely, color the vertices of $T_0$ by $A,B,C$ and give each new vertex the color of the vertex it replaces; then every $T_i$ has one vertex of each color, the new vertex is joined to the latest vertices of the other two colors, and the strip is $G_W$ for the word of colors of the replaced vertices. Theorem~\ref{thm-words-graph} gives the isomorphism.
\end{proof}
\subsection{The Hilbert basis of a chordal graph cone}\label{sec-words-chordal}
A pair of triangles $T,T'$ of a graph $G$ is \emph{exceptional} if they are vertex disjoint and no edge of $G$ joins them; put $h_{T,T'}:=(-\mathbf e_{T\cup T'}\,,3)$, of degree three.

\begin{lemma}[Half-integral minimum covers]\label{lem-words-cover}
Let $(y,k)\in C_G\cap(\Z^{V(G)}\times\Z)$. There is a nonnegative edge weighting $a$ with $\deg_a(v)\ge\max(0,-y_v)$ for every $v$ and $\sum_ea_e\le k$, whose non-integral edges form vertex-disjoint odd cycles with weight in $\frac12+\Z$ on every edge of each of them; moreover $y+\deg_a$ is integral and nonnegative.
\end{lemma}
\begin{proof}
Membership $(y,k)\in C_G$ means $y=u-\sum_ea_e(\mathbf e_{u_e}+\mathbf e_{v_e})$ with $u\ge0$ and $\sum_ea_e=k$, so the polyhedron of weightings with the stated degree lower bounds contains one of total weight $k$. Minimize $\sum_ea_e$ over it and take an optimal vertex. In each component of its positive support, the incidence columns restricted to the tight degree rows are independent, so the component is a tree or contains exactly one cycle, which is odd (an even cycle gives an alternating dependence). In a tree at most one vertex is not tight, and peeling the other leaves shows that all weights are integral. In a unicyclic component every degree row is tight; peeling leaves the odd-cycle equations with integer right-hand sides, whose solution is half-integral, and either integral or in $\frac12+\Z$ on every cycle edge. Every vertex of a fractional cycle meets exactly two of its edges, so all weighted degrees are integral.
\end{proof}

The following self-contained chordal specialization of \cite[Section~4]{BL} records irredundancy and exact clearing sets.
\begin{theorem}[Hilbert basis of a chordal graph cone]\label{thm-words-chordal}
If $G$ is chordal, then
\begin{multline*}
 \Hilb(C_G)=\{(\mathbf e_v\,,0):v\in V(G)\}\cup\{(-\mathbf e_u-\mathbf e_v\,,1):uv\in E(G)\}\\ \cup\{h_{T,T'}:\{T,T'\}\text{ exceptional}\},
\end{multline*}
and every $h_{T,T'}$ has exact clearing set $2\Z_{>0}$ with respect to $S_G\,$.
\end{theorem}
\begin{proof}
If $G$ has no edges, $C_G$ is the orthant $\R^{V(G)}_{\ge0}\times\{0\}$; assume otherwise, so that $(0,1)=(-\mathbf e_u-\mathbf e_v\,,1)+(\mathbf e_u\,,0)+(\mathbf e_v\,,0)\in S_G$ for any edge $uv$. Take $a$ as in Lemma~\ref{lem-words-cover} and remove its integer parts. A half-weight odd cycle of length at least four has a chord, because $G$ is chordal; the chord splits it into a shorter odd cycle and an even cycle through the chord, and the half-weight original cycle has the same degree vector and total weight as the half-weight shorter odd cycle plus the alternating perfect matching of the interior vertices of the remaining path. Iterating replaces every fractional cycle by a half-weight triangle, the triangles remaining vertex disjoint. Pair the triangles arbitrarily. If $T,T'$ are joined by an edge $uv$, replace their half weights by $uv$ and the edges opposite $u$ in $T$ and $v$ in $T'$: three integer edges with the same degree vector and the same total. If no edge joins them, their contribution is exactly $h_{T,T'}$ plus units. If one triangle remains, cover it by two of its edges, which overcovers one vertex and costs $\frac12$ more; because the number of triangles is then odd, $\sum_ea_e\in\frac12+\Z$, and the integrality of $k$ supplies the extra half. In all cases the integer edge count is at most $k$, the surplus is supplied by copies of $(0,1)$, and positive unit generators make up the nonnegative remaining vertex coordinates. This proves that the displayed set generates the saturated monoid.

For an exceptional pair, $h_{T,T'}$ is half the sum of the six edge generators of $T\cup T'$. An expression of $qh_{T,T'}$ in $S_G$ uses exactly $3q$ edge generators, with degree at least $q$ at each of the six vertices; the $6q$ edge ends force degree exactly $q$ there and no other endpoints, so every edge lies in $T$ or in $T'$, as no edge joins them, and each triangle then has total degree $3q$, which is even. Hence $q$ is even, and even $q$ clear by scaling. All generators of $S_G$ have degree one, and by the generation argument every lattice point of degree at most two is in $S_G$; so each $h_{T,T'}\,$, being of degree three and outside $S_G\,$, is in the Hilbert basis. Distinct exceptional pairs give distinct vectors.
\end{proof}

\begin{lemma}[Exceptional state triangles are bad subwords]\label{lem-words-exc}
For $0\le i<j\le \ell$, the state triangles $T_i\,,T_j$ form an exceptional pair if and only if $(i,j)\in\mathcal B(W)$.
\end{lemma}
\begin{proof}
If a letter does not occur in $W_{(i,j]}\,$, its latest vertex is common to $T_i$ and $T_j\,$. Otherwise the triangles are disjoint. For letters $X\ne Y$ with first and last occurrences $f_X\,,l_Y$ in the subword, the old $X$ vertex of $T_i$ and the terminal $Y$ vertex of $T_j$ are adjacent exactly when $l_Y<f_X$: the terminal $Y$ vertex is joined at its birth to the then latest $X$ vertex, and acquires later neighbors only among newer vertices. Vertices of the same letter are never adjacent. Hence a joining edge exists exactly when two of the three occurrence intervals are disjoint, and three intervals on a line meet pairwise exactly when all three meet.
\end{proof}

\subsection{The ambient classification}\label{sec-words-proof}
The graph representation of Theorem~\ref{thm-words-graph} and the Hilbert basis of chordal graph cones in Theorem~\ref{thm-words-chordal} now give the classification.
\begin{theorem}[Ambient Hilbert basis and the initial facet]\label{thm-words}
For every word $W$,
\begin{equation}\label{eq-words-HB}
 \Hilb(K_W)=V_W\cup\{h_{ij}:(i,j)\in\mathcal B(W)\},
\end{equation}
and every $h_{ij}$ has degree three and exact clearing set $\mathcal C_W(h_{ij})=2\Z_{>0}\,$. The initial face $K_F=K_W\cap\{x_g=0\}$ is a facet, and
\begin{equation}\label{eq-words-completion}
 \overline S_W=\bigl(K_F\cap\Z^{d+3}\bigr)+S_W \,.
\end{equation}
Consequently $\mathcal D_W$ has IDP if and only if its initial facet has IDP, and this holds exactly when
\begin{equation}\label{eq-words-idp}
 I_A(0,\ell)\cap I_B(0,\ell)\cap I_C(0,\ell)=\varnothing .
\end{equation}
\end{theorem}
\begin{proof}
By Theorems~\ref{thm-words-graph} and~\ref{thm-words-chordal} and Lemma~\ref{lem-words-exc}, $\Hilb(K_W)$ consists of $V_W$ and the preimages of $h_{T_i,T_j}$ for $(i,j)\in\mathcal B(W)$, with degree three and clearing set $2\Z_{>0}\,$. The inverse formulas give $\lambda=(1,1,1)$, because each letter has two of the six vertices, and $x_g=2\cdot3-6=0$. For a position $p$ carrying the $m$th occurrence of a letter $X$, $x_p$ counts the selected $X$ vertices of index at least $m$; the selected ones are the latest $X$ vertices at times $i$ and $j$, and their indices are at least $m$ exactly when $p\le i$, respectively $p\le j$. This is~\eqref{eq-hij}, and proves~\eqref{eq-words-HB}.

The form $x_g$ is nonnegative on $V_W\,$, and the face $x_g=0$ contains the three empty prefixes and, because every position is available in two constituents before the guard, consecutive prefix differences giving every unit vector of a position; its rank is $\ell+3=\dim K_W-1$, so it is a facet. Every indecomposable hole lies on it, which gives~\eqref{eq-words-completion}; by Lemma~\ref{lem-face} the Hilbert basis of the facet consists of its lifted vertices and the same $h_{ij}\,$. If the three global intervals do not meet, no subword is bad, because occurrence intervals of a subword lie in the global ones, and $\mathcal D_W$ has IDP. If they meet, $(0,\ell)$ is bad and $h_{0\ell}$ is an indecomposable hole on the initial facet.
\end{proof}

\subsection{Consequences}\label{sec-words-cons}
We record consequences of Theorem~\ref{thm-words}: clearing at two, a test using five positions, the sharp maximum number of indecomposable holes, a unimodular facet that coexists with them, and the minimal non-IDP faces.
\begin{corollary}[Twofold clearing]\label{cor-words-two}
Every $y\in\overline S_W$ satisfies $2y\in S_W$; in particular every $y\in\overline S_W\setminus S_W$ has least clearing integer two.
\end{corollary}
\begin{proof}
Double a Hilbert basis expression of $y$ and use $2h_{ij}\in S_W\,$.
\end{proof}

\begin{corollary}[Five-position criterion]\label{cor-words-five}
$\mathcal D_W$ has IDP if and only if $W$ contains neither $XYZXY$ nor $XYZYX$ as a (not necessarily consecutive) subsequence, for distinct letters $X,Y,Z$.
\end{corollary}
\begin{proof}
Either subsequence puts an occurrence of $Z$ inside the occurrence intervals of $X$ and $Y$. Conversely, if the three intervals meet, let $Z$ be the letter whose first occurrence is latest; the other two letters occur both before and after it, and their first and last occurrences give one of the two patterns.
\end{proof}

\begin{theorem}[The sharp quadratic family]\label{thm-words-count}
For $d\ge6$, the maximum number of indecomposable holes of a word disjunction in spatial dimension $d$ is $\binom{d-4}2$. It is attained by the prefix of length $d-1$ of $ABCABC\cdots$, whose Hilbert basis has $3(d+1)+\binom{d-4}2$ elements; all its indecomposable holes have degree three, multiplicities $(1,1,1)$ and clearing set $2\Z_{>0}\,$. For $d\ge7$ these periodic words are, up to relabeling the letters, the only maximizers; for $d=6$ the words $XYZXY$ and $XYZYX$ both attain the maximum.
\end{theorem}
\begin{proof}
By Corollary~\ref{cor-words-five} every bad pair has $j-i\ge5$, and the profile~\eqref{eq-hij} determines $(i,j)$; there are $\sum_{r=5}^{\ell}(\ell-r+1)=\binom{\ell-3}2$ pairs with $j-i\ge5$, and $\ell=d-1$. In the periodic word every window of length five is $XYZXY$ up to relabeling, so every pair with $j-i\ge5$ is bad. Conversely, a maximizer has every length-five window bad, hence of the form $XYZXY$ or $XYZYX$ with distinct $X,Y,Z$. A window $XYZYX$ is not followed by a bad window, because the next window begins with $YZY$, and is not preceded by one, because the previous window reads $\ast XYZY$, whose fourth and fifth letters both differ from its second, whereas both patterns repeat the second letter in the fourth or fifth position; and $XYZXY$ can only be followed by $Z$. For $\ell\ge6$ this forces periodicity.
\end{proof}
\begin{corollary}[A unimodular facet under $\binom{d-4}2$ indecomposable holes]\label{cor-words-unimodular}
For every word $W$, the face $x_g=t$ of $K_W$ is a unimodular simplex facet, exposed by the primitive form $t-x_g\,$, which takes only the values $0,1$ on $V_W\,$, and every $h_{ij}$ has exposing slack three. In particular, for every $d\ge6$ the periodic word gives three staircase simplices in $\R^d$ with a unimodular simplex facet and $\binom{d-4}2$ ambient Hilbert basis elements of degree three, slack three and clearing set $2\Z_{>0}\,$, the maximum within the family.
\end{corollary}
\begin{proof}
Under $\Psi$ the form $t-x_g$ is the coordinate $k$, which is $1$ on the edge generators and $0$ on the unit generators; its coefficient of $t$ is one, so it is primitive. The face $k=0$ of $C_{G_W}$ is spanned by the $\ell+3$ unit vectors $(\mathbf e_v\,,0)$, which are part of a lattice basis; so it is a unimodular simplex cone of codimension one. Finally $k(h_{T_i,T_j})=3$, and Theorem~\ref{thm-words-count} gives the count.
\end{proof}

\begin{theorem}[Minimal non-IDP faces]\label{thm-words-minface}
For every bad pair $(i,j)$, the six edge generators of $T_i\cup T_j$ correspond to the vertices of a five-dimensional simplex face $\Pi_{ij}$ of $\mathcal D_W\,$, of normalized volume two in its induced lattice, whose proper faces are unimodular and whose only indecomposable hole is $h_{ij}\,$. The faces of $\mathcal D_W$ that fail IDP and are minimal with this property are exactly the $\Pi_{ij}\,$. Consequently every face of $\mathcal D_W$ of dimension at most four has IDP, and $\mathcal D_W$ has IDP if and only if every five-dimensional simplex face of $\mathcal D_W$ is unimodular.
\end{theorem}
\begin{proof}
In graph coordinates, the linear form $(y,k)\mapsto\sum_vz_vy_v+4k$ with $z_v=2$ on the six vertices of $T_i\cup T_j$ and $z_v=1$ on every other vertex takes positive values on all unit generators and the value $4-z_u-z_v$ on the edge generator of $uv$; this vanishes exactly on the six triangle edges, because no edge joins the triangles. The six edge generators are independent, and the torus kernel of their coefficient matrix is computed triangle by triangle: the three vertex rows force the three coefficients of a triangle to a common value in $\{0,\frac12\}$, and the $k$ row forces the two values to agree. So the induced-lattice index is two, every proper subset has index one, and the only nonzero box point is the half-sum $h_{ij}\,$. A face failing IDP has an indecomposable hole, which is an ambient one by Lemma~\ref{lem-face}, hence some $h_{ij}$; as a positive half-sum of the six rays of $\Pi_{ij}\,$, it forces the face to contain $\Pi_{ij}\,$. Finally, a lattice simplex face whose only lattice points are its vertices has IDP exactly when it is unimodular, and every simplex face of $\mathcal D_W$ is of this kind by Lemma~\ref{lem-full}.
\end{proof}

A grouped word $A^aB^bC^c$ with $a,b,c\ge1$ satisfies~\eqref{eq-words-idp}, so the family also contains arbitrarily large IDP examples with three distinct constituents.

\subsection{Exact clearing sets}\label{sec-words-clearing}
Theorem~\ref{thm-words} gives clearing set $2\Z_{>0}$ for every indecomposable hole of a word disjunction. For an arbitrary lattice point an odd multiple may also split; we determine all clearing sets, for every graph cone $C_G$ of Theorem~\ref{thm-words-graph}, with vertex monoid $S_G$ generated by the units $(\mathbf e_v\,,0)$ and the edge vectors $(-\mathbf e_u-\mathbf e_v\,,1)$; the units and one edge generate the lattice $\Z^{V(G)}\times\Z$. Let $\nu(G)$ be the largest number of pairwise disjoint odd cycles of $G$ with no edge joining two of them; induced odd cycles give the same number.

\begin{theorem}[Exact clearing sets and seminormality]\label{thm-words-clearing}
Let $G$ have an edge and $h=(y,k)$ be a lattice point of $C_G\,$. Put $b_v=\max(0,-y_v)$, $B=\sum_vb_v\,$, $b(A)=\sum_{v\in A}b_v\,$, and $H=G[\{v:b_v>0\}]$. For $X\subseteq V(H)$ let $I_X$ be the set of isolated vertices of $H-X$, let $o_X$ be the number of its components $A$ with at least two vertices and $b(A)$ odd, and put $M_X=2k-B-b(I_X)+b(X)$. Then every $M_X\ge0$, every even multiplier clears $h$, and an odd $m$ clears $h$ if and only if $mM_X\ge o_X$ for all $X$. Hence
\[
 \mathcal C(h)=2\Z_{>0}\cup\{m\text{ odd}:m\ge\tau(h)\}
\]
for some $\tau(h)\in\{1,3,5,\ldots\}\cup\{\infty\}$. The largest finite $\tau(h)$ is the largest odd integer at most $\max\{1,\nu(G)\}$, and $S_G$ is seminormal if and only if $\nu(G)\le2$.
\end{theorem}
\begin{proof}
\emph{Covers and matchings.} The point $mh$ lies in $S_G$ exactly when some multiset of $mk$ edges has degree at least $mb_v$ at every $v$: the units then supply the remaining coordinates, and a smaller cover can be padded with any edge. Replace each $v\in V(H)$ by $mb_v$ copies and each edge of $H$ by all edges between the copy sets. A least cover has $mB-\mu_m$ edges, where $\mu_m$ is the maximum matching size of this copy graph: matched pairs use one edge for two demand units, and every other unit uses one edge at its vertex, which exists because $h\in C_G$; conversely, assigning the demand units of a cover to distinct edge ends gives a matching on the edges that receive two units. By the Tutte--Berge formula~\cite{Schrijver}, $mB-2\mu_m$ is the maximum over vertex sets $U$ of the copy graph of the number of odd components of the copy graph minus $U$, less $|U|$. Such a maximizer may be taken to consist of whole copy groups: if a group is partly deleted, restoring one deleted copy lowers $|U|$ by one, and, because a surviving copy of the same vertex has the same neighbors, the restored copy joins at most one existing component, so the number of odd components drops by at most one. For $U$ the copies of $X$, the isolated vertices of $H-X$ give $mb(I_X)$ singleton components, and each component $A$ of $H-X$ with at least two vertices gives one component with $mb(A)$ vertices. Thus
\[
 mB-2\mu_m=\max_X\bigl\{m\bigl(b(I_X)-b(X)\bigr)+\mathbf 1_{m\text{ odd}}\,o_X\bigr\},
\]
and $mB-\mu_m\le mk$ is equivalent to $mM_X\ge\mathbf 1_{m\text{ odd}}\,o_X$ for all $X$. By Lemma~\ref{lem-words-cover}, $2h\in S_G\,$, so every $M_X\ge0$, every even multiplier clears $h$, and the odd clearing multipliers form an up-set.

\emph{The largest threshold.} Suppose some odd multiplier clears $h$, fix $X$, and let $p$ and $r$ count the odd-weight components of $H-X$ with at least two vertices that are bipartite and nonbipartite, respectively. Because $M_X=2k-2b(I_X)-\sum_Ab(A)$ over all components $A$ of $H-X$ with at least two vertices, $M_X\equiv p+r\pmod2$. Adding to $X$ the side of smaller demand of each odd-weight bipartite component isolates the other side and lowers $M$ by at least one per component, while the $r$ nonbipartite components survive; call the enlarged set $X'$. Then $M_X\ge p+M_{X'}$ and $M_{X'}\ge0$, so $M_X\ge p$. If $r>0$, the odd-clearing criterion applied to $X'$ gives $mM_{X'}\ge o_{X'}\ge r>0$ for the odd clearing multiplier $m$, so $M_{X'}\ge1$ by integrality and $M_X\ge p+1$. The $r$ nonbipartite components contain pairwise disjoint odd cycles with no joining edge, so $r\le\nu(G)$. Let $a$ be the largest odd integer at most $\max\{1,\nu(G)\}$. If $r=0$, then $aM_X\ge p$. If $1\le r\le a$, then $aM_X\ge r(p+1)\ge p+r$. Otherwise $r=\nu(G)$ is even and $a=r-1$: for $r=2$ parity gives $M_X\ge p+2$; for $r\ge4$ and $p=0$ parity gives $M_X\ge2$ and $aM_X\ge2(r-1)\ge r$; for $r\ge4$ and $p\ge1$, $(r-1)(p+1)\ge p+r$. So $a$ clears every point that has an odd clearing multiplier.

\emph{Attainment and seminormality.} If $a=1$, a vertex vector has $\tau=1$. If $a\ge3$, take $a$ pairwise disjoint induced odd cycles with no joining edge, unit demand on their vertices and zero elsewhere, $y=-b$, and $k=(B+1)/2$. Half weights on the cycle edges together with $\frac12(0,1)\in C_G$ show $h\in C_G\,$. For odd $m$, a cycle of length $L$ needs at least $(mL+1)/2$ edges, attained by $(m-1)/2$ copies of each cycle edge and a minimum cover of the remaining unit demand, and no edge serves two cycles; so the least cover has $(mB+a)/2$ edges, which is at most $mk$ exactly when $m\ge a$. Hence $\tau(h)=a$. The monoid $S_G$ is seminormal if $2h,3h\in S_G$ imply $h\in S_G$ for every lattice point $h$. If $\nu(G)\le2$, then $a=1$ and every point with an odd clearing multiplier splits. If $\nu(G)\ge3$, the construction with three cycles has $2h,3h\in S_G$ and $h\notin S_G\,$.
\end{proof}

\begin{corollary}[Odd clearing in word disjunctions]\label{cor-words-clearing}
For every odd $q\ge3$, the least spatial dimension of a word disjunction with a lattice point $h$ of $\tau(h)=q$ is $5q-4$, and the least degree of such a point is $3(q+1)/2$, attained in that dimension. In particular the vertex monoid of every word disjunction with $d\le10$ is seminormal, while in dimension $d=11$ some word disjunction has a vertex monoid that is not seminormal, with a witness of degree six.
\end{corollary}
\begin{proof}
By the $2$-tree construction of Subsection~\ref{sec-words-graph}, $G_W$ is chordal, so its induced odd cycles are triangles, and its triangles are the state triangles $T_0\,,\ldots,T_\ell\,$, which occur in word order. Two of them are disjoint with no joining edge exactly when they bound a bad substring (Lemma~\ref{lem-words-exc}); badness persists when a substring is enlarged and needs at least five letters. Hence $\nu(G_W)\le1+\lfloor\ell/5\rfloor=1+\lfloor(d-1)/5\rfloor$, with equality for the prefixes of $ABCABC\cdots$, whose windows of length five are bad. Theorem~\ref{thm-words-clearing} gives $d\ge5q-4$, and its construction with $q$ triangles gives equality, with $B=3q$, $k=(3q+1)/2$ and degree $-B+3k=3(q+1)/2$ by Theorem~\ref{thm-words-graph}.

For the degree bound, let $\tau(h)=q\ge3$. Some $X$ violates the inequality at $q-2$, and $M_X\ge1$, so by parity $o_X\ge(q-2)M_X+2\ge q$. Every odd-weight component with at least two vertices has demand at least three, so $B\ge3q$; and $2k\ge B$, because every edge has two ends. Hence $\deg h=\sum_vy_v+3k\ge-B+3k\ge B/2$, which gives $\deg h\ge3(q+1)/2$ when $B=3q$ (then $k\ge(3q+1)/2$) or $B\ge3q+2$, because $\deg h$ is an integer and $q$ is odd. If $B=3q+1$ and $2k>B$, then $2k\ge B+2$, as $B$ is even, so $\deg h\ge-B+3k\ge B/2+3$. In the remaining case $B=3q+1$ and $2k=B$, the bound $B\ge3o_X$ forces $M_X=1$ and $o_X=q$. If all $q$ components were bipartite, then $q=p\le M_X=1$; so $r>0$, hence $M_X\ge p+1$ and $p=0$. A nonbipartite component of demand three consists of a triangle of unit demand, because $H$ contains only vertices of positive demand. The remaining unit of demand gives $b(X)-b(I_X)=M_X=1$, so $X=\{v\}$ with $b_v=1$ and $I_X=\varnothing$. Also $M_\varnothing=-b(I_\varnothing)\ge0$ gives $M_\varnothing=0$, and because some odd multiplier clears $h$, $o_\varnothing=0$; a triangle not adjacent to $v$ would be an odd-weight component of $H$, so $v$ is adjacent to all $q$ triangles. In a triangle strip each vertex is the latest vertex of its letter during an interval of states, two vertices are adjacent exactly when these intervals overlap, and for each state triangle the union of the intervals of its vertices is an interval containing its state; for the $q$ triangles these unions are pairwise disjoint. The interval of $v$ meets the unions of the first and the last triangle, so it contains the state of a middle triangle, at which only that triangle's vertices are latest; so $v$ would be a vertex of that triangle, although $v\in X$. The case $q=3$ gives $d=11$ and degree six.
\end{proof}

\subsection{Grouped letters}\label{sec-grouped}
Finally, we restate the forbidden-subsequence criterion in terms of grouped letters, which gives the linear-time test for IDP mentioned in the introduction.
\begin{proposition}[Grouped letters]\label{prop-grouped}
A word $W$ avoids the subsequences $XYZXY$ and $XYZYX$, for distinct letters $X,Y,Z$, if and only if, for some two letters, every occurrence of one precedes every occurrence of the other.
\end{proposition}
\begin{proof}
Use the graph $G_W$ of Theorem~\ref{thm-words-graph}: its edges between letters $X$ and $Y$ are $X_iY_j$ for the states $(i,j)$ visited while reading $W|_{XY}\,$, and every odd cycle uses an edge of each of the three letter pairs, because two pair types alone form a bipartite graph. By Theorem~\ref{thm-words} and Lemma~\ref{lem-words-exc}, $W$ contains one of the subsequences exactly when two state triangles $T_i\,,T_j$ are disjoint and joined by no edge. If all $X$'s precede all $Y$'s, the states of $W|_{XY}$ are $(i,0)$ and then $(m_X\,,j)$, so every $XY$ edge contains $Y_0$ or $X_{m_X}\,$, which are adjacent; two disjoint triangles use different ones of these two vertices and are therefore joined by an edge. If no two letters are grouped, every letter occurs, so the initial triangle $T_0$ and the final triangle $T_\ell$ are disjoint, and an edge between them would be $X_0Y_{m_Y}$ or $X_{m_X}Y_0\,$, i.e.\ the state $(0,m_Y)$ or $(m_X\,,0)$ of $W|_{XY}\,$, which would group $X$ and $Y$. So $(0,\ell)$ is bad.
\end{proof}

\subsection*{Acknowledgments}
This work was supported in part by the Office of Naval Research under grant N00014-24-1-2694. During the preparation of this work, the authors used Claude and ChatGPT to assist in developing and checking the work. The authors take full responsibility for the content of the article.

\bibliographystyle{amsplaindoi}
\bibliography{arithmetic_arxiv}

\end{document}